\documentclass[11pt]{article}
\pdfoutput=1 

\usepackage{amssymb}
\usepackage[margin=1in]{geometry}

\usepackage{lineno}
\usepackage{graphicx}
\usepackage{moreverb,url}
\usepackage[ruled,vlined]{algorithm2e}
\usepackage[colorlinks,bookmarksopen,bookmarksnumbered,citecolor=red,urlcolor=red]{hyperref}

\usepackage{xcolor,tabularx,nccmath}
\usepackage{amsmath,amsthm,enumitem}
\usepackage{booktabs,caption,subcaption,mathtools,multirow}

\newcommand{\bigM}{\mathbb{M}}
\DeclareMathAlphabet{\mathdutchcal}{U}{dutchcal}{m}{n}
\SetMathAlphabet{\mathdutchcal}{bold}{U}{dutchcal}{b}{n}
\DeclareMathAlphabet{\mathdutchbcal}{U}{dutchcal}{b}{n}

\newcommand{\sB}{\mathdutchcal{B}}

\newcommand{\sT}{\mathdutchcal{T}}

\newcommand{\sI}{\mathdutchcal{I}}
\newcommand{\sL}{\mathdutchcal{L}}

\newcommand{\sR}{\mathdutchcal{R}}
\newcommand{\sC}{\mathdutchcal{C}}

\usepackage[nameinlink,capitalise]{cleveref}
\hypersetup{colorlinks={true},allcolors={blue}}

\crefformat{section}{\S#2#1#3}
\crefformat{subsection}{\S#2#1#3}
\crefformat{equation}{#2Equation~\textup{(#1)}#3}
\crefformat{cons}{#2Constraint~\textup{(#1)}#3}
\labelcrefformat{cons}{\textup{#2(#1)#3}}
\crefformat{consset}{#2Constraints~\textup{(#1)}#3} 
\labelcrefformat{consset}{\textup{#2(#1)#3}}
\crefformat{func}{#2Function~\textup{(#1)}#3}
\labelcrefformat{func}{#2Function~\textup{(#1)}#3}
\crefformat{algo}{#2Algorithm~\textup{#1}#3}
\crefformat{objfunc}{#2Objective function~\textup{(#1)}#3}
\labelcrefformat{objfunc}{\textup{#2(#1)#3}}
\crefformat{ch}{#2Chapter~\textup{#1}#3}
\crefformat{sec}{#2Section~\textup{#1}#3}
\crefformat{tab}{#2Table~\textup{#1}#3}
\crefformat{fig}{#2Figure~\textup{#1}#3}
\crefdefaultlabelformat{#2#1#3} 
\Crefname{figure}{Figure}{Figures}
\crefformat{eq}{#2equation~\textup{(#1)}#3}
\labelcrefformat{eq}{\textup{(#2#1#3)}}
\crefformat{theo}{#2Theorem~\textup{#1}#3}
\crefformat{lem}{#2Lemma~\textup{#1}#3}
\crefformat{def}{#2Definition~\textup{#1}#3}
\crefformat{app}{#2\textup{\!#1}#3}
\crefformat{scenario}{#2Scenario~\textup{#1}#3}
 \AtBeginDocument{%
    \crefname{figure}{Figure}{figures}%
}
\let\origref\cref
\def\cref#1{\origref{#1}}

\makeatletter
\let\oldappendix\appendix
\renewcommand{\appendix}{%
  \oldappendix
  \renewcommand{\thesection}{Appendix \Alph{section}}%
}
\makeatother

\crefformat{app}{#2\textup{#1}#3} 

\expandafter\def\expandafter\UrlBreaks\expandafter{\UrlBreaks%
  \do\a\do\b\do\c\do\d\do\e\do\f\do\g\do\h\do\i\do\j%
  \do\k\do\l\do\m\do\n\do\o\do\p\do\q\do\r\do\s\do\t%
  \do\u\do\v\do\w\do\x\do\y\do\z\do\A\do\B\do\C\do\D%
  \do\E\do\F\do\G\do\H\do\I\do\J\do\K\do\L\do\M\do\N%
  \do\O\do\P\do\Q\do\R\do\S\do\T\do\U\do\V\do\W\do\X%
  \do\Y\do\Z}

\usepackage{placeins}
\usepackage{etoolbox}

\newcommand{\reducedisplayskips}{%
  \setlength{\abovedisplayskip}{1pt}
  \setlength{\belowdisplayskip}{1pt}
  \setlength{\abovedisplayshortskip}{0pt}
  \setlength{\belowdisplayshortskip}{0pt}
}

\AtBeginEnvironment{equation}{\reducedisplayskips}

\AtBeginEnvironment{align}{\reducedisplayskips}
\AtBeginEnvironment{gather}{\reducedisplayskips}
\AtBeginEnvironment{flalign}{\reducedisplayskips}
\AtBeginEnvironment{multline}{\reducedisplayskips}

\usepackage{stackengine}

\usepackage{algpseudocode}

\DeclarePairedDelimiter\floor{\lfloor}{\rfloor}

\usepackage[normalem]{ulem}

\usepackage{rotating}

\crefname{defi}{definition}{definitions}
\Crefname{defi}{Definition}{Definitions}
\crefname{lemma}{lemma}{lemmas}
\Crefname{lemma}{Lemma}{Lemmas}
\crefname{assumption}{assumption}{assumptions}
\Crefname{assumption}{Assumption}{Assumptions}

\usepackage{multicol}

\usepackage{authblk}
\providecommand{\keywords}[1]{\noindent\textbf{Keywords:} #1}
\usepackage{natbib}

\makeatletter
\AtBeginDocument{
  \labelcrefformat{equation}{\textup{(#2#1#3)}}
  
  \let\old@label@in@display\label@in@display
  \def\label@in@display{\@ifnextchar[\strip@opt@label@display\old@label@in@display}
  \def\strip@opt@label@display[#1]#2{\old@label@in@display{#2}}
}
\makeatother

\begin{document}

\title{Integrated Optimization of Scheduling and Flexible Charging in Mixed Electric-Diesel Urban Transit Bus Systems}

\author[a]{Sadjad Bazarnovi}
\author[b]{Taner Cokyasar\thanks{Corresponding author: \texttt{tcokyasar@anl.gov}}}
\author[b]{Omer Verbas}
\author[a]{Abolfazl (Kouros) Mohammadian}

\affil[a]{University of Illinois Chicago, 842 W. Taylor St., Chicago, IL, 60607 USA}
\affil[b]{Argonne National Laboratory, 9700 S. Cass Avenue, Lemont, IL, 60439 USA}

\maketitle

\begin{abstract}
The transition of transit fleets to alternative powertrains offers a potential pathway to reducing the cost of mobility. However, the limited range and long charging durations of battery electric buses (BEBs) introduce significant operational complexities, necessitating innovative scheduling and charging strategies. This study proposes an integrated mixed-integer linear programming model to optimize vehicle scheduling and charging strategies for mixed fleets of BEBs and diesel buses. Unlike existing models, which often assume a fixed BEB fleet size or restrict charging to a single charger type, our approach simultaneously determines the optimal fleet composition, scheduling, and flexible partial charging strategy incorporating both slow and fast chargers at garages and terminal stations. The model minimizes combined fleet purchase and operational costs. A queuing strategy is introduced, departing from traditional first-come, first-served methods by dynamically allocating waiting and charging times based on operational priorities and resource availability, improving overall scheduling efficiency. To overcome computational complexities arising from numerous variables, a column generation framework is developed, facilitating scalable solutions for large-scale transit networks. Numerical experiments using real-world transit data from the Chicago Transit Authority and the Pace suburban bus systems demonstrate the model's effectiveness. Results indicate that while a full transition to alternative powertrains results in a modest cost increase, optimal mixed-fleet configurations can actually reduce total system costs. Furthermore, sensitivity analyses reveal that restricting charging to garages significantly increases fleet size and operational costs, underscoring the potential of distributed opportunistic charging.
\end{abstract}


\section{Introduction}\label[sec]{sec:introduction}
Public transportation networks are the backbone of urban mobility, and as cities expand, the demand for efficient and reliable public transit continues to rise in order to alleviate congestion and improve accessibility. However, many transit agencies face challenges with aging Diesel Bus (DB) fleets that require frequent maintenance, compounded by the volatility of fuel prices, which increases operational costs. In response, technological advancements in alternative-fuel-powered mobility, such as enhanced battery efficiency and fast-charging infrastructure, have positioned battery electric buses (BEBs) as a promising alternative to DBs. BEBs offer lower operating and maintenance costs due to their simpler powertrains \citep{ALVO2021102528, NAJAFI2025104664, HE2023128227, XIE2023103551}. Recognizing these advantages, major U.S. transit agencies, including the Los Angeles County Metropolitan Transportation Authority (LA Metro), New York City's Metropolitan Transportation Authority (MTA), and Chicago Transit Authority (CTA), have committed to proliferating the powertrain choices of their fleets in the near future \citep{HE2023103653, LIU2023120483}.

Despite their advantages, BEBs present several challenges, including limited battery range, high initial capital expenditures, and operational complexities related to scheduling and charging infrastructure \citep{TANG2023103652, PERUMAL2021105268}. One critical challenge is optimizing vehicle scheduling to ensure the efficient operation of the urban transit system. The Vehicle Scheduling Problem (VSP) involves assigning a set of timetabled trips to a fleet of vehicles while minimizing costs such as fleet size and operational expenses \citep{PERUMAL2021105268}. Each bus can operate only one trip at a time, and every revenue trip, consisting of a sequence of stops, must be assigned to a single vehicle \citep{ALVO2021102528}. When incorporating BEBs into transit fleets, VSP becomes significantly more complex due to additional constraints, including limited battery capacity, charging station availability, and energy consumption patterns \citep{TANG2023103652}. These challenges require advanced optimization strategies and efficient scheduling frameworks to facilitate the successful integration of BEBs into public transportation networks \citep{XIE2023103551, HE2023103653}.

This study addresses these complexities by optimizing bus and charging schedules for a mixed fleet of BEBs and DBs. We propose a Mixed-Integer Linear Programming (MILP) model that minimizes the total fixed and variable fleet costs. The model optimally schedules BEB charging by determining the start time, duration, and location of charging events while ensuring full coverage of daily trips. Unlike many existing studies that assume a fixed BEB fleet size \citep{HE2023103653}, our model determines the optimal number of both BEBs and DBs. Additionally, while most previous studies consider only a single charger type \citep{XIE2023103551}, our model incorporates both slow and fast chargers across multiple station locations (garages and terminal stops), adopting the infrastructure topology from \citet{bazarnovi2024problem}.

Given the NP-hard nature of the VSP with resource constraints (e.g., battery limitations) and the added dimensionality of charging constraints, exact solution methods struggle to scale for large real-world networks. To overcome this, we develop a specialized Column Generation (CG) framework. This approach decomposes the problem into a master problem, which selects optimal schedules, and pricing subproblems that generate feasible vehicle routes. Specifically, we model the BEB subproblem as a Shortest Path Problem with Resource Constraints (SPPRC) to strictly enforce battery and charging rules, and we implement stabilization mechanisms to mitigate algorithmic stalling caused by degeneracy.

\begin{figure}[!ht]
    \centering
    \includegraphics[width=0.5\linewidth]{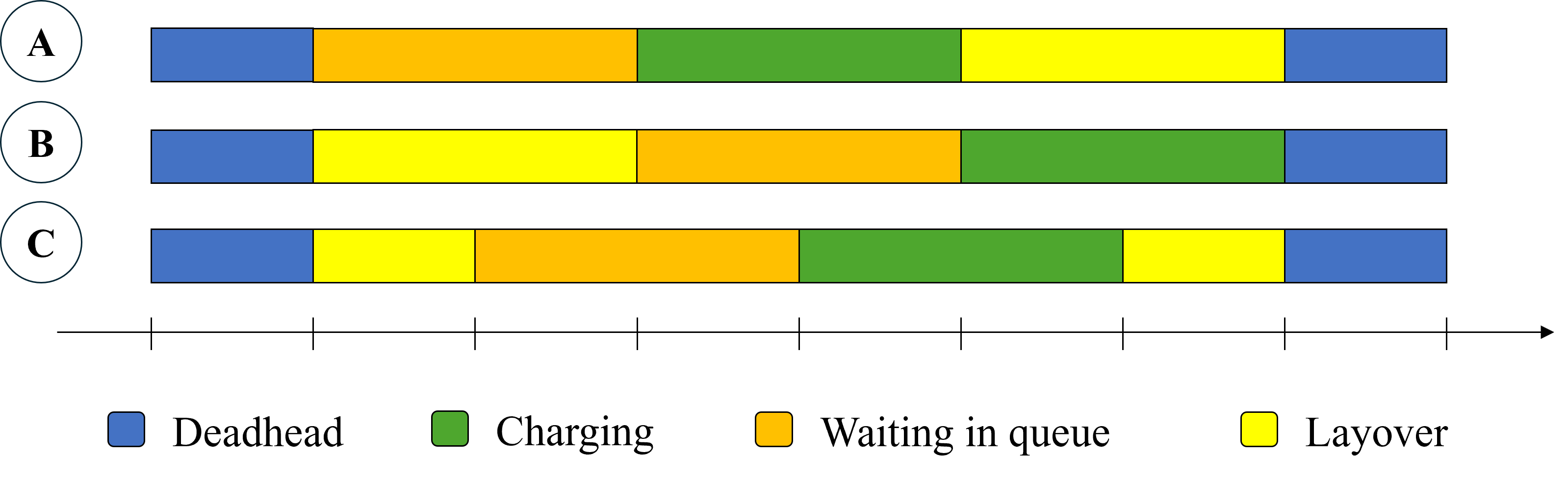}
    \caption{Queuing strategies for BEBs at charging stations.}
    \label[fig]{charging_strategies}
\end{figure}

The specific contributions of this paper are as follows: (i) To the best of our knowledge, this is the first study to simultaneously optimize vehicle scheduling, fleet composition (DB vs. BEB), and charging decisions for a mixed fleet, considering both garage and on-route charging with multiple charger types. (ii) We introduce a \textit{flexible charging} strategy that allows for partial recharging and eliminates rigid First-Come, First-Served (FCFS) constraints. As illustrated in \cref{charging_strategies}, our strategy (Approach C) dynamically allocates \textit{layover} time, defined as the idle time at terminals between trips. This allows buses with tighter schedules to prioritize charging slots, reducing congestion compared to traditional approaches. (iii) We propose a tailored CG algorithm capable of solving large-scale instances that are intractable for standard solvers. The framework includes handling of the BEB-SPPRC and a stalling criterion to ensure convergence in degenerate solution spaces.
    

The remainder of this paper is structured as follows: \cref{literature} reviews the relevant literature. \cref{problem} presents the problem formulation. \cref{solution} details the solution approaches, including the CG heuristic. \cref{experiments} describes the case study and computational analysis, followed by \cref{conclusion} which summarizes key findings.

\section{Literature review}\label[sec]{literature}
The integration of BEBs into public transportation networks has been the subject of extensive research, particularly concerning scheduling and charging strategies to optimize their operational efficiency. Various studies have explored different approaches to address the challenges associated with BEB deployment, focusing on minimizing costs, optimizing scheduling, and managing charging infrastructure. \citet{perumal2022electric} provided a thorough review on BEB scheduling problems. \cref{literature_comparison} presents a summary of related studies. To the best of the authors' knowledge, this is the first study that simultaneously optimizes vehicle and charging schedules for a mixed fleet of BEBs and DBs, determines the optimal fleet composition, and incorporates both garages and terminal stations as charging locations. The model also accounts for multiple recharging options, station capacity constraints, various charger types (fast and slow), and the feasibility of partial charging. Additionally, the proposed CG heuristic enhances scalability, enabling the efficient solution of large-scale case studies.

\begin{table*}[!ht]\centering
\caption{Summary of related studies.}\label[tab]{literature_comparison}
\scriptsize
\begin{tabularx}{\textwidth}{>{\raggedright\arraybackslash}X ccccccccccc X X}
\toprule
\textbf{} & \multicolumn{4}{c}{\textbf{Problem}} & \multicolumn{6}{c}{\textbf{Charging}} & \multicolumn{3}{c}{\textbf{Approach}}\\
\cmidrule(lr){2-5} \cmidrule(lr){6-11} \cmidrule(lr){12-14}
\textbf{Study} & \textbf{VS} & \textbf{CS} & \textbf{NB} & \textbf{BF} & \textbf{ST} & \textbf{CL} & \textbf{MO} & \textbf{SC} & \textbf{PC} & \textbf{CT} & \textbf{MH} & \textbf{Model} & \textbf{CY} \\
\midrule
\citet{HU20219564617} &\checkmark &\checkmark &\checkmark &BEB &S & &\checkmark &\checkmark &\checkmark & &GA &MILP &Sydney, Australia \\
\citet{PERUMAL2021105268} & \checkmark &  &  & BEB & G &  &  &  &  & & ALNS & MINLP & Denmark \& Sweden \\
{\citet{Jovanovic202114206610}} &\checkmark &\checkmark &\checkmark &BEB &G & & & & & &GRASP &MIP &Florida, USA \\
\citet{ALVO2021102528} &\checkmark &\checkmark &\checkmark &Mixed &G & &\checkmark &\checkmark &\checkmark && - &two-stage MIP &Santiago, Chile \\
\citet{HU2022103732} & &\checkmark & &BEB &S &\checkmark &\checkmark & &\checkmark &&- &MILP &Sydney, Australia \\
\citet{HE2023103587} &\checkmark &\checkmark &\checkmark &BEB &G, S &\checkmark &\checkmark & &\checkmark &&- &MINLP &Salt Lake City, UT \\
\citet{XIE2023103551} &\checkmark &\checkmark &\checkmark &BEB &G & & & &\checkmark & \checkmark&- &MINLP &Beijing, China \\
\citet{LIU2023120483} & & &\checkmark &Mixed &S &\checkmark &\checkmark &\checkmark & & \checkmark&GA &Simulation &Dallas, TX \\
\citet{BAO2023120512} & &\checkmark & &BEB &G, S & &\checkmark &\checkmark &\checkmark &\checkmark &- &MLIP &Shanghai, China \\
\citet{HE2023103653} & &\checkmark & &BEB &G & & & & & \checkmark& GA &MILP &Small network \\
\citet{TANG2023103652} & &\checkmark &\checkmark &BEB &G & & & & &&GA &MINLP &Dandong, China \\
\citet{GAIROLA2023103697} & &\checkmark & &BEB &G &\checkmark & &\checkmark &\checkmark &&- &MILP &New Delhi, India \\
\citet{HE2023128227} & &\checkmark & &Mixed &G, S &\checkmark &\checkmark &\checkmark & &\checkmark&- &bi-obj. MILP &Salt Lake City, UT \\
\citet{NAJAFI2025104664} &\checkmark &\checkmark &\checkmark &BEB &G, S &\checkmark &\checkmark &\checkmark &\checkmark &\checkmark&- &MINLP &Skövde, Sweden \\ 
\hline
\textbf{This study} &\checkmark &\checkmark &\checkmark &Mixed &G, S & &\checkmark &\checkmark &\checkmark &\checkmark & CG &MILP &Chicago, IL \\
\bottomrule
\multicolumn{14}{p{450pt}}{\emph{VS: Vehicle scheduling, CS: Charging scheduling, NB: Number of buses, BF: Bus fleet, ST: Charging station, CL: Charging location, MO: Multiple options, SC: Station capacity, CT: Charger type, PC: Partial charging, MH: Meta-heuristic, CY: Case study, G: Garages. S: Stops}}\\
\end{tabularx}
\end{table*}

Several studies focused on optimizing bus scheduling and charging strategies. 
\citet{HE2023103587} proposed a simplified MILP model to optimize bus charging while ensuring feasibility and cost-effectiveness. 
Similarly, \citet{HE2023103653} and \citet{GAIROLA2023103697} proposed mixed-integer nonlinear programming (MINLP) models to jointly optimize vehicle scheduling and charging management, with a focus on minimizing total operational costs. 
\citet{HU2022103732} explored fast charging at bus stops, minimizing passenger delay costs, and accounting for time-varying electricity prices.
\citet{BAO2023120512} coordinated charging schedules across multiple routes and terminals using a Lagrangian relaxation framework.
\citet{TANG2023103652} integrated a skip-stop operation strategy into single-line electric bus scheduling.
\citet{XIE2023103551} extended the scheduling problem by incorporating driver schedules and considering different charging modes.

Furthermore, existing optimization frameworks can be categorized based on the depot configuration and trip assignment logic. Several studies addressed the Multi-Depot Vehicle Scheduling Problem (MDVSP), where schedules are generated jointly for the entire network without pre-partitioning trips to specific garages \citep{desaulniers1998multi, dell1993heuristic, gkiotsalitis2023exact}. In contrast, our study adopts a single-depot VSP (SDVSP) framework, aligning with the modeling approaches found in \citet{cokyasar2023electric, cokyasar2023solving, freling2001models, rinaldi2020mixed}.

The operational challenge of limited driving range and long charging times has spurred research into fast-charging and en-route charging approaches.
\citet{HU20219564617, HU2022103732} explored strategies ranging from depot plug‐in charging to fast charging at bus stops.
Similarly, \citet{GAIROLA2023103697} investigated fast-charging configurations by balancing the reduced charging duration against increased electricity costs.
Furthermore, \citet{bazarnovi2024problem} focused on the placement and allocation of charging equipment, showing that advancements in charging technology are more impactful when they result in higher charger power and shorter charging times compared to lower charger costs.

The integration of BEBs with conventional DB fleets presents additional layers of complexity, particularly in terms of scheduling and charging infrastructure deployment. 
While most studies focused on fully electric fleets \citep{Jovanovic202114206610, LIU2023120483, NAJAFI2025104664}, \citet{ALVO2021102528} explored operations within a mixed fleet of BEBs and DBs. 
This contrasts with the work of \citet{BAO2023120512}, which exclusively examined a BEB fleet, focusing on the impact of time-of-use tariffs, station load capacities, and battery capacities on system performance. 
\citet{NAJAFI2025104664} further extended this discussion by integrating renewable energy resources with charging scheduling, highlighting how coordinated planning not only minimizes operational costs but also addresses technical constraints of power distribution networks. 
The charging strategies also vary, with some studies emphasizing en-route fast charging \citep{HE2023103587, BAO2023120512, NAJAFI2025104664}, while others focusing on garage charging \citep{HE2023103653, TANG2023103652, GAIROLA2023103697}. 
These studies collectively suggested that a holistic view, encompassing both mixed fleet dynamics and infrastructure considerations, is crucial for achieving resilient transit operations.

The complexity of BEB scheduling and charging planning has led to the adoption of various optimization techniques. 
\citet{ALVO2021102528} leveraged a Benders’ decomposition method to dynamically inject feasibility cuts into the optimization process, improving solution efficiency. 
\citet{Jovanovic202114206610} formulated the problem as a mixed-integer program (MIP) and employ a Greedy Randomized Adaptive Search Procedure (GRASP) to efficiently approximate solutions for large-scale examples. 
Moreover, the adaptive large neighborhood search (ALNS) approach demonstrated by \citet{PERUMAL2021105268} offers additional flexibility in exploring various scheduling configurations while accommodating the operational constraints of BEB systems. 
Complementing this, \citet{HU20219564617} and \citet{HE2023103653} adopted GA-based heuristics to navigate the expansive solution space inherent in multi-dimensional scheduling problems. 

A range of real-world case studies underscores the practical relevance and adaptability of optimization models in managing BEB operations. 
For example, \citet{PERUMAL2021105268} validated their integrated scheduling framework using data from Danish and Swedish transit systems.
Similarly, \citet{ALVO2021102528} demonstrated the efficacy of their decomposition approach in a mixed fleet setting through a case study from Santiago, Chile, while \citet{LIU2023120483} applied a simulation-based framework to an airport shuttle system at Dallas-Fort Worth International Airport. 
Additionally, \citet{HU2022103732, HU20219564617} tested their methods on bus routes in Salt Lake City, UT, demonstrating the practicality of en-route charging. 
\citet{HE2023103653, HE2023128227} applied their framework to a sub-transit network in Salt Lake City, comparing multiple operational scenarios. 
\citet{BAO2023120512} analyzed a regional BEB fleet in Shanghai, providing insights into the impacts of charging infrastructure configurations and energy pricing policies. 

\section{Problem definition}\label[sec]{problem}
In this section, we formulate an SDVSP and the charging scheduling problem as an MILP model for a mixed fleet of BEBs and DBs. To ensure clarity and consistency, we use specific notation throughout the formulation. Calligraphic letters, such as \( \sI \), represent sets. Lowercase Roman letters, like \( x_{ijk} \), are used for variables and indices, while Greek letters, such as \( \delta \), appear as superscripts modifying parameters and variables. Uppercase Roman letters, including \( T_{ij} \), denote parameters. Blackboard bold letters, such as \( \mathbb{R} \), are reserved for domains. Notations specific to small and large numbers, such as \( \epsilon \) for an infinitesimal value and \( \bigM \) for a sufficiently large number, are excluded from this convention. For a complete list of all sets, parameters, and variables used in this mathematical formulation, please refer to \cref{set}, \cref{param}, and \cref{vars}.

\begin{table}[!ht]
  \scriptsize
  \caption{Sets used in the MILP.}\label[tab]{set}
  \begin{tabularx}{\linewidth}{lX}
  \toprule
    \textbf{Set} & \textbf{Definition}\\
    \midrule
    $\sC$ & Set of charging stations including the garage $s$\\
    $\sC_i$ & Subset of charging stations in $\sC$ that a BEB can visit after completing trip $i \in \sI$\\
    $\sC_{ij}$ & Subset of charging stations in $\sC_i$ that a BEB can visit after completing trip $i \in \sI$ and before starting trip $j \in \sI_{ic}^\sigma$\\
    $\sI$ & Set of timetabled trips\\
    $\sI^\tau_i$ & Set of all feasible trips connecting to trip $i \in \sI$, i.e., $\{j \in \sI | T_{ji}^\delta \leq T_i^\alpha - T_j^\beta \leq G\} \cup \{s\}$, where $s$ denotes the garage\\
    $\sI^\omega_i$ & Set of all feasible trips connecting from trip $i \in \sI$, i.e., $\sI^\omega_i = \sI_i \cup \{s\}$, where $s$ denote the garage that buses return to\\
    $\sI_i$ & Set of trips in $\sI$ that are compatible with trip $i$, i.e., $\{j \in \sI | T_{ij}^\delta \leq T_j^\alpha - T_i^\beta \leq G \land i \neq j\}$\\
    $\sI_{i}^\lambda$ & Subset of trips in $\sI_i$ where layover between trips does not exceed $L$, i.e., $\{j \in \sI_i | T_j^{\alpha} - T_i^{\beta} - T_{ij}^\delta \leq L\}$ \\
    $\sI_{ic}^\sigma$ & Subset of trips in $\sI_i$ where a visit to $c \in \sC_i$ between trips is possible, i.e., $\{j \in \sI_i | T_j^{\alpha} - T_i^{\beta} - T_{ic}^\delta - T_{cj}^\delta \geq W\}$\\
    $\sI_{ict}^\sigma$ & Subset of trips in $\sI_{ic}^\sigma$ where a visit to $c \in \sC_i$ between trips is possible at time step $t \in \sT_{ijc}$\\
    {$\sR$} & Set of all feasible arcs connecting two compatible trips, $\sR = \bigcup_{i \in \sI} \left(i \times \sI_i\right) \cup (s \times \sI) \cup (\sI \times s)$, where $s$ denotes garage\\
    $\sT$ & Set of time steps\\
    $\sT_{ijc}$ & Subset of time steps at the beginning of which a BEB can visit $c \in \sC_i$, i.e., $\{t \in \sT | T_i^\beta + T_{ic}^\delta \leq t T^\Delta \leq T_j^\alpha - T_{cj}^\delta\}$\\
    \bottomrule
  \end{tabularx}
\end{table}

\begin{table}[!ht]
  \scriptsize
  \caption{Parameters used in the MILP.}\label[tab]{param}
  \begin{tabularx}{\linewidth}{lX}
  \toprule
    \textbf{Parameter} & \textbf{Definition}\\
    \midrule
    $A^\nu$ & Minimum required proportion of BEBs in the total fleet, $A^\nu\in [0,1]$\\
    $A^\tau$ & Minimum required share of total revenue trip time served by BEBs, $A^\tau\in [0,1]$\\
    $B_i$ & Energy consumption measured in time units of trip $i \in \sI$\\
    $B^\iota$ & Battery level measured in time units at the beginning of a run\\
    $\overline{B}$ & Maximum battery level measured in time units\\
    $\underline{B}$ & Minimum battery level measured in time units\\
    $D_t$ & Actual time associated with beginning of the time step $t \in \sT$, that is $D_t = T^\Delta \times t$\\
    $F^\varepsilon$ & The cost incurred per time unit for operating a BEB\\
    $F^\kappa$ & The cost incurred per time unit for operating a DB\\
    $G$ & Maximum acceptable gap measured in time units between two compatible trips\\
    $L$ & Maximum acceptable layover time measured in time units between two connected trips \\
    $\bigM$ & A sufficiently large number\\
    $P$ & The penalty cost for each unit of shortfall in the required minimum BEB share in the fleet\\
    $R_c$ & Rate of charging at the charging station $c \in \sC$\\
    $T_i^{\alpha}$ & Start time of trip $i \in \sI$\\
    $T_i^{\beta}$ & End time of trip $i \in \sI$\\
    $T_{ij}^\delta$ & Deadheading time from $i$ to $j$\\
    $T^\Delta$ & Duration of time steps\\
    $\overline{T}$ & End of planning horizon in time units\\
    $V^\varepsilon$ & The cost of a BEB per time unit, amortized over the vehicle's lifespan\\
    $V^\kappa$ & The cost of a DB per time unit, amortized over the vehicle's lifespan\\
    $W$ & The minimum amount of time a vehicle must spend visiting the garage/charging station between its trips\\
    $Z_c$ & Number of charger plugs at the charging station $c \in \sC$\\
    $\epsilon$ & An infinitesimal number\\
    \bottomrule
  \end{tabularx}
\end{table}

A set of timetabled trips $\sI$ is provided with known origin and destination locations and start and end times. For trip $i \in \sI$, let $O_i$ and $D_i$ be its origin and destination locations, $T_i^{\alpha}$ and $T_i^{\beta}$ be its start and end times, and let $B_i$ be its duration. Each trip originates from a single garage, indexed by \( s \), located at \( L_s \). The deadheading time from trips $i$ to $j$ is calculated as the distance between $D_i$ and $O_j$ divided by the average speed of the vehicle ($\bar{V}$), i.e. $T_{ij}^\delta = \text{dist}(D_i, O_j)/ \bar{V}$. The trip-to-garage and garage-to-trip deadheading times are calculated as $T_{si}^\delta = \text{dist}(L_s, O_i)/ \bar{V}$ and $T_{is}^\delta = \text{dist}(D_i, L_s)/ \bar{V}$, respectively. Also, there is a set of charging stations, $\sC$, each with a known location and limited capacity. For each $c \in \sC$, let $L_c$, $R_c$ and $Z_c$ be its location, charging power per plug, and the number of plugs, respectively. We consider multiple charger types with varying charging speeds; however, all plugs at a given charging station are identical. For locations with multiple charger types, we model them as separate stations at the same exact location, each equipped with only one type of charger plug. We assume that garages house charging equipment, i.e., $s \in \sC$. 

BEBs begin their daily schedule with a State-of-Charge (SoC) level of \( B^\iota \). At the end of each run, BEBs must have sufficient time to recharge and restore their SoC to \( B^\iota \) in preparation for the following day. To mitigate range anxiety and extend battery life, BEBs should maintain an SoC above a minimum threshold, \( \underline{B} \), at all times during operation. Additionally, their SoC should not exceed an upper limit, \( \overline{B} \), to preserve battery health. We also assume that the charging and discharging behavior of batteries is linear within this interval \citep{bazarnovi2024problem}. On the other hand, DBs are assumed to have sufficient fuel to complete any feasible run without constraints. Finally, we assume that the transit agency aims to ensure that at least \( A^\nu \)\% of its fleet consists of BEBs and that at least \( A^\tau \)\% of the total revenue service time is operated by BEBs. Note that many of these assumptions do not obstruct hard boundaries on the definition of the model, and some can be easily relaxed.

We define a trip pair $(i,j)$ as compatible if the following conditions are satisfied: (1) trip $j \in \sI$ starts after trip $i \in \sI$; (2) there is sufficient time to deadhead from $D_i$ to $O_j$ (i.e., $T_j^\alpha \geq T_i^\beta + T_{ij}^\delta$); and (3) the time gap between $T_i^\beta$ and $T_j^\alpha$ does not exceed the maximum threshold $G$ (i.e., $T_j^{\alpha} - T_i^{\beta} \leq G$). Based on these criteria, for each trip $i \in \sI$, we define $\sI_i$ as the subset of trips compatible with $i$. Let $\sC_i$ denote the subset of charging stations $c \in \sC$ that a BEB can access after completing trip $i$. We define the subset $\sI_{ic}^\sigma$ as all trips $j \in \sI_i$ where the interval between trips $i$ and $j$ allows the bus to visit station $c \in \sC_i$ and remain there for at least a minimum duration $W$ (i.e., $T_i^\beta + T_{ic}^\delta + T_{cj}^\delta + W \leq T_j^\alpha$). Additionally, we impose a maximum acceptable idle time $L$ ($L < G$) at trip terminals; if the layover between two trips exceeds $L$, the bus must relocate to a garage or charging station. Consequently, the subset $\sI_i^\lambda$ contains trips $j \in \sI_i$ where the layover time does not exceed $L$ (i.e., $T_j^\alpha - T_i^\beta - T_{ij}^\delta \leq L$).

The binary variable $y_{ij}^\varepsilon$ equals 1 if the compatible trip pair $(i, j)$ is served by a BEB, and 0 otherwise. Similarly, $y_{ij}^\kappa$ is defined for trips served by a DB. To indicate whether a BEB visits a charging station $c \in \sC_i$ between trips $i$ and $j$, we define the binary variable $q_{ijc}^\varepsilon$. Such a visit may occur either to recharge or to satisfy the layover constraint for intervals exceeding $L$. For DBs, the binary variable $q_{ij}^\kappa$ equals 1 if the bus visits the garage during a layover exceeding $L$, and 0 otherwise. The continuous variable $b_i$ represents the SoC at the start of trip $i \in \sI$. For all trips, we ensure the BEB has sufficient SoC to complete the trip and reach either the next trip's origin or a charging station $c \in \sC_i$. Note that while BEBs may recharge fully or partially, the SoC must always remain within the operational bounds $[\underline{B}, \overline{B}]$.

\begin{table}[!ht]
  \scriptsize
  \caption{Variables used in the MILP.}\label[tab]{vars}
  \begin{tabularx}{\linewidth}{lX}
  \toprule
    \textbf{Variable} & \textbf{Definition}\\
    \midrule
    $b_i$ & SoC at the beginning of trip $i \in \sI$\\
    $q_{ijc}^\varepsilon$ & $\begin{cases}
        1 & \text{\parbox[t][][t]{0.80\textwidth}{if the BEB visits charging station $c \in \sC_{ij}$ between the compatible trips $i \in \sI$ and $j \in \sI_{ic}^\sigma$}}\\
        0 & \text{otherwise}\\
        \end{cases}$\\
    $q_{ij}^\kappa$ & $\begin{cases}
          1 & \text{\parbox[t][][t]{0.80\textwidth}{if the DB visits the garage $s$, between the compatible trips $i \in \sI$ and $j \in \sI_{is}^\sigma$}}\\
          0 & \text{otherwise}\\
          \end{cases}$\\
    $s_{ijc}$ & The start time of recharging at charging station $c \in \sC_{ij}$ between the trip $i \in \sI$ and its compatible trip $j \in \sI_{ic}^\sigma$\\
    $t_i^\eta$ & Elapsed time from the start of the run to the beginning of trip $i \in \sI$\\
    $u_{ijc}$ & The charging time at charging station $c \in \sC_{ij}$ between the compatible trips $(i,j)$\\
    $v$ & The violation of the required minimum BEB share (i.e., $A^\nu \sum_{i \in \sI}(y^\kappa_{is}+y^\varepsilon_{is})$), representing the shortfall in the number of BEBs in the fleet\\
    $v^\prime$ & The violation of the required minimum BEB-conducted revenue trip time, representing the shortfall in the amount of revenue trip time conducted by BEBs.\\
    $x_{ijct}$ & $\begin{cases}
        1 & \text{\parbox[t][][t]{0.80\textwidth}{if the BEB servicing the compatible trips $(i,j): i \in \sI, j \in \sI_{ic}^\sigma$ is recharging at $c \in \sC_{ij}$ where $t \in \sT_{ijc} \cap \floor{s_{ijc}/T^\Delta} \leq t \leq \floor{(s_{ijc} + u_{ijc})/T^\Delta}$}}\\
        0 & \text{otherwise}\\
        \end{cases}$\\
    $x_{ijct}^\alpha$ & $\begin{cases}
        1 & \text{\parbox[t][][t]{0.80\textwidth}{if the BEB servicing the compatible trips $(i,j): i \in \sI, j \in \sI_{ic}^\sigma$ is waiting to be charged at $c \in \sC_{ij}$ where $t \in \sT_{ijc} \cap t < \floor{s_{ijc}/T^\Delta}$}}\\
        0 & \text{otherwise}\\
        \end{cases}$\\
    $x_{ijct}^\beta$ & $\begin{cases}
        1 & \text{\parbox[t][][t]{0.80\textwidth}{if the BEB servicing the compatible trips $(i,j): i \in \sI, j \in \sI_{ic}^\sigma$ is waiting or recharging at $c \in \sC_{ij}$ where $t \in \sT_{ijc} \cap t \leq \floor{(s_{ijc} + u_{ijc})/T^\Delta}$}}\\
        0 & \text{otherwise}\\
        \end{cases}$\\
    $y_{ij}^\varepsilon$ & $\begin{cases}
          1 & \text{\parbox[t][][t]{0.80\textwidth}{if both trips in the compatible trip pair $(i, j) \in \sR$ are serviced by a BEB}}\\
          0 & \text{otherwise}\\
          \end{cases}$\\
    $y_{ij}^\kappa$ & $\begin{cases}
          1 & \text{\parbox[t][][t]{0.80\textwidth}{if both trips in the compatible trip pair $(i, j) \in \sR$ are serviced by a DB}}\\
          0 & \text{otherwise}\\
          \end{cases}$\\
    \bottomrule
  \end{tabularx}
\end{table}

In this problem, the planning horizon $\overline{T}$ is discretized into steps of duration $T^\Delta$. Let $\sT$ represent the set of time indices, where each $t \in \sT$ corresponds to the timestamp $D_t = t T^\Delta$. A vehicle may commence operations at any point within $\overline{T}$ and follows a periodic schedule. The continuous variable $t_i^\eta$ tracks the elapsed time from the moment a bus departs the garage for its first trip until the start of trip $i \in \sI$. For BEBs, the continuous variables $s_{ijc}$ and $u_{ijc}$ define the start time and duration of a recharging event at station $c$ between trips $i$ and $j$, respectively. Consequently, the recharging event concludes at $s_{ijc} + u_{ijc}$.

To model charging occupancy, we introduce the binary variable $x_{ijct}$. This variable equals 1 if a BEB covers trip pair $(i,j)$ using station $c \in \sC_i$ and is actively charging during time step $t \in \sT$ (i.e., $D_t \in [s_{ijc}, s_{ijc} + u_{ijc}]$); otherwise, $x_{ijct} = 0$. We strictly enforce capacity constraints: the total number of buses charging at station $c$ at any time $t$ cannot exceed the number of available plugs $Z_c$. If demand exceeds capacity, buses must wait. Importantly, this system does not enforce an FCFS discipline. Instead, the optimization model dynamically allocates plugs to maximize overall efficiency, allowing vehicles with tighter schedules or shorter dwell times to be prioritized.

For each compatible trip pair $(i, j)$ and reachable station $c \in \sC_i$, let $\sT_{ijc}$ denote the subset of time steps falling within the maximum possible dwell interval (i.e., $T_i^\beta + T_{ic}^\delta \leq D_t \leq T_j^\alpha - T_{cj}^\delta$). To linearize the charging constraints, similar to \citet{DAVATGARI2024953}, we employ the auxiliary binary variables $x_{ijct}^\alpha$ and $x_{ijct}^\beta$. The variable $x_{ijct}^\alpha$ equals 1 for all time steps strictly preceding the start of charging ($0 \leq D_t < s_{ijc}$). Similarly, $x_{ijct}^\beta$ equals 1 for all time steps up to and including the end of the charging session ($0 \leq D_t \leq s_{ijc} + u_{ijc}$). By imposing the relationship $x_{ijct} = x^\beta_{ijct} - x^\alpha_{ijct}$, we ensure that $x_{ijct} = 1$ exclusively during the active charging window. This formulation provides a rigorous method for tracking charger occupancy across continuous time intervals within a discrete framework.

The objective \labelcref{objective_function} is to minimize the total system cost, comprising fixed capital expenditures and variable operational expenses, while satisfying the minimum BEB fleet share requirements. Fixed costs include the acquisition of BEBs and DBs. Operational costs aggregate the expenses for revenue service, deadheading (to trips, garages, or stations), fuel, and maintenance. The variable $v$ quantifies the shortfall in the BEB fleet count relative to the target share $A^\nu$, and the parameter $P$ applies a penalty cost per unit of violation. The complete MILP formulation is presented below:

{\begin{equation} \label[objfunc]{objective_function}
\begin{split}
    & \begin{multlined}
    \min \sum_{i\in\sI} \left( V^\varepsilon y_{is}^\varepsilon
    + V^\kappa y_{is}^\kappa \right) + P (v + v^\prime) \\
    + \sum_{i \in \sI} \sum_{j \in \sI_i^\omega} \left(T_i^\beta - T_i^\alpha\right) \left(F^\varepsilon y^\varepsilon_{ij} + F^\kappa y^\kappa_{ij}\right) + \sum_{\substack{{\left(i,j\right)}\in \sR}}  T_{ij}^\delta \left( F^\varepsilon y_{ij}^\varepsilon + F^\kappa y_{ij}^\kappa \right) \\
    + \sum_{i \in \sI} {\sum_{\substack{c\in \sC_i}}} \sum_{j \in \sI_{ic}^\sigma} F^\varepsilon \left( -T_{ij}^\delta + T_{ic}^\delta + T_{cj}^\delta \right){q_{ijc}^\varepsilon} 
    + \sum_{i \in \sI} \sum_{j \in \sI_{is}^\sigma} F^\kappa \left(-T_{ij}^\delta + T_{is}^\delta + T_{sj}^\delta \right) q_{ij}^\kappa \\
    \end{multlined}
\end{split}
\end{equation}}
\noindent subject to,

{\begin{equation}\label[consset]{single_vehicle_1}
    \sum_{j \in \sI_i^\omega} \left(y_{ij}^\kappa + y_{ij}^\varepsilon\right) = 1 \qquad \forall i\in\sI
\end{equation}}

{\begin{equation}\label[consset]{single_vehicle_2}
    \sum_{i\in \sI_j^\tau} \left(y_{ij}^\kappa + y_{ij}^\varepsilon \right) = 1 \qquad \forall j\in\sI
\end{equation}}

{\begin{equation}\label[consset]{vehicle_balance}
    \sum_{i\in \sI_j^\tau}y_{ij}^\kappa = \sum_{k \in \sI_j^\omega}y_{jk}^\kappa \qquad \forall j\in\sI
\end{equation}}

{\begin{equation}\label[cons]{share_of_ev}
    {v \geq \sum_{i \in \sI} \left(A^\nu (y_{is}^\kappa + y_{is}^\varepsilon) - y_{is}^\varepsilon\right)}
\end{equation}}

{\begin{equation}\label[cons]{share_of_ev_times}
     {v^\prime \geq A^\tau \sum_{i \in \sI} \left(T_i^\beta - T_i^\alpha\right) - \sum_{i \in \sI} \sum_{j \in \sI_i^\omega}\left(T_i^\beta - T_i^\alpha\right) y^\varepsilon_{ij}}
\end{equation}}

{\begin{equation}\label[consset]{sij}
   s_{ijc} \geq (T_i^\beta + T_{ic}^\delta) q_{ijc}^\varepsilon \qquad \forall i \in \sI, c \in \sC_i, j \in \sI_{ic}^\sigma
\end{equation}}

{\begin{equation}\label[consset]{eij}
   {s_{ijc} + u_{ijc}} \leq (T_j^\alpha - T_{cj}^\delta) q_{ijc}^\varepsilon \qquad \forall i \in \sI, c \in \sC_i, j \in \sI_{ic}^\sigma
\end{equation}}

{\begin{equation}\label[consset]{xalpha1}
   s_{ijc} \leq D_t + T^\Delta - \epsilon + \bigM(1 - q_{ijc}^\varepsilon + x_{{ij}ct}^\alpha) \qquad \forall i \in \sI, c \in \sC_i, j \in \sI_{ic}^\sigma, t \in \sT_{ijc}
\end{equation}}

{\begin{equation}\label[consset]{xalpha2}
   s_{ijc} \geq D_t + T^\Delta - \bigM(2 - q_{ijc}^\varepsilon - x_{{ij}ct}^\alpha) \qquad \forall i \in \sI, c \in \sC_i, j \in \sI_{ic}^\sigma, t \in \sT_{ijc}
\end{equation}}

{\begin{equation}\label[consset]{xbeta1}
   {s_{ijc} + u_{ijc}} \leq D_t - \epsilon + \bigM(1 - q_{ijc}^\varepsilon + x_{{ij}ct}^\beta) \qquad \forall i \in \sI, c \in \sC_i, j \in \sI_{ic}^\sigma, t \in \sT_{ijc}
\end{equation}}

{\begin{equation}\label[consset]{xbeta2}
   {s_{ijc} + u_{ijc}} \geq D_t - \bigM(2 - q_{ijc}^\varepsilon - x_{{ij}ct}^\beta) \qquad \forall i \in \sI, c \in \sC_i, j \in \sI_{ic}^\sigma, t \in \sT_{ijc}
\end{equation}}

{\begin{equation}\label[consset]{xijct}
   x_{{ij}ct} = x_{{ij}ct}^\beta - x_{{ij}ct}^\alpha  \qquad \forall i \in \sI, c \in \sC_i, j \in \sI_{ic}^\sigma, t \in \sT_{ijc}
\end{equation}}

{\begin{equation}\label[consset]{uij_lower}
    u_{ijc} \geq T^\Delta (\sum_{t \in \sT_{ijc}} x_{{ij}ct} - 2) + \epsilon \qquad \forall i \in \sI, c \in \sC_i, j \in \sI_{ic}^\sigma
\end{equation}}

{\begin{equation}\label[consset]{SoC_for_diesel}
    b_i \geq \overline{B} - \bigM(1 - \sum_{j \in \sI^\omega_i}y^\kappa_{ij}) \qquad \forall i \in \sI
\end{equation}}

{\begin{equation}\label[consset]{energy_first_trip_from_garage1}
   b_i \leq B^\iota - T_{si}^\delta + (1 - y_{si}^\varepsilon) (\overline{B} - B^\iota + T_{si}^\delta) \qquad \forall i\in\sI
\end{equation}}

{\begin{equation}\label[consset]{energy_first_trip_from_garage2}
   b_i \geq B^\iota - T_{si}^\delta + (1 - y_{si}^\varepsilon) (\underline{B} - B^\iota + T_{si}^\delta) \qquad \forall i\in\sI
\end{equation}}

{\begin{equation}\label[consset]{soc_trip1}
   b_j \leq b_i - B_i - T_{ij}^\delta - \sum_{c \in \sC_{ij}} \left((T_{ic}^\delta + T_{cj}^\delta - T_{ij}^\delta) q_{ijc}^\varepsilon - {R_c \times} u_{ijc} \right) + \bigM (1 - y_{ij}^\varepsilon) \qquad \forall i \in \sI, j \in \sI_i
\end{equation}}

{\begin{equation}\label[consset]{soc_trip2}
   b_j \geq b_i - B_i - T_{ij}^\delta - \sum_{c \in \sC_{ij}} \left((T_{ic}^\delta + T_{cj}^\delta - T_{ij}^\delta) q_{ijc}^\varepsilon - {R_c \times} u_{ijc} \right) - \bigM (1 - y_{ij}^\varepsilon) \qquad \forall i \in \sI, j \in \sI_i
\end{equation}}

{\begin{equation}\label[consset]{energy_gained_at_station}
    b_i - B_i - \sum_{c \in \sC_i} \sum_{j \in \sI_{ic}^\sigma} \left(T_{ic}^\delta q_{ijc}^\varepsilon - {R_c \times} u_{ijc} \right) \leq \overline{B} \qquad \forall i\in\sI
\end{equation}}

{\begin{equation}\label[consset]{trip_to_garage_energy}
   b_i \geq \underline{B} + B_i + \sum_{c \in \sC_i} \sum_{j \in \sI_{ic}^\sigma} T_{ic}^\delta q_{ijc}^\varepsilon - \bigM \sum_{j \in I_i^\omega} y_{ij}^\kappa   
   \qquad \forall i\in\sI
\end{equation}}

{\begin{equation}\label[consset]{energy_last_trip_before_garage}
   b_i \geq (\underline{B} + B_i + T_{is}^\delta)  y_{is}^\varepsilon \qquad \forall i\in\sI
\end{equation}}

{\begin{equation}\label[consset]{first_trip_time}
   t_i^\eta - T_{si}^\delta \leq \bigM (1 - y_{si}^\varepsilon - y_{si}^\kappa) \qquad \forall i\in\sI
\end{equation}}

{\begin{equation}\label[consset]{last_trip_time}
   t_i^\eta + B_i + T_{is}^\delta \leq \overline{T} + \bigM (1 - y_{is}^\varepsilon - y_{is}^\kappa) \qquad \forall i\in\sI
\end{equation}}

{\begin{equation}\label[consset]{time_consecutive_trips}
   t_j^\eta - t_i^\eta \geq T_j^{\alpha} - T_i^{\alpha} - \bigM (1 - y_{ij}^\varepsilon - y_{ij}^\kappa) \quad \forall i \in \sI, j \in \sI_i
\end{equation}}

{\begin{equation}\label[consset]{charge_time}
   b_i - B_i - T_{is}^\delta + R_s \left(\overline{T} - (t_i^\eta + B_i + T_{is}^\delta)\right) \geq B^\iota - \bigM (1 - y_{is}^\varepsilon) \qquad \forall i\in\sI
\end{equation}}

{\begin{equation}\label[consset]{charger_max_cap}
    \sum_{i \in \sI} \sum_{j \in \sI_{ict}^\sigma} x_{{ij}ct} \leq Z_c \qquad \forall c \in \sC, t \in \sT
\end{equation}}

{\begin{equation}\label[consset]{qij_lower}
    x_{ijct} \leq q_{ijc}^\varepsilon \qquad \forall i \in \sI, c \in \sC_i, j \in \sI_{ic}^\sigma, t \in \sT_{ijc}
\end{equation}}

{\begin{equation}\label[consset]{garage_visit_vars_CV}
  q_{ij}^\kappa \leq y_{ij}^\kappa \qquad \forall i \in \sI, j \in \sI_{is}^\sigma
\end{equation}}

{\begin{equation}\label[consset]{charger_visit_vars_EV}
  \sum_{c \in \sC_{ij}} q_{ijc}^\varepsilon \leq y_{ij}^\varepsilon \qquad \forall i \in \sI, j \in \sI_i
\end{equation}}

{\begin{equation}\label[consset]{garage_visit_vars_2_CV}
  q_{ij}^\kappa \geq y_{ij}^\kappa \qquad \forall i \in \sI, j \in \sI_{is}^\sigma \setminus \sI_{i}^\lambda
\end{equation}}

{\begin{equation}\label[consset]{charger_visit_vars_2_EV}
  \sum_{c \in \sC_{ij}} q_{ijc}^\varepsilon \geq y_{ij}^\varepsilon \qquad \forall i \in \sI, j \in \sI_{i} \setminus \sI_{i}^\lambda
\end{equation}}

{\begin{equation}\label[consset]{vars_y}
   x_{{ij}ct}, x_{{ij}ct}^\alpha, x_{{ij}ct}^\beta \in \{0,1\}~\forall i \in \sI, c \in \sC_i, j \in \sI_{ic}^\sigma, t \in \sT_{ijc};\quad y_{ij}^\varepsilon, y_{ij}^\kappa \in \{0,1\}~\forall (i,j) \in \sR
\end{equation}}

{\begin{equation}\label[consset]{vars_time}
    {v, v^\prime, t_i^\eta, b_i, s_{ijc},} u_{ijc} \in \mathbb{R}^+  \qquad \forall i \in \sI, c \in \sC_i, j \in \sI_{ic}^\sigma
\end{equation}}

{\begin{equation}\label[consset]{vars_q_e}
   q_{ijc}^\varepsilon \in \{0,1\}~\forall i \in \sI, c \in \sC_i, j \in \sI_{ic}^\sigma; \quad q_{ij}^\kappa \in \{0,1\}~\forall i \in \sI, j \in \sI_{is}^\sigma
\end{equation}}

Constraints \labelcref{single_vehicle_1} and \labelcref{single_vehicle_2} ensure that each trip is assigned to exactly one vehicle. Constraints \labelcref{vehicle_balance} guarantee that all trips within a feasible run are serviced by the same vehicle type. Constraint \labelcref{share_of_ev} defines the variable \( v \), representing the shortfall in the number of BEBs needed to meet the desired \( A^\nu \)\% fleet share, thereby enforcing a minimum BEB requirement. Similarly, constraint \labelcref{share_of_ev_times} defines the variable \( v^\prime \), denoting the shortfall in the percent of service time needed to meet the desired \( A^\tau \)\% of the total revenue trip time to be conducted by BEBs.

Constraints \labelcref{sij} define the start time of recharging, ensuring that the bus has sufficient time to reach the charging station after completing its trip. Constraints \labelcref{eij} determine the charging duration, ensuring that, after recharging, the bus has enough time to deadhead to the starting point of its next trip. Constraints \labelcref{sij} and \labelcref{eij} collectively ensure that the charging duration remains within the bus's layover time at the charging station. Constraints \labelcref{xalpha1} and \labelcref{xalpha2} define \( x_{ijct}^\alpha \), while constraints \labelcref{xbeta1} and \labelcref{xbeta2} define \( x_{ijct}^\beta \). Constraint \labelcref{xijct} links \(x_{ijct}\) with \(x_{ijct}^\alpha\) and \(x_{ijct}^\beta\). Additionally, constraint \labelcref{uij_lower} synchronizes the continuous charging duration $u_{ijc}$ with the discrete time-step counters, ensuring the reserved time slots do not exceed the actual charged energy.

Constraints \labelcref{SoC_for_diesel} ensure that the SoC level for all DBs is maintained at its maximum level ($\overline{B}$). Constraints \labelcref{energy_first_trip_from_garage1} and \labelcref{energy_first_trip_from_garage2} specify the SoC when the bus departs from the garage, setting it to $B^\iota$. Constraints \labelcref{soc_trip1} and \labelcref{soc_trip2} determine the SoC at the start of trips that are not the first in the chain, while constraints \labelcref{energy_gained_at_station} ensure that the energy gained at the charging station does not exceed the battery’s capacity. Constraints \labelcref{trip_to_garage_energy} ensure the SoC always remains above the minimum threshold, and constraints \labelcref{energy_last_trip_before_garage} set it at the end of the last trip in a BEB run. Constraints \labelcref{first_trip_time} initialize the run time $t_i^\eta$ when a bus leaves the garage to start its run. Constraints \labelcref{last_trip_time} limit the total run time of a vehicle to be less than \( \overline{T} \). Constraints \labelcref{time_consecutive_trips} relate the elapsed time between two consecutive trips. Constraints \labelcref{charge_time} ensure a BEB has adequate time within its run to recharge to achieve an SoC level of \( B^\iota \).

Constraints \labelcref{charger_max_cap} ensure that the total number of buses charging at any given time does not exceed the facility's capacity. Constraints \labelcref{qij_lower} prevent buses from recharging at a station unless they actually visit it. Constraints \labelcref{garage_visit_vars_CV} are logical constraints ensuring that if a DB does not cover a trip pair, it cannot visit the garage between those trips. Similarly, constraints \labelcref{charger_visit_vars_EV} maintain the same logic for BEBs visiting charging stations, additionally ensuring that BEBs don't recharge at multiple charging stations. Lastly, constraints \labelcref{garage_visit_vars_2_CV} and \labelcref{charger_visit_vars_2_EV} ensure that if a vehicle covers a trip pair with a sufficiently long inter-trip layover, greater than $L$, it must visit either a garage or a charging station during the layover. Constraints \labelcref{vars_y}--\labelcref{vars_q_e} indicate variable domains.

The proposed model minimizes complexity by excluding explicit decision variables for waiting times and layovers occurring immediately before or after charging events. Instead, these operational metrics are derived via post-processing. Consider a BEB traversing the arc $(i,j)$ and visiting station $c$ (i.e., $q^\varepsilon_{ijc} = 1$). The vehicle arrives at the station at time $T^\beta_{i} + T^\delta_{ic}$. The interval between arrival and the start of charging, defined as $[T^\beta_{i} + T^\delta_{ic}, s_{ijc}]$, constitutes the total pre-charging dwell time. Within this window, any time steps where the station capacity is fully utilized (i.e., $\sum x_{ijct} = Z_c$) are classified as \textit{waiting time} due to congestion. Conversely, intervals where plugs are available are categorized as \textit{pre-charging layover}. Upon completion of the charging process at $s_{ijc} + u_{ijc}$, the vehicle remains at the station in a \textit{post-charging layover} state until the mandatory departure time $T^\alpha_j - T^\delta_{cj}$, ensuring it reaches the origin of trip $j$ on schedule. This approach allows for the accurate reporting of service quality metrics without increasing the computational burden of the optimization model.

\section{Solution method}\label[sec]{solution}
The proposed model exhibits significant complexity due to the large number of decision variables associated with integrating bus scheduling and charging strategies. Much of this complexity arises from the intricate charging-related variables and the extensive set of time intervals \(\sT\), substantially expanding the problem size and making exact optimization challenging, particularly for large-scale case studies. To efficiently address large problem instances, a CG method is proposed. The CG is commonly utilized in the literature to solve BEB scheduling and charging optimization problems \citep{lin2016column,sundin2018scheduling, DUAN2023104175, GERBAUX2025106848, XU2025111138}. 

The key idea in CG is to decompose the problem into a Restricted Master Problem (RMP) and one or more Subproblems (SPs). The RMP solves a reduced version of the original problem, considering only a limited subset of feasible vehicle schedules (referred to as columns). The SPs are then used to identify potentially beneficial columns, that is, bus schedules with negative reduced cost, that could improve the RMP solution. This iterative process continues until no improving columns can be found, at which point optimality is achieved or the algorithm terminates because of a predefined time limit.

In this study, a column represents a complete, feasible daily schedule for a single bus, beginning and ending at the garage. Each schedule comprises a sequence of assigned trips, potential deadheading between trips, and visits to the garage or charging stations as required. We distinguish between two types of columns: DB columns represent DB schedules and may include deadhead movements and layovers at the garage. BEB columns represent BEB schedules, incorporating not only trip sequences and deadheading but also charging activities and compliance with battery SoC constraints.


\subsection{Restricted master problem}\label[sec]{sec:RMP}
The RMP is formulated as a set-covering model with four key goals: (1) to minimize the total number of buses, while minimizing overall operational costs, (2) to select a combination of bus schedules that collectively cover all trips in $\sI$, (3) to ensure compliance with charging capacity constraints, and (4) to meet the minimum required share of BEBs in the fleet.

Let $\sB$ denote the set of all feasible columns. Similarly, $\sB^\kappa$ and $\sB^\varepsilon$ define the subset of feasible DB and BEB columns, respectively. The RMP is solved over a limited subset $\sB' \subset \sB$, seeking the best solution within the current pool of available schedules. The operational cost of a given column includes the cost of acquiring a vehicle to serve the schedule, the travel time costs of assigned trips, and the deadhead costs to garages or charging stations. The operational costs of a DB column, denoted by $C_r^\kappa$, and a BEB column, denoted by $C_r^\varepsilon$, are defined in \labelcref{cost_of_diesel_block} and \labelcref{cost_of_electric_block}, respectively. In these formulations, the binary parameter $Y_{rij}^\kappa = 1$ if DB column $r$ includes a deadhead movement from trip $i$ to trip $j$, and binary parameter $Q_{rij}^\kappa = 1$ if the vehicle assigned to column $r$ visits the garage between trips $i$ and $j$. Parameters $Y_{rij}^\varepsilon$ follow a similar definition for BEBs. Parameter $Q_{rijc}^\varepsilon = 1$ if the vehicle assigned to column $r$ visits charging station $c$ between trips $i$ and $j$. Let $S_{ri}^\varepsilon$ be a binary parameter equal to 1 if BEB column $r$ services trip $i \in \sI$, and 0 otherwise. Its DB version is defined as $S_{ri}^\kappa$.

\begin{equation}\label[eq]{cost_of_diesel_block}
\begin{split}
    & \begin{multlined}
    C^\kappa_r = V^\kappa + F^\kappa \left[ \sum_{i \in \sI} S_{ri}^\kappa (T^\beta_i - T^\alpha_i) + \sum_{(i, j) \in \sR} Y^\kappa_{rij} T^\delta_{ij} +
    \sum_{i \in \sI} \sum_{j \in \sI_{is}^\sigma} Q_{rij}^\kappa (-T^\delta_{ij} + T^\delta_{is} + T^\delta_{sj}) \right]
    \end{multlined}
\end{split}
\end{equation}

\begin{equation}\label[eq]{cost_of_electric_block}
\begin{split}
    & \begin{multlined}
    C^\varepsilon_r = V^\varepsilon + F^\varepsilon \left[ \sum_{i \in \sI} S_{ri}^\varepsilon (T^\beta_i - T^\alpha_i) + \sum_{(i, j) \in \sR} Y^\varepsilon_{rij} T^\delta_{ij} +
    \sum_{i \in \sI} \sum_{c\in \sC_i} \sum_{j \in \sI_{ic}^\sigma} Q_{rijc}^\varepsilon (-T^\delta_{ij} + T^\delta_{ic} + T^\delta_{cj}) \right]
    \end{multlined}
\end{split}
\end{equation}

We define the binary decision variable $z_r^\varepsilon$ that equals 1 if BEB column $r$ is selected and 0 otherwise. We also define its DB counterpart as $z_r^\kappa$. Finally, let $H_{rct} = 1$ if BEB schedule $r$ utilizes a charger at station $c \in \sC$ during time step $t \in \sT$, and 0 otherwise. The RMP formulation is as follows:

\begin{equation}\label[objfunc]{master_problem}
    \min \sum_{r\in\sB^{\varepsilon}} C^{\varepsilon}_r z^{\varepsilon}_r + \sum_{r\in\sB^{\kappa}} C^{\kappa}_r z^{\kappa}_r + Pv
\end{equation}

\noindent subject to,
\begin{equation}\label[consset]{trip_coverage}
    \sum_{r \in \sB^\varepsilon} S^\varepsilon_{ri}z^\varepsilon_r + \sum_{r \in \sB^\kappa} S^\kappa_{ri}z^\kappa_r = 1 \qquad \forall i\in\sI \quad (\pi_i)
\end{equation}

\begin{equation}\label[cons]{BEB_minimum_share_fleet}
    A^\nu (\sum_{r \in \sB^\varepsilon} z^\varepsilon_r + \sum_{r \in \sB^\kappa} z^\kappa_r) - \sum_{r \in \sB^\varepsilon} z^\varepsilon_r \leq v \quad (\alpha)
\end{equation}

\begin{equation}\label[consset]{charger_plug_capacity}
    \sum_{r \in \sB^\varepsilon} H_{rct} z^\varepsilon_r \leq Z_c \quad \forall c \in \sC, t \in \sT \quad (\phi_{ct})
\end{equation}

\begin{equation}\label[consset]{CG_binary_vars}
    v \geq 0, z^\varepsilon_r, z^\kappa_r \in \{0,1\} \qquad \forall r \in \sB^\varepsilon, r \in \sB^\kappa
\end{equation}

Objective function \labelcref{master_problem} minimizes the total daily cost of operating the selected DB and BEB column, while imposing a penalty for any shortfall from the target BEB fleet share. Constraints \labelcref{trip_coverage} define the trip coverage constraints, ensuring that each scheduled trip is performed exactly once. The corresponding dual variables are denoted by $\pi_i$. Constraints \labelcref{BEB_minimum_share_fleet} introduce a penalty variable that activates when the number of BEBs in the fleet falls below the target share $A^\nu$; the associated dual variable is $\alpha$. Note that we dropped $v^\prime$ from the model because $v$ readily functions in a similar fashion, and dropping this variable relaxes the problem to some extent. One can embed $v^\prime$ into the model similar to the way $v$ is used but should also acknowledge the fact that it adds another dimensionality to the difficulty of the problem. Constraints \labelcref{charger_plug_capacity} enforce charging station capacity limits by ensuring that the number of BEBs charging at any station at a given time step does not exceed the station’s available plug capacity. The dual variables associated with these constraints are denoted by $\phi_{ct}$. Constraints \labelcref{CG_binary_vars} enforce the binary nature of the variables $z^\varepsilon_r$ and $z^\kappa$.

While the RMP is ultimately an integer program, the CG requires a modification to its structure. The binary nature of the decision variables, $z^\kappa_r$ and $z^\varepsilon_r$, must be temporarily relaxed, allowing them to take continuous values between 0 and 1, as defined by the constraints \labelcref{CG_relaxed_consset}. This step is essential because it enables the RMP to be solved as a Linear Program (LP). The solution to this LP provides the crucial dual values associated with the RMP's constraints. These duals are then used as pricing signals in the SPs to identify columns with negative reduced costs. After the iterative CG process terminates, the integrality constraints \labelcref{CG_binary_vars} are re-imposed. The final RMP, now containing a rich set of profitable routes, is solved one last time as an MILP to find the optimal integer solution.

\begin{equation}\label[consset]{CG_relaxed_consset}
    0 \leq z^\varepsilon_r, z^\kappa_r \leq 1 \qquad \forall r \in \sB^\varepsilon, r \in \sB^\kappa
\end{equation}

\subsection{Pricing subproblems}\label[sec]{sec:pricing_sub_problems}
The pricing SPs identify valid schedules with negative reduced costs. The reduced cost of a column, denoted by $\bar{c}_r^\kappa$ for DBs and $\bar{c}_r^\varepsilon$ for BEBs, is calculated as its total operational cost minus the sum of the dual values corresponding to the RMP constraints it satisfies. If a column with negative reduced cost is found, it is added to the RMP, which is then re-optimized. Due to the presence of two distinct fleet types, separate pricing SPs are formulated: one for DBs called DB-SP and a more complex one for BEBs called BEB-SP. Each SP searches its respective solution space for feasible column that, if incorporated into the RMP, could potentially improve the overall objective.

\subsubsection{Diesel bus subproblem}
The DB-SP is formulated as a shortest path problem that seeks to find the DB column with the minimum reduced cost. The objective function of the DB-SP minimizes the sum of a DB's operational costs minus the dual values associated with the trips it serves. The decision variables in this model primarily determine the sequence of trips, represented by binary variables $y^\kappa_{ij}$ that indicate travel from one trip to another, and $q^\kappa_{ij}$ for intermediate garage visits. The model is subject to a set of constraints that ensure the generated column is a valid sequence of trips, starting and ending at the main garage, with proper flow conservation between trips. Note in the below formulation that (1) values to $S_{ri}^\kappa$ and $Y_{rij}^\kappa$ in the RMP are deduced from $y^\kappa_{ij}$ variables, (2) values to $Q_{rij}^\kappa$ in the RMP are obtained from $q^\kappa_{ij}$ values, and (3) objective function \labelcref{DB-SP-OBJ} is to minimize the reduced cost $\bar{c}_r^\kappa = C_r^\kappa-\sum_{i\in\sI} S_{ri}^\kappa \pi_i - A^\nu \alpha$ and is written along with its subjected constraints to clearly define the DB-SP with its own set of decision variables.

\begin{equation}\label[objfunc]{DB-SP-OBJ}
\begin{multlined}
    \min \bar{c}_r^\kappa = V^\kappa + F^\kappa [ \sum_{i \in \sI} \sum_{j \in \sI_i^\omega} y_{ij}^\kappa (T^\beta_i - T^\alpha_i) + \sum_{(i, j) \in \sR} y^\kappa_{ij} T^\delta_{ij} \\ +
    \sum_{i \in \sI} \sum_{j \in \sI_{is}^\sigma} q_{ij}^\kappa (-T^\delta_{ij} + T^\delta_{is} + T^\delta_{sj}) ] -
    \sum_{i\in \sI} \pi_i \sum_{j \in \sI_i^\omega} y_{ij}^\kappa - A^\nu \alpha
\end{multlined}
\end{equation}

\noindent subject to \labelcref{vehicle_balance}, \labelcref{garage_visit_vars_CV}, \labelcref{garage_visit_vars_2_CV}, \labelcref{vars_y}, \labelcref{vars_q_e} and
\begin{equation}\label[cons]{DB-SP-flow_conservation1}
    \sum_{j \in \sI} y^\kappa_{sj} = 1; \quad \sum_{i \in \sI} y^\kappa_{is} = 1
\end{equation}
\subsubsection{Battery electric bus subproblem via SPPRC}
The BEB-SP is modeled as an SPPRC. To represent the subproblem that minimizes $\bar{c}^\varepsilon_r$, we define variables: Binary variables $q_{ijc}^\varepsilon=1$ if a BEB visits charging station $c$ between trips $i$ and $j$ and determines the values for $Q_{rijc}^\varepsilon$ in the RMP. Binary variables $h_{ct}=1$ if a BEB uses charger at station $c$ during time step $t$ and sets values for $H_{rct}$ in the RMP. The reduced cost for a BEB column is $\bar{c}^\varepsilon_r = C_r^\varepsilon - \sum_{i \in \sI} S_{ri}^\varepsilon \pi_i - (A^\nu - 1) \alpha - \sum_{c \in \sC} \sum_{t \in \sT} H_{rct} \phi_{ct}$.

\begin{equation}\label[objfunc]{BEB-SP-OBJ}
\begin{multlined}
    \min \bar{c}_r^\varepsilon = V^\varepsilon + F^\varepsilon [ \sum_{i \in \sI} \sum_{j \in \sI_i^\omega} y_{ij}^\varepsilon (T^\beta_i - T^\alpha_i) + \sum_{(i, j) \in \sR} y^\varepsilon_{ij} T^\delta_{ij} \\+
    \sum_{i \in \sI} \sum_{c\in \sC_i} \sum_{j \in \sI_{ic}^\sigma} q_{ijc}^\varepsilon (-T^\delta_{ij} + T^\delta_{ic} + T^\delta_{cj}) ] -
    \sum_{i\in \sI} \pi_i \sum_{j \in \sI_i^\omega} y_{ij}^\varepsilon - (A^\nu - 1) \alpha - \sum_{c \in \sC} \sum_{t \in \sT} h_{ct} \phi_{ct}
\end{multlined}
\end{equation}

\noindent subject to \labelcref{sij}, 
\labelcref{eij}, 
\labelcref{energy_first_trip_from_garage1} --
\labelcref{charge_time}, 
\labelcref{charger_visit_vars_EV}, 
\labelcref{charger_visit_vars_2_EV} --
\labelcref{vars_q_e} and

\begin{equation}\label[cons]{BEB-SP-flow_conservation1}
    \sum_{j \in \sI} y^\varepsilon_{sj} = 1; \quad \sum_{i \in \sI} y^\varepsilon_{is} = 1
\end{equation}

\begin{equation}\label[cons]{BEB-SP-flow_conservation3}
    \sum_{i \in \sI^\tau_j} y^\varepsilon_{ij} = \sum_{k \in \sI^\omega_j} y^\varepsilon_{jk} \qquad \forall j \in \sI
\end{equation}

Since $\phi_{ct}$ varies by time step, the cost of traversing a charging arc in the SPPRC
network is time-dependent. This requires the algorithm to track time as a resource and
evaluate the cost of charging based on the exact arrival and duration at the station. While
SPPRC provides an exact method for finding the column with the minimum reduced cost,
solving it to optimality becomes computationally prohibitive for the large-scale networks.

\subsubsection{Heuristic pricing for large-scale instances}\label[sec]{heuristic_pricing}
To address the computational challenges of the exact SPPRC on large instances, we implement a two-stage heuristic pricing approach. This method decomposes the subproblem into schedule generation and charging optimization, allowing for the rapid discovery of high-quality columns with negative reduced costs. In the first stage, a greedy constructive heuristic generates a candidate sequence of trips. Starting from the garage, the next trip in the sequence is selected probabilistically. The selection weights are derived from the dual values ($\pi_i$) and temporal compatibility, prioritizing trips that are valuable to the RMP and fit within the schedule. In the second stage, we solve a subordinate MILP to determine the optimal charging strategy for a fixed trip chain generated in Stage 1. This model decides where and for how long the bus should charge to ensure feasibility with respect to battery limits ($\underline{B}, \overline{B}$) while minimizing the charging cost component $\sum_{c\in \sC, t\in\sT} H_{rct} \phi_{ct}$. This reduced model is significantly smaller than the full SPPRC, as the sequence of visits is fixed, making it solvable in negligible time. If a feasible charging schedule with a negative reduced cost is found, the column is added to the RMP.

\subsection{Column generation procedure}\label[sec]{sec:CG_algorithm}
The CG algorithm is an iterative procedure that cycles between solving the RMP and the SPs. Because the final solution requires an integer number of buses for each schedule, this process is embedded within a branch-and-bound search, forming a complete branch-and-price algorithm. The core CG procedure executed at each node of the branch-and-bound tree is detailed in \cref{cg_algorithm}. 

\begin{algorithm}[!ht]
\caption{Column generation for mixed-fleet bus scheduling}\label[algo]{cg_algorithm}
\SetArgSty{textnormal}
\scriptsize
\SetKwInOut{Input}{Input}
\SetKwInOut{Output}{Output}
\Input{~Set of trips $\sI$, Stations $\sC$, Parameters.}
\Output{~Integer schedule solution $\boldsymbol{z^{\kappa*}}$ and $ \boldsymbol{z^{\varepsilon*}}$.}
\SetAlgoLined
\SetKwFunction{FSolveRMP}{\textproc{SolveRMP}}
\SetKwFunction{FSolveDieselSP}{\textproc{SolveDB-SP}}
\SetKwFunction{FSolveElectricSP}{\textproc{SolveBEB-SP}}
\SetKwFunction{FSolveElectricHeuristic}{\textproc{SolveBEB-SPHeuristic}}

    ${\sB^\kappa} \gets$ Generate an initial set of feasible DB schedules covering all trips in $\sI$\;
    
    ${\sB^\varepsilon} \gets$ Generate an initial set of feasible BEB schedules that cover all trips in $\sI$, ensuring that buses have sufficient SoC to complete their assigned trips\;
    \tcp{We initialize $\sB^\kappa$ and $\sB^\varepsilon$ with singleton schedules in our implementation.}
    
    $converged \gets$ \textbf{false}; \Comment{Convergence criteria: (i) no columns, (ii) no improvement in solution, (iii) reached time limit.}
    
    \While{\textbf{not} $converged$}{
        $\pi, \alpha, \phi \gets$ \FSolveRMP{$\sB^\kappa, \sB^\varepsilon$ (Relaxed)}; \Comment{Relax integrality constraints in the RMP: $0 \leq z^\kappa, z^\varepsilon \leq 1$.}
        
        $r^\kappa, \bar{c}^\kappa \gets$ \FSolveDieselSP{$\pi, \alpha$};\Comment{Solve DB-SP to get DB column $r^\kappa$ and reduced cost $\bar{c}^\kappa$.}\
        
        \If{$\bar{c}^\kappa < -\epsilon$}{
            $\sB^\kappa \gets \sB^\kappa \cup \{r^\kappa\}$;\Comment{Add new DB schedule to the RMP}\ 
        }
        \tcp{If BEB-SP formulation is used.}
        \If{BEB-SP}{$r^\varepsilon, \bar{c}^\varepsilon \gets$ \FSolveElectricSP{$\pi, \alpha$, $\phi$}; \Comment{Solve BEB-SP to get BEB column $r^\varepsilon$ and reduced cost $\bar{c}^\varepsilon$.}\
        
        \If{$\bar{c}^\varepsilon < -\epsilon$}{
            $\sB^\varepsilon \gets \sB^\varepsilon \cup \varepsilon$; \Comment{Add new BEB schedule to the RMP}\
        }}
        \tcp{If BEB-SP heuristic is used.}
        \Else{        
        \tcp{Run heuristic multiple times to generate diverse columns}
        \For{$iter \in 1 \dots K$}{
            $r^\varepsilon, \bar{c}^\varepsilon \gets$ \FSolveElectricHeuristic{$\pi, \alpha, \phi$}; \Comment{Solve BEB-SP heuristic to get $r^\varepsilon$ and $\bar{c}^\varepsilon$.}\
            
            \If{$\bar{c}^\varepsilon < -\epsilon$}{
                $\sB^\varepsilon \gets \sB^\varepsilon \cup \{r^\varepsilon\}$; \Comment{Add new BEB schedule to the RMP}\
            }
        }}
        
        \If{No columns added \textbf{or} Stalling criterion met \textbf{or} Time limit reached}{
            $converged \gets$ \textbf{true}\;
        }
    }    
    
    $z^\kappa, z^\varepsilon \in \{0, 1\}$; \Comment{Set variables $\boldsymbol{z^\kappa}$ and $\boldsymbol{z^\varepsilon}$ in \FSolveRMP to be binary.}\

    $\boldsymbol{z^{\kappa*}}$, $\boldsymbol{z^{\varepsilon*}} \gets $ \FSolveRMP($\sB^\kappa, \sB^\varepsilon$ (Integer)); \Comment{Solve final integer RMP to get $\boldsymbol{z^{\kappa*}}$ and $\boldsymbol{z^{\varepsilon*}}$}\
    
    \Return $\boldsymbol{z^{\kappa*}}$ and $\boldsymbol{z^{\varepsilon*}}$\
\end{algorithm}

A key consideration in the implementation of the CG procedure was the handling of columns that are repeatedly generated by the SPs. Initially, a standard check was included to prevent the addition of a column identical to one already present in the RMP. However, this practice was found to cause the algorithm to terminate prematurely, even when further improvement to the solution was possible. The root cause of this behavior is degeneracy in the RMP's LP relaxation. The dual values obtained from the RMP are determined by the basis used by the simplex algorithm. This can lead to a situation where, even after a new profitable column is added, the solver pivots to a new basis that yields the exact same set of duals as the previous iteration. Consequently, the subproblem is fed identical pricing information and generates the same \textit{best} column again. When the algorithm enforced uniqueness, it would block this re-generated column and terminate.

Therefore, to overcome this cyclic behavior, the check for duplicate columns was deliberately removed from the procedure. By allowing identical columns to be added to the RMP, we provided the flexibility needed to pivot to a different basis. This change in basis produces new duals, \textit{un-sticking} the algorithm and enabling the search for other, structurally different routes to continue. While allowing duplicate columns helps mitigate some effects of degeneracy, a more persistent form of stalling was still observed during the implementation. In this particular case, the RMP objective value would plateau and fail to improve, yet the subproblems would continue to find columns with a negative reduced cost. This occurs when a newly added profitable column is not utilized in the subsequent RMP solution. This results in an unchanged basis, stagnant dual values, and consequently, the repeated generation of the same improving columns in an unproductive loop. To ensure robust convergence and prevent this infinite cycling, a direct stalling criterion was implemented in the final algorithm. This mechanism tracks the RMP's objective value across iterations. If the objective fails to improve by a predefined tolerance, for a set number of consecutive iterations, the CG is terminated. This strategy directly addresses the symptom of non-improvement and provides a practical termination condition.

\subsection{Post-processing heuristic}\label[sec]{post_processing}
Upon obtaining an integer solution from the CG procedure, a post-processing local search heuristic is applied to further improve the solution quality. The primary objective of this phase is to consolidate the schedule by dissolving inefficient, short vehicle chains and merging their trips into longer, existing chains. This step is particularly effective in reducing the total fleet size required. The procedure, outlined in \cref{post_proc}, begins by sorting all generated BEB columns ($\boldsymbol{z^{\varepsilon*}}$) in ascending order based on the number of trips they cover. The algorithm iteratively attempts to empty the shortest chains (candidates for removal) by moving their trips to longer, ``target'' chains.

\begin{algorithm}[!ht]
\caption{Post-processing chain consolidation heuristic}\label[algo]{post_proc}
\SetArgSty{textnormal}
\scriptsize
\SetKwInOut{Input}{Input}
\SetKwInOut{Output}{Output}
\Input{~Integer solution schedules $\boldsymbol{z^{\varepsilon*}}$.}
\Output{~Consolidated schedules $\boldsymbol{z^{\varepsilon**}}$.}
\SetAlgoLined
\SetKwFunction{Sort}{Sort}
\SetKwFunction{SolveMILP}{SolveChargingMILP}
\SetKwFunction{CheckCap}{CheckCapacity}

    $\sL \gets$ \Sort($\boldsymbol{z^{\varepsilon*}}$ by trip count ascending)\;
    $idx^{\alpha} \gets 0$\;
    
    \While{$idx^{\alpha} < |\sL|$ \textbf{and} Time Limit not reached}{
        $r^{\alpha} \gets \sL[idx^{\alpha}]$; \Comment{Candidate chain to dissolve}\
        
        $merged \gets$ \textbf{false}\;
        
        \For{trip $i \in r^{\alpha}$}{
            \tcp{Try to merge into longest chains first}
            \For{$idx^{\beta} \gets |\sL|-1$ \KwTo $idx^{\alpha} + 1$}{
                $r^{\beta} \gets \sL[idx^{\beta}]$\;
                $r^{\beta \prime} \gets r^{\beta} \cup \{i\}$; \Comment{Tentative insertion}\
                
                \tcp{Check feasibility using the Pricing Heuristic's MILP}
                $feasible, \text{cost}, \text{slots} \gets$ \SolveMILP($r^{\beta \prime}$)\;
                
                \If{$feasible$}{
                    \If{\CheckCap($\text{slots}$, $\sL \setminus \{r^{\alpha}, r^{\beta}\}$) \textbf{is true}}{
                        \If{Cost($r^{\beta \prime}$) $\leq$ Cost($r^{\beta}$) + Cost($r^{\alpha}$)}{
                            Update $r^{\beta} \gets r^{\beta \prime}$\;
                            Remove $i$ from $r^{\alpha}$\;
                            $merged \gets$ \textbf{true}\;
                            \textbf{break}\;
                        }
                    }
                }
            }
        }
        
        \If{$merged$}{
            \If{$r^{\alpha}$ is empty}{
                Remove $r^{\alpha}$ from $\sL$\;
            }
            $idx^{\alpha} \gets 0$; \Comment{Restart search from the new shortest chain}\
        }
        \Else{
            $idx^{\alpha} \gets idx^{\alpha} + 1$; \Comment{Move to next candidate}\
        }
    }
    \Return $\sL$\;
\end{algorithm}

For a selected trip $i$ in a candidate short chain $r^{\alpha}$ and a target long chain $r^{\beta}$, the algorithm checks if $i$ can be feasibly inserted into $r^{\beta}$. This feasibility check is rigorous: it re-solves the charging schedule optimization (described in \cref{heuristic_pricing}) for the proposed merged chain. This ensures that the new chain, containing trip $i$, still satisfies all battery SoC ($B_i$) and time window constraints. Furthermore, the algorithm explicitly checks that the new charging times determined by the MILP do not violate station capacity constraints ($Z_c$) when considered alongside all other unchanged schedules. If the merge is feasible and results in a cost reduction (or neutral cost change with fleet reduction), the move is accepted, effectively densifying the schedule.

\section{Numerical experiments}\label[sec]{experiments}
This section presents a comprehensive numerical evaluation of the proposed optimization framework. 
All computations were performed on a high-performance computing system equipped with an Intel\textsuperscript{\textregistered} Xeon\texttrademark{} Gold 6248R processor (3.00 GHz, 24 cores, 48 logical processors) and 354 GiB of RAM, running Ubuntu 22.04. The models were implemented in Python 3.10.12 using the Gurobi 12.0.1 solver \citep{gurobi}.

\subsection{Design of experiments}\label[sec]{design_experiments}
This section outlines the input data, operational constraints, and parameter settings that form the baseline for the experimental analysis. Unless explicitly stated otherwise in subsequent sections, the values defined here constitute the default experimental environment. Notably, to rigorously evaluate the solution methods under the most demanding operational conditions, the baseline fleet composition is set to fully electric ($A^\nu = 100\%$). This configuration maximizes computational complexity by strictly enforcing all battery management and charging scheduling constraints, effectively stress-testing the optimization framework. The study focuses on the Chicago, IL metropolitan region, utilizing data from its two primary bus operators: the CTA and Pace Suburban Bus. This dual-agency scope captures a diverse range of operational environments, from the high-density urban core served by CTA to the suburban and inter-municipal routes served by Pace. The combined network comprises 10 CTA garages and 7 Pace garages.

Timetabled trip data, including start and end locations and times ($T^\alpha_i, T^\beta_i$), is sourced from the publicly available General Transit Feed Specification (GTFS) feeds \citep{GTFS}. The baseline instance consists of 25,999 daily scheduled trips connecting 1,565 terminal stops. To estimate deadhead movements between terminal stops and garages, we calculate Euclidean distances and assume a constant non-revenue speed of 20 mph.

Since the public data does not explicitly link trips to specific home garages, we employed a spatial assignment heuristic to distribute the workload. For each trip, we calculated the Euclidean midpoint between the start and end coordinates and assigned the trip to the nearest garage belonging to the respective agency. To ensure operational realism, assignments were balanced to maintain approximately equal trip counts across facilities within each agency. This process resulted in an average assignment of 2,650 trips to each CTA garage and 730 trips to each Pace garage. The planning horizon ($\overline{T}$) is set to 24 hours. For the exact model formulation, time is discretized into 5-minute intervals ($T^\Delta$). To ensure operational feasibility, we impose specific compatibility constraints for trip pairing: the maximum allowable layover time ($L$) is limited to 30 minutes, the maximum acceptable gap between trips ($G$) is set to 6 hours, and the minimum required dwell time at a station ($W$) is fixed at 30 minutes.

The charging network topology is derived from prior strategic planning studies \citep{bazarnovi2024problem}. \ref{app:charger_network_layout} illustrates the spatial distribution of charging facilities, distinguishing between garage-based chargers and on-route (opportunity) chargers. The number of plugs ($Z_c$) and the power output ($R_c$) at each location are fixed inputs for the operational problem, and all chargers in the baseline are assumed to be fast chargers. BEB parameters are standardized to reflect current market technology. We assume a battery capacity of 440 kWh with a usable SoC window between 20\% ($\underline{B}$) and 80\% ($\overline{B}$) to mitigate degradation and account for weather-related variances. Buses begin the operational day with an initial charge ($B^\iota$) of 72\%.

In the mathematical model, energy parameters are normalized by time. The energy consumption for a trip $i$ ($B_i$) is calculated based on trip length, a consumption rate of 2.8 kWh/mi, and the average network speed. Similarly, the charging rate ($R_c$) is converted into a time-based replenishment rate (minutes of operation gained per minute of charging). This rate varies by charger type, with slow chargers (125 kW) providing a (unitless) rate of $R^\text{slow}_c \approx 2.23$ min/min, and fast chargers (450 kW) delivering a rate of $R^\text{fast}_c \approx 8.03$ min/min. In practical terms, this indicates that a BEB recovers approximately 8.03 minutes of operational range for every minute of fast charging. \cref{param_value_model} summarizes the key scalar parameters used in the baseline model.

\begin{table}[!ht]
  \scriptsize
  \centering
  \caption{Parameter values used in the baseline scenario.}\label[tab]{param_value_model}
  \begin{tabularx}{1.0\textwidth}{lX|lX|lX}
  \toprule
    \textbf{Parameter} & \textbf{Value} & \textbf{Parameter} & \textbf{Value} & \textbf{Parameter} & \textbf{Value}\\
    \midrule
    $A^\nu$ & 1 (i.e., 100\%) & $F^\kappa$ & 52 \$/hr & $T^\Delta$ & 5 minutes\\
    $B^\iota$ & 5.66 hours & $G$ & 6 hours & $\overline{T}$ & 24 hours\\
    $\overline{B}$ & 6.29 hours & $L$ & 30 minutes & $V^\varepsilon$ & 386 \$/day\\
    $\underline{B}$ & 1.59 hours & $R^\text{fast}_c$ & 8.03 min/min & $V^\kappa$ & 251 \$/day\\
    $F^\varepsilon$ & 40 \$/hr & $R^\text{slow}_c$ & 2.23 min/min & $W$ & 30 minutes\\
    \bottomrule
  \end{tabularx}
\end{table}

To accurately compare BEBs against DBs, we normalize all capital and operating expenses into daily costs. Input data for these calculations are sourced from industry reports and technical literature, as detailed in \cref{param_values1}.

\begin{table}[!ht]
\centering
  \scriptsize
  \caption{Input data used to compute financial parameters.}\label[tab]{param_values1}
    \begin{tabular}{lll lll}
    \toprule
        \textbf{Definition} & \textbf{Value} & \textbf{Source} & \textbf{Definition} & \textbf{Value} & \textbf{Source} \\
        \midrule
        Battery capacity & 440 (kWh) & \citet{authority2022charging} & DB maintenance cost & 1 to 2.5 (\$/mi) & \citet{authority2022charging} \\
        Battery range & 20\% - 80\% & \citet{HU2022103732} & DB purchase cost & 650,000 (\$) & \citet{authority2022charging} \\
        BEB energy consumption & 2.8 (kWh/mi) & \citet{HE2023128227} & Diesel cost & 3.7 (\$/gal) & \citet{eiadieselprice} \\
        BEB maintenance cost & 1 to 2.5 (\$/mi) & \citet{authority2022charging} & Electricity cost & 0.08 (\$/kWh) & \citet{authority2022charging} \\
        BEB purchase cost & 1,000,000 (\$) & \citet{authority2022charging} & Fast charger power & 450 (kW) & \citet{authority2022charging} \\
        Bus avg. speed & 20 (mph) & Assumed & Inflation rate & 4.3 \% & \citet{blsinflation} \\
        Bus lifetime & 14 (years) & \citet{authority2022charging} & Slow charger power & 125 (kW) & \citet{authority2022charging} \\
        DB fuel consumption & 3.59 (mpg) & \citet{authority2022charging} & & & \\
        \bottomrule
    \end{tabular}
\end{table}

We employ the Capital Recovery Factor (CRF) to amortize the initial purchase price of vehicles over their 14-year lifespan ($n=14$) assuming an inflation rate ($r$) of 4.3\%:
\begin{equation}\label[eq]{CRF}
    CRF = \frac{r(1+r)^n}{(1+r)^n - 1}
\end{equation}
The resulting annual capital cost is divided by 250 operational days to derive the daily vehicle cost.

Variable costs, such as energy and fuel, are adjusted using the Present Value of Annuity (PVA) method to account for price escalation over the project horizon:
\begin{equation}\label[eq]{PVA}
    PVA = \frac{1 - (1+r)^{-n}}{r}
\end{equation}
Maintenance costs require a specific adjustment; based on CTA reports \citep{authority2022charging}, maintenance costs grow at a higher rate of 6.76\% annually. Consequently, we apply the PVA formula with $r = 6.76\%$ for maintenance components, while using the standard inflation rate for energy components.

\subsection{Computational experiments}\label[sec]{computational_experiments}
We evaluate the computational performance of two solution approaches: an \textit{exact} method and the CG heuristic. The exact method solves the MILP formulation using the Gurobi commercial solver \citep{gurobi}, employing standard branch-and-bound. The CG is compared against the exact approach to assess scalability, solution quality, and computational efficiency.

Instance generation utilizes garage and trip data from CTA and Pace. We define the problem size by the number of trips, $|\mathcal{I}| \in \{5, 10, 25, 100, 250, 500, 1000\}$. For each size, a garage is randomly chosen via uniform distribution, and a subset of trips is sampled. To ensure statistical power, the number of instances generated decreases as problem complexity increases: 1,000, 500, 250, 100, 50, 10, and 5 instances, respectively. Computation time limits are scaled to the problem size: 10, 30, 120, 600, 1,200, 2,400, and 3,600 seconds. To account for the stochastic nature of the CG, we solve each instance 10 times (referring to each run as a \textit{replica}). In total, 1,915 instances are solved using the exact method, and the CG is used 19,150 times.

\cref{quality_comparison} details the instance parameters and optimality results. The table reports the number of instances where the exact method and the CG method (best of 10 replicas) achieved optimality. A key metric, \textit{CG Impr.}, counts instances where the CG solution strictly dominated the solution found by the exact method when the latter failed to converge. The results delineate a clear operational boundary: the exact method is highly effective for small instances ($|\mathcal{I}| \leq 100$) but suffers from tractability issues as $|\mathcal{I}|$ grows. Conversely, while CG struggles to prove optimality in medium instances, it consistently outperforms the exact method in large-scale scenarios ($|\mathcal{I}| \geq 250$).

\begin{table}[!ht]
\centering
  \scriptsize
    \caption{Comparison of optimality counts and solution improvements.}
    \label[tab]{quality_comparison}
    \begin{center}
    \begin{tabular}{lrrrrr} 
    \toprule
    \multicolumn{3}{c}{Instance Parameters} & \multicolumn{2}{c}{\# Optimal} & \multicolumn{1}{c}{CG Impr.} \\
    \cmidrule(lr){1-3} \cmidrule(lr){4-5} \cmidrule(lr){6-6}
    $|\mathcal{I}|$ & Time Limit (s) & \# Inst. & Exact & CG & Count (Total) \\
    \midrule
    5 & 10 & 1,000 & 1,000 & 703 & - \\
    10 & 30 & 500 & 500 & 61 & - \\
    25 & 120 & 250 & 250 & 0 & - \\
    100 & 600 & 100 & 56 & 0 & 0/44 \\
    250 & 1,200 & 50 & 0 & - & 50/50 \\
    500 & 2,400 & 10 & 0 & - & 10/10 \\
    1,000 & 3,600 & 5 & 0 & - & 5/5 \\
    \bottomrule
    \end{tabular}
    \end{center}
    \scriptsize
    \vspace{-10pt}
    \emph{Note:} CG Impr. denotes the count of instances where the best CG solution was strictly better than the exact solution, calculated only for cases where the exact method did not prove optimality.
\end{table}

\cref{gap_comparison} quantifies solution quality. For the exact method, we define the lower and upper bounds reported by the solver as $\text{lb}^{\text{exact}}$ and $\text{ub}^{\text{exact}}$, respectively. The optimality gap is calculated as $100(1-\text{lb}^{\text{exact}}/\text{ub}^{\text{exact}})$. For the CG method, denoting the best objective value across replicas as $\text{ub}^{\text{CG}}$, we report the gap relative to the exact lower bound: $100(1-\text{lb}^{\text{exact}}/\text{ub}^{\text{CG}})$. Crucially, for instances where the exact method failed to converge, we introduce the \textit{CG Better} metric: $100(1-\text{ub}^{\text{CG}}/\text{ub}^{\text{exact}})$. This metric reveals a trend of increasing returns: as problem difficulty increases, the relative value of the CG heuristic grows. Specifically, the CG solution improves upon the feasible exact solution by 60.1\% at $|\mathcal{I}| = 250$, rising to 76.0\% at $|\mathcal{I}| = 1000$. This confirms that for large-scale operations, the CG heuristic provides significantly more realistic schedules than a time-truncated exact approach.

\begin{table}[!ht]
\centering
\scriptsize
\caption{Comparison of solution gaps. Values reported as Mean (SD) and Range [Min, Max].}
\label[tab]{gap_comparison}
\begin{center}
\begin{tabular}{ccccccc}
\toprule
Scenario & \multicolumn{2}{c}{Exact Gap (\%)} & \multicolumn{2}{c}{CG Gap (\%)} & \multicolumn{2}{c}{CG Better (\%)} \\
\cmidrule(lr){1-1} \cmidrule(lr){2-3} \cmidrule(lr){4-5} \cmidrule(lr){6-7}
$|\mathcal{I}|$ & Mean (SD) & Range & Mean (SD) & Range & Mean (SD) & Range \\
\midrule
5 & 0.0 (0.0) & [0.0, 0.0] & 0.3 (2.0) & [0.0, 28.8] & - & - \\
10 & 0.0 (0.0) & [0.0, 0.0] & 4.1 (12.2) & [0.0, 82.9] & - & - \\
25 & 0.0 (0.0) & [0.0, 0.0] & 50.7 (12.0) & [0.9, 74.9] & - & - \\
100 & 1.2 (2.9) & [0.0, 21.3] & 63.1 (5.3) & [44.8, 73.6] & - & - \\
250 & 97.9 (3.5) & [74.1, 99.6] & 95.0 (6.8) & [50.5, 99.4] & 60.1 (5.9) & [39.0, 73.2] \\
500 & 98.4 (0.9) & [97.2, 99.4] & 95.2 (2.6) & [91.4, 98.2] & 66.8 (4.4) & [56.5, 71.9] \\
1,000 & 98.4 (0.7) & [97.3, 99.1] & 93.3 (2.7) & [89.1, 96.0] & 76.0 (1.5) & [74.2, 77.5] \\
\bottomrule
\end{tabular}
\end{center}
\end{table}

\begin{figure}[!ht]
    \centering
    \begin{subfigure}[t]{0.49\textwidth}
        \centering
        \includegraphics[width=\linewidth]{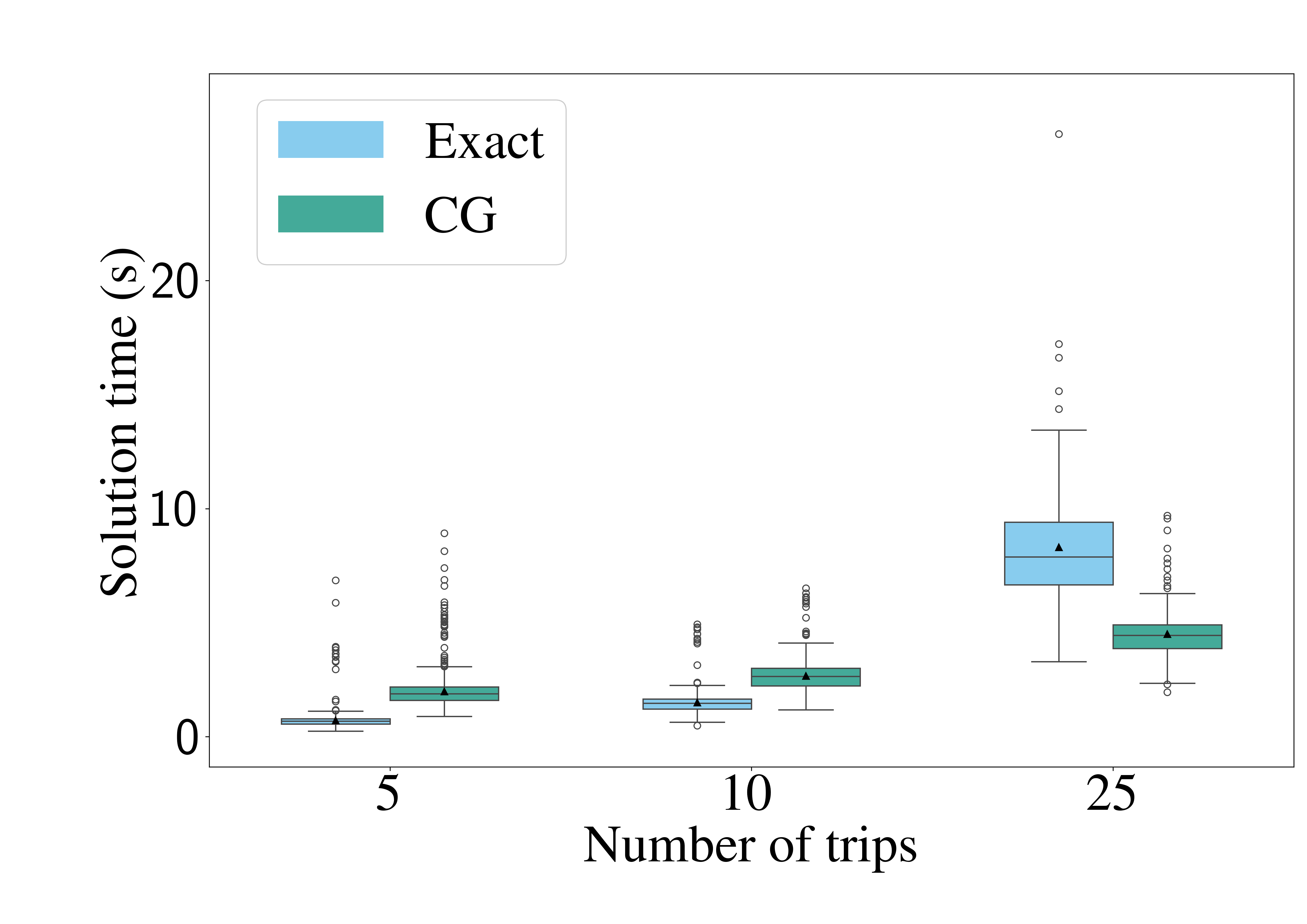}
        \caption{Scenarios with $|\mathcal{I}| \in \{5, 10, 25\}$}
        \label[fig]{comp_time_small}
    \end{subfigure}
    \hfill
    \begin{subfigure}[t]{0.49\textwidth}
        \centering
        \includegraphics[width=\linewidth]{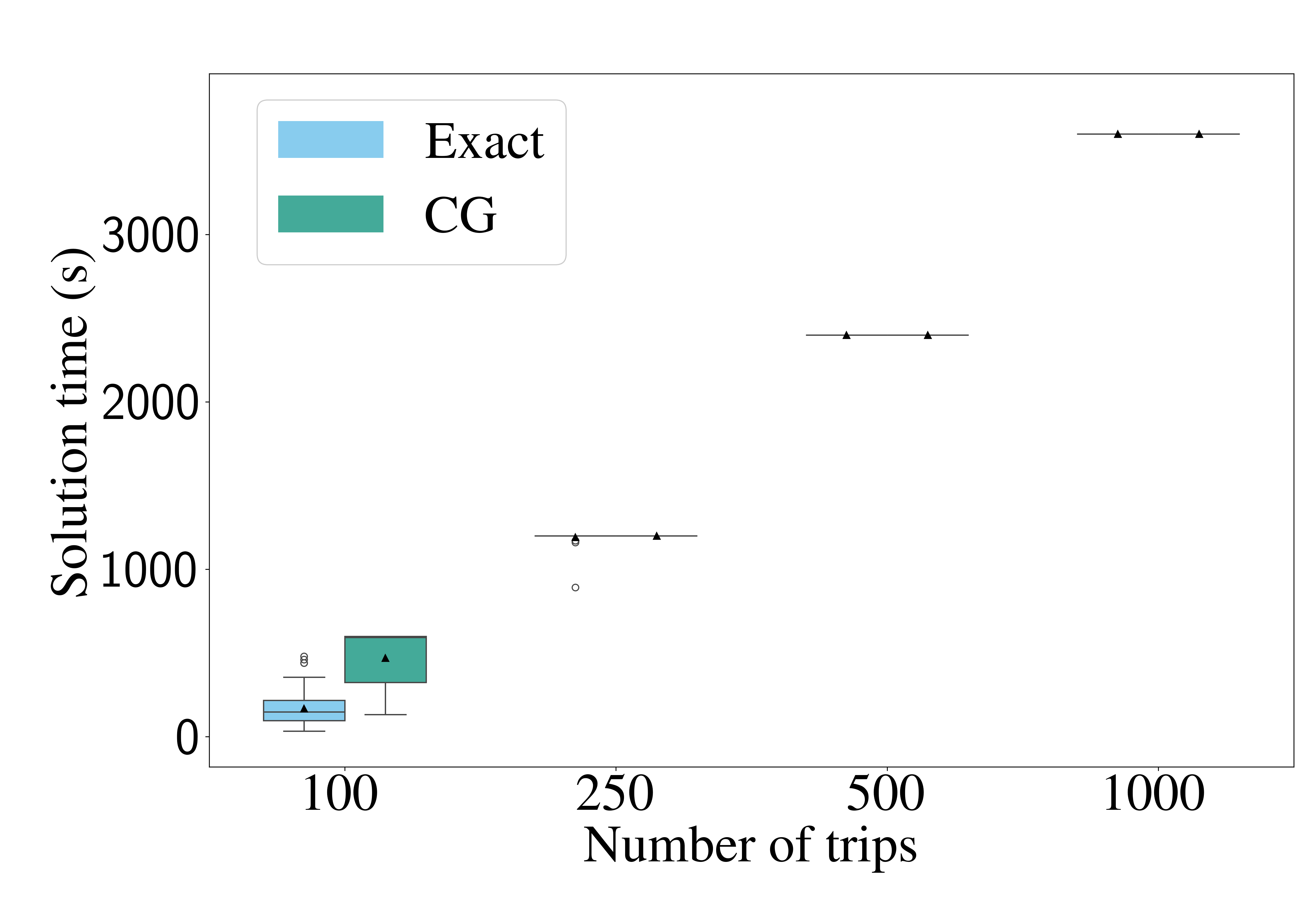}
        \caption{Scenarios with $|\mathcal{I}| \in \{100, 250, 500, 1000\}$}
        \label[fig]{comp_time_large}
    \end{subfigure}
    \caption{Comparison of computational time for exact and CG methods.}
    \label[fig]{comp_time}
\end{figure}

\cref{comp_time} compares computational times. \cref{comp_time_small} highlights a notable crossover point at $|\mathcal{I}| = 25$. While the exact method (blue) is faster for very small instances, it exhibits high variance and significant outliers at $|\mathcal{I}| = 25$. In contrast, the CG method (green) maintains a tighter interquartile range and a lower median solution time. This suggests that CG offers superior stability even before the problem size prevents the exact method from finding an optimal solution. For large instances (\cref{comp_time_large}), the distributions collapse to the time limit (represented by flat lines at 2,400s and 3,600s), indicating that both methods fully utilize the computational budget. \ref{app:more_computational} provides further analyses into CG stability and the impact of $T^\Delta$ on computational complexity. 

\subsection{Case study}\label[sec]{case_study}

In this section, we apply the proposed methodology to the complete transit network described in \cref{design_experiments}. Moving beyond the smaller-scale instances used for computational benchmarking, this analysis focuses on the full operational capabilities of the CTA and Pace bus networks. We evaluate the optimal charging schedules, fleet sizing requirements, and infrastructure utilization patterns that emerge when solving for the entire region's demand on a garage-by-garage basis.

We define eight scenarios based on the minimum required BEB fleet share, denoted by $A^\nu \in \{0, 1, 5, 10, 25, 50, 75, 100\}$ percent. The scenario where $A^\nu=0$ represents a baseline close to current operations, primarily driven by DBs, whereas $A^\nu=100$ represents a fully electric fleet. Through these scenarios, we analyze variations in total system cost, fleet composition, charging activities, and the temporal distribution of bus activities (e.g., revenue service, deadheading, and charging).

\begin{figure}[!ht]
    \centering
    \begin{subfigure}[t]{.49\textwidth}
        \centering
        \includegraphics[width=\linewidth]{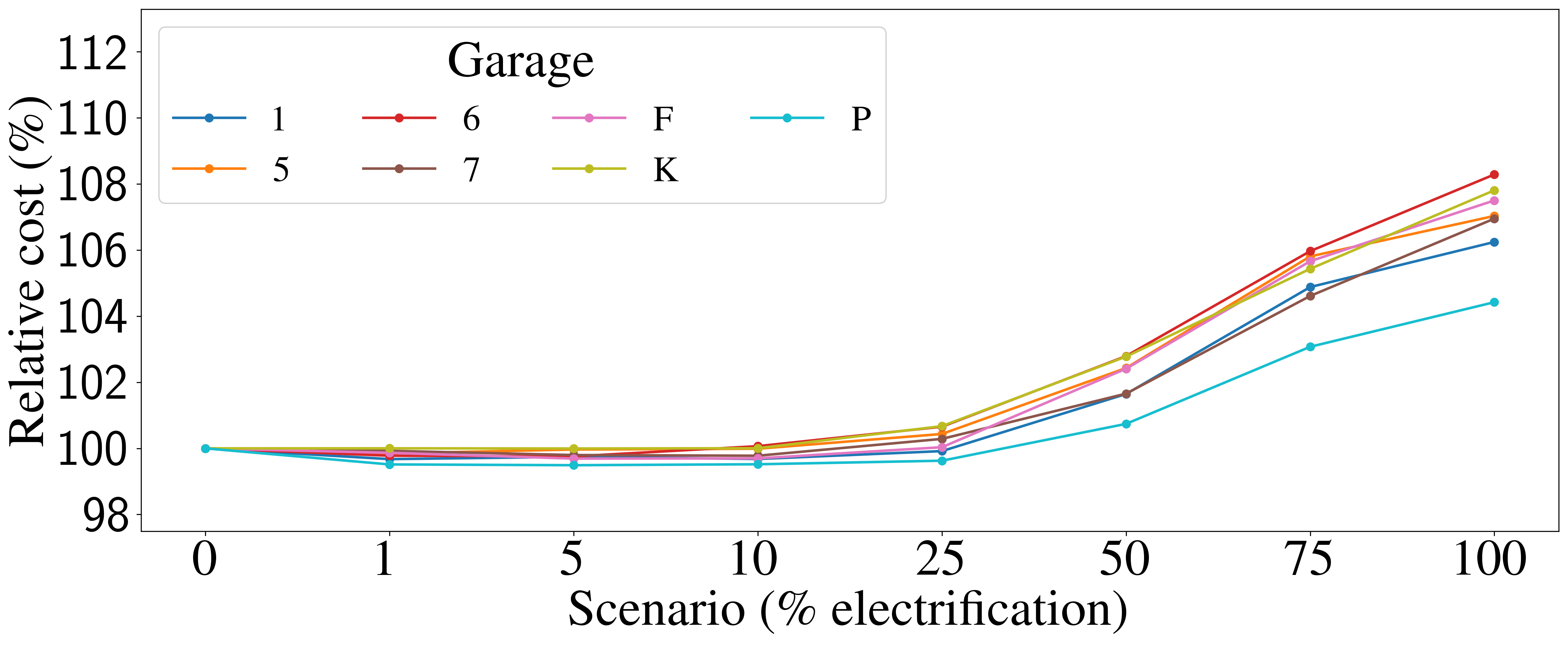} 
        \caption{CTA garages}
        \label[fig]{cost_CTA}
    \end{subfigure}
    \begin{subfigure}[t]{.49\textwidth}
        \centering
        \includegraphics[width=\linewidth]{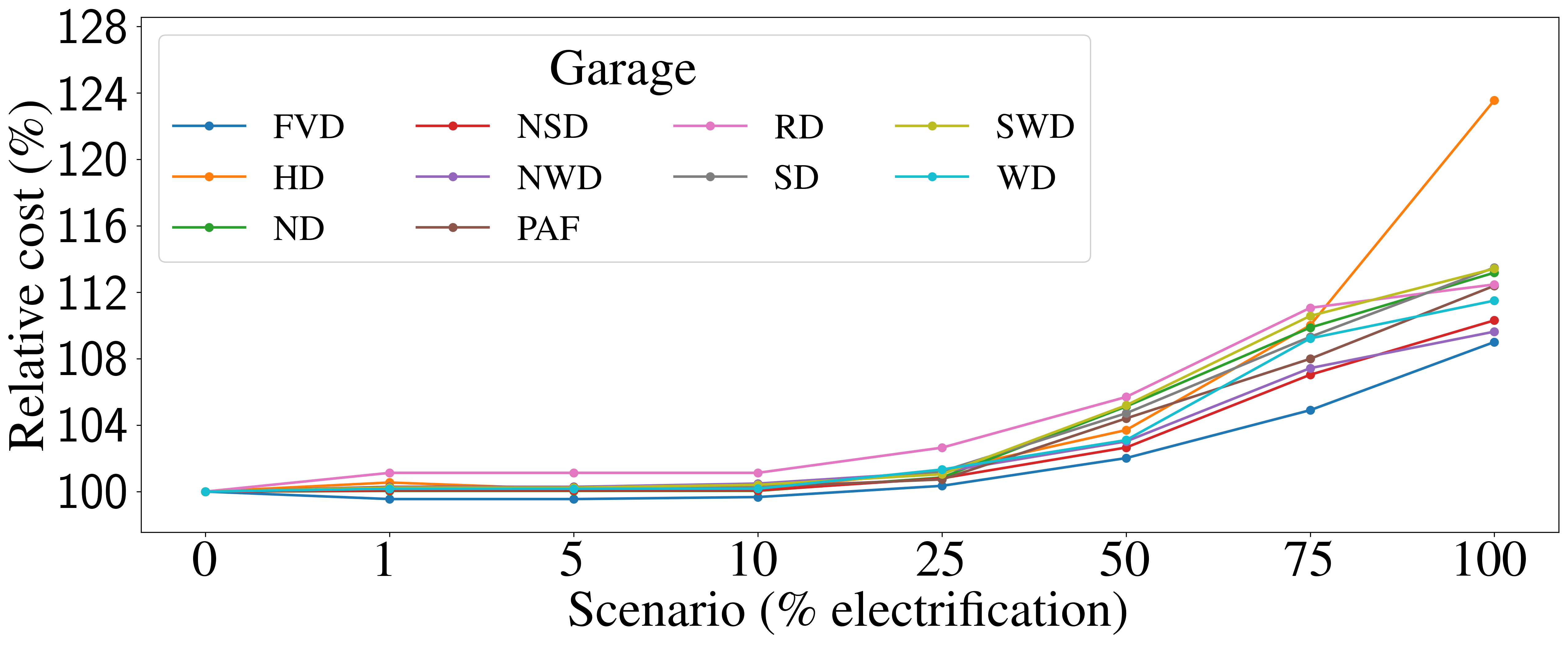} 
        \caption{Pace garages}
        \label[fig]{cost_Pace}
    \end{subfigure}
    \caption{Change in the total system cost relative to the baseline scenario ($A^\nu=0$).}
    \label[fig]{cost_change}
\end{figure}

\cref{cost_change} illustrates the percentage change in total system cost relative to the baseline scenario ($A^\nu=0$) for individual CTA and Pace garages. As the adoption of BEBs increases, a non-linear cost trajectory is observed. For many garages of CTA, the total system cost initially decreases (by up to 0.5\%) at lower adoption levels (e.g., 10\%). This suggests that the lower marginal operating costs of BEBs (specifically energy savings) outweigh their higher capital acquisition costs when they are selectively deployed on the most energy-intensive or operationally suitable routes. However, as $A^\nu$ exceeds 50\%, the total costs generally rise, culminating in an increase of approximately 3-24\% under the full transition scenario. This increase is driven by the significant capital expenditure required to replace the DB fleet with BEBs, which currently cost approximately \$1,000,000 compared to \$650,000 for a DB.

We observe distinct trends in \cref{cost_change} driven by the interplay between fleet size and operational savings. For instance, in CTA's Garage P (\cref{cost_CTA}), the total number of buses remains constant at 217 across scenarios 0 to 25, even as the BEB share increases. Because the total fleet size does not change, the overall costs decrease due to the savings in operational expenses, which outweigh the additional deadheading costs required for charging. However, as the BEB share increases to 50, 75, and 100\%, the total fleet size rises to 219, 223, and 222, respectively. At these higher penetration levels, the operational savings from BEBs are insufficient to offset the capital costs of the additional vehicles, causing the total system cost to increase.

A notable outlier is observed in Pace's Garage HD at the 100\% BEB share level (\cref{cost_Pace}), which exhibits a sharp cost spike. This anomaly is attributed to specific non-electrifiable trips within the garage's schedule that exceed the range capabilities of current BEBs. Consequently, the optimization model is forced to retain two DBs in the solution even in the full electrification scenario to ensure service feasibility.

\cref{agency_level_changes} presents the aggregate solution statistics for both transit agencies across all scenarios. First, comparing the baseline (Scenario 0) against real-world fleet sizes provides context for the model's calibration. CTA reportedly operates over 1,800 buses \citep{authority2022charging}, and Pace operates approximately 733 buses \citep{pace}. In contrast, our optimization model determines optimal fleet sizes of 1,463 and 679 for CTA and Pace, respectively. This difference is expected, as the model calculates the theoretical minimum fleet required to satisfy schedule and charging constraints, without accounting for crew scheduling restrictions, maintenance spares, or other agency-specific operational buffers. Consequently, these results serve as a lower-bound operational baseline.

\begin{table}[!ht]
    \centering
    \scriptsize
    \caption{Aggregate solution statistics across scenarios with cost change (\%) relative to Scenario 0.}\label[tab]{agency_level_changes}
    \begin{tabular}{llcccc}
        \toprule
        Agency & Scenario ($A^\nu$) & Cost Change (\%) & Fleet Size & \# BEBs & \# DBs \\
        \midrule
        \multirow{8}{*}{CTA} 
          & 0   & -     & 1,463 & 0     & 1,463 \\
          & 1   & -0.20 & 1,468 & 307   & 1,161 \\
          & 5   & -0.22 & 1,467 & 308   & 1,159 \\
          & 10  & -0.18 & 1,469 & 306   & 1,163 \\
          & 25  & 0.24  & 1,478 & 461   & 1,017 \\
          & 50  & 2.07  & 1,522 & 807   & 715 \\
          & 75  & 5.07  & 1,581 & 1,246 & 335 \\
          & 100 & 6.91  & 1,594 & 1,594 & 0 \\
        \cmidrule{1-6}
        \multirow{8}{*}{Pace} 
          & 0   & -     & 679 & 0   & 679 \\
          & 1   & 0.22  & 681 & 81  & 600 \\
          & 5   & 0.18  & 680 & 81  & 599 \\
          & 10  & 0.24  & 680 & 93  & 587 \\
          & 25  & 1.09  & 683 & 203 & 480 \\
          & 50  & 3.93  & 704 & 371 & 333 \\
          & 75  & 8.71  & 739 & 581 & 158 \\
          & 100 & 12.99 & 761 & 759 & 2 \\
        \bottomrule
    \end{tabular}
\end{table}

Regarding the electrification scenarios, we observe a general increase in both total system cost and fleet size as $A^\nu$ rises. This reflects the operational trade-offs required to accommodate the range and charging limitations of electric buses, often necessitating a larger fleet to maintain the same level of service. Compared to CTA's 6.91\%, Pace presents a higher cost increase at $A^\nu=100$\% with 12.99\%.

\begin{figure}[!ht]
    \centering
    \begin{subfigure}[t]{.49\textwidth}
        \centering
        \includegraphics[width=\linewidth]{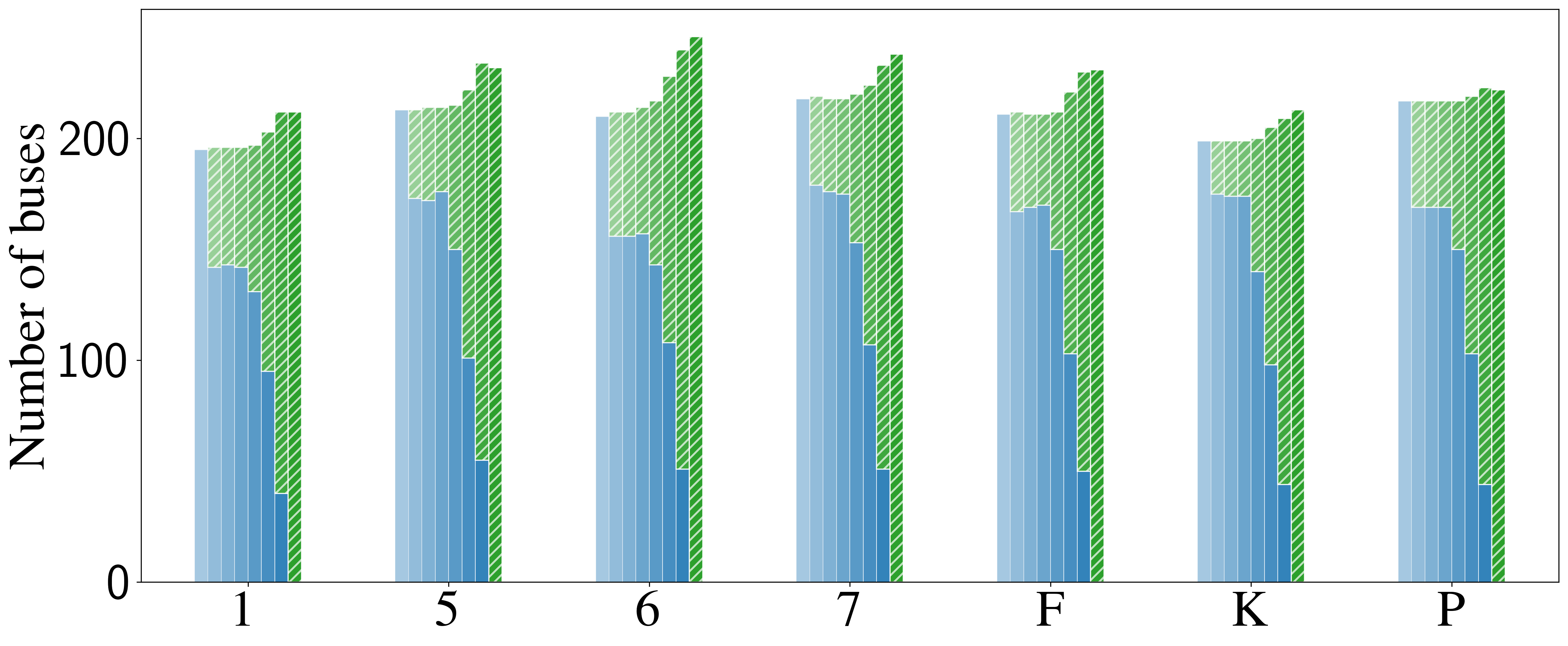} 
        \caption{CTA garages}
        \label[fig]{fleet_CTA}
    \end{subfigure}
    \hfill
    \begin{subfigure}[t]{.49\textwidth}
        \centering
        \includegraphics[width=\linewidth]{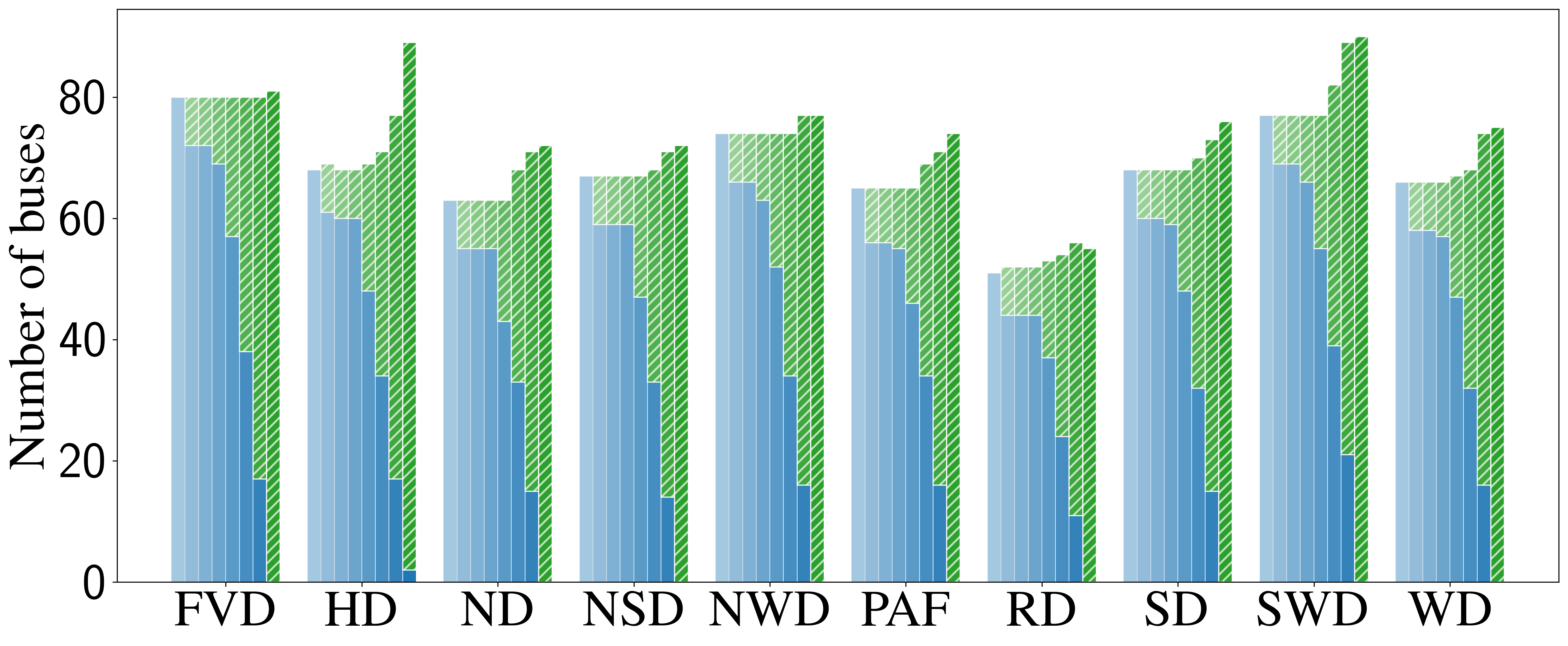} 
        \caption{Pace garages}
        \label[fig]{fleet_Pace}
    \end{subfigure}
    \begin{subfigure}[t]{\textwidth}
        \centering
        \includegraphics[width=\linewidth]{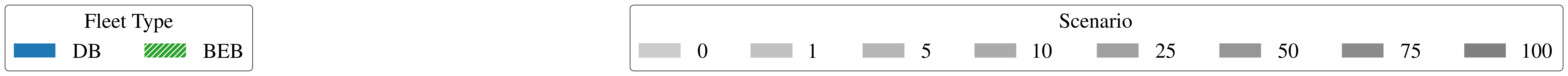}
    \end{subfigure}
    \caption{Fleet decomposition under different BEB penetration rates.}
    \label[fig]{fleet_decomposition}
\end{figure}

\cref{fleet_decomposition} details the optimal fleet size and decomposition for each garage. The results indicate distinct characteristics between the two agencies. For Pace garages (\cref{fleet_Pace}), the total fleet size displays sharp increases in high BEB penetration scenarios, particularly for Garages HD and SWD. Conversely, CTA garages (\cref{fleet_CTA}) generally depict a smoother, gradual increase in fleet size as $A^\nu$ rises. This trend suggests that replacing DBs with BEBs on CTA's typically shorter, dense urban routes requires fewer additional vehicles compared to Pace's longer suburban service trips, where range limitations are more acute.

We also observe some fluctuations in \cref{fleet_decomposition}, specifically in CTA's Garages 5 and P, and Pace's Garage RD. For example, moving from the $A^\nu=75$\% to the $A^\nu=100$\% scenario, there is a slight, counter-intuitive drop in the total fleet size for Garage P. Analyzing the results in detail reveals that this is due to the discrete nature of the optimization problem. To strictly satisfy the 75\% share constraint in Garage P, the model allocates 179 BEBs and 44 DBs. However, in the 100\% scenario, the model returns an optimal solution of 222 BEBs. This implies that in the 75\% case, an additional bus was generated specifically to satisfy the ratio constraint and avoid penalties, a requirement that relaxed naturally in the full electrification scenario.

\begin{figure}[!ht]
    \centering
    \includegraphics[width=.5\linewidth]{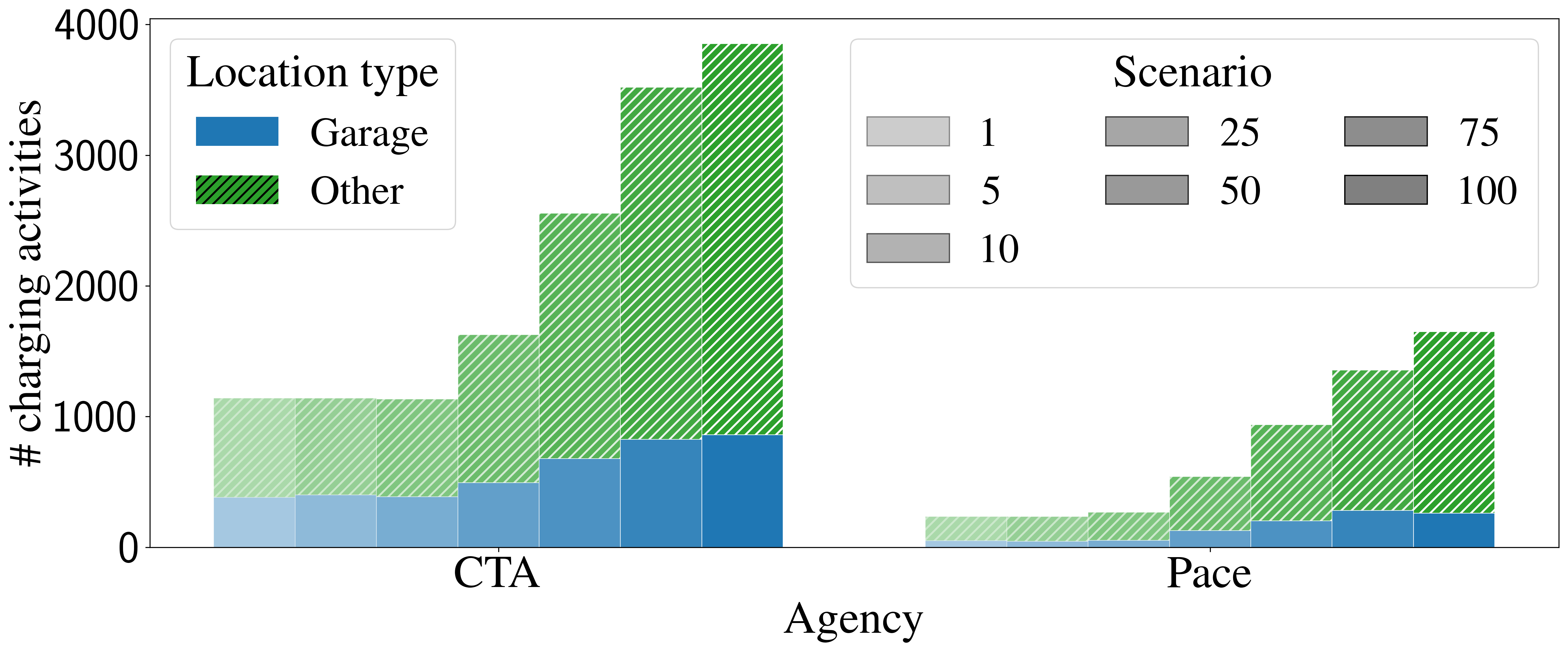} 
    \caption{Spatial distribution of charging activities (Garage vs. Other locations).}
    \label[fig]{charging_activity_locations}
\end{figure}

\cref{charging_activity_locations} presents the spatial distribution of charging events. A striking finding is the predominance of charging activities at non-garage locations (designated as \textit{Other}). Across all scenarios, on-route and terminal chargers handle the vast majority of energy replenishment. This indicates that the optimization model strongly favors opportunistic charging, utilizing dwell times at terminals between trips over the deadheading required to return to a garage for charging. This validates the importance of a flexible charging strategy that incorporates fast chargers at terminal nodes to maintain operational efficiency.

\begin{figure}[!ht]
    \centering
    \begin{subfigure}[t]{.49\textwidth}
        \centering
        \includegraphics[width=\linewidth]{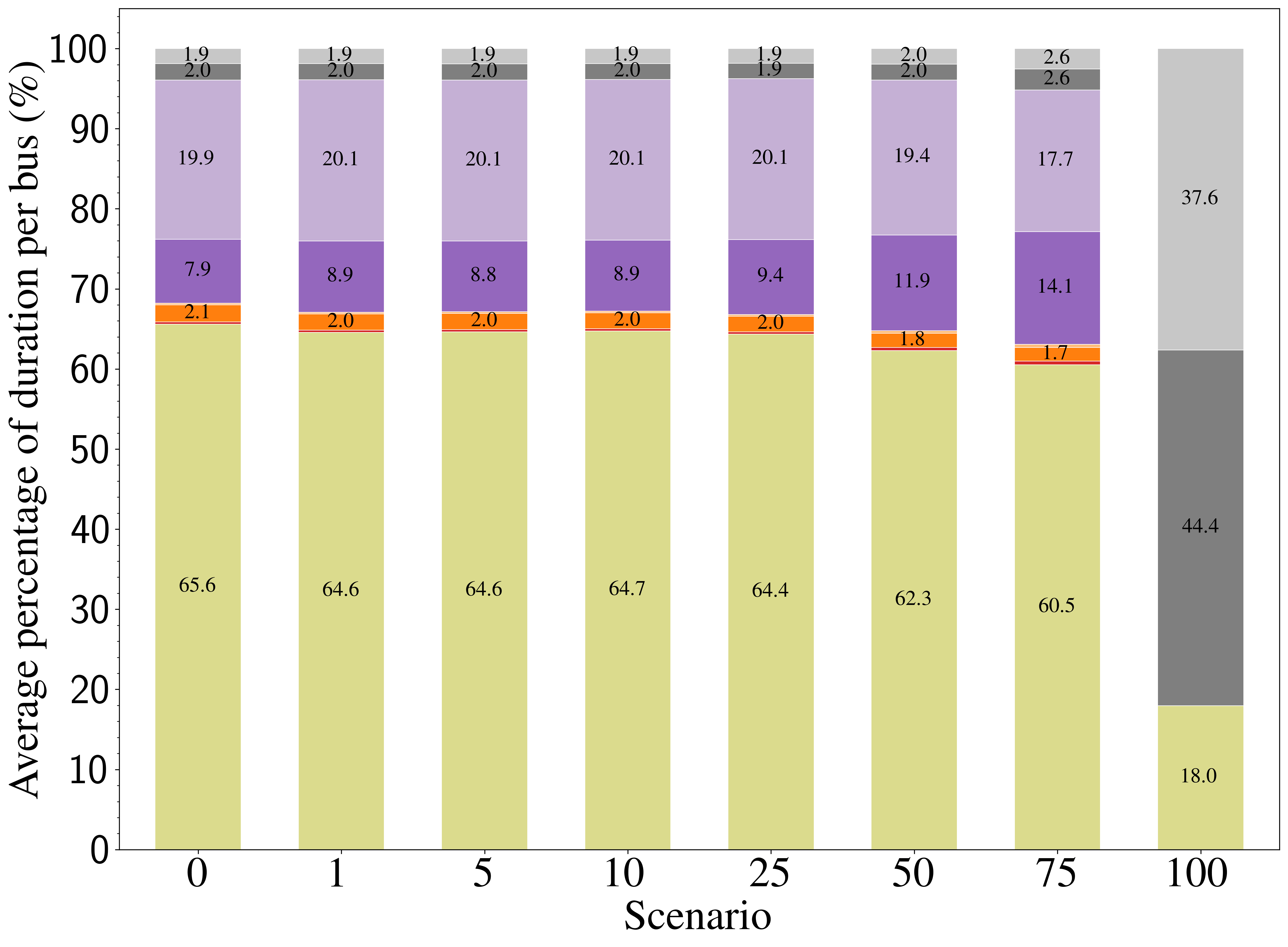} 
        \caption{DBs}
        \label[fig]{activity_DB}
    \end{subfigure}
    \hfill
    \begin{subfigure}[t]{.49\textwidth}
        \centering
        \includegraphics[width=\linewidth]{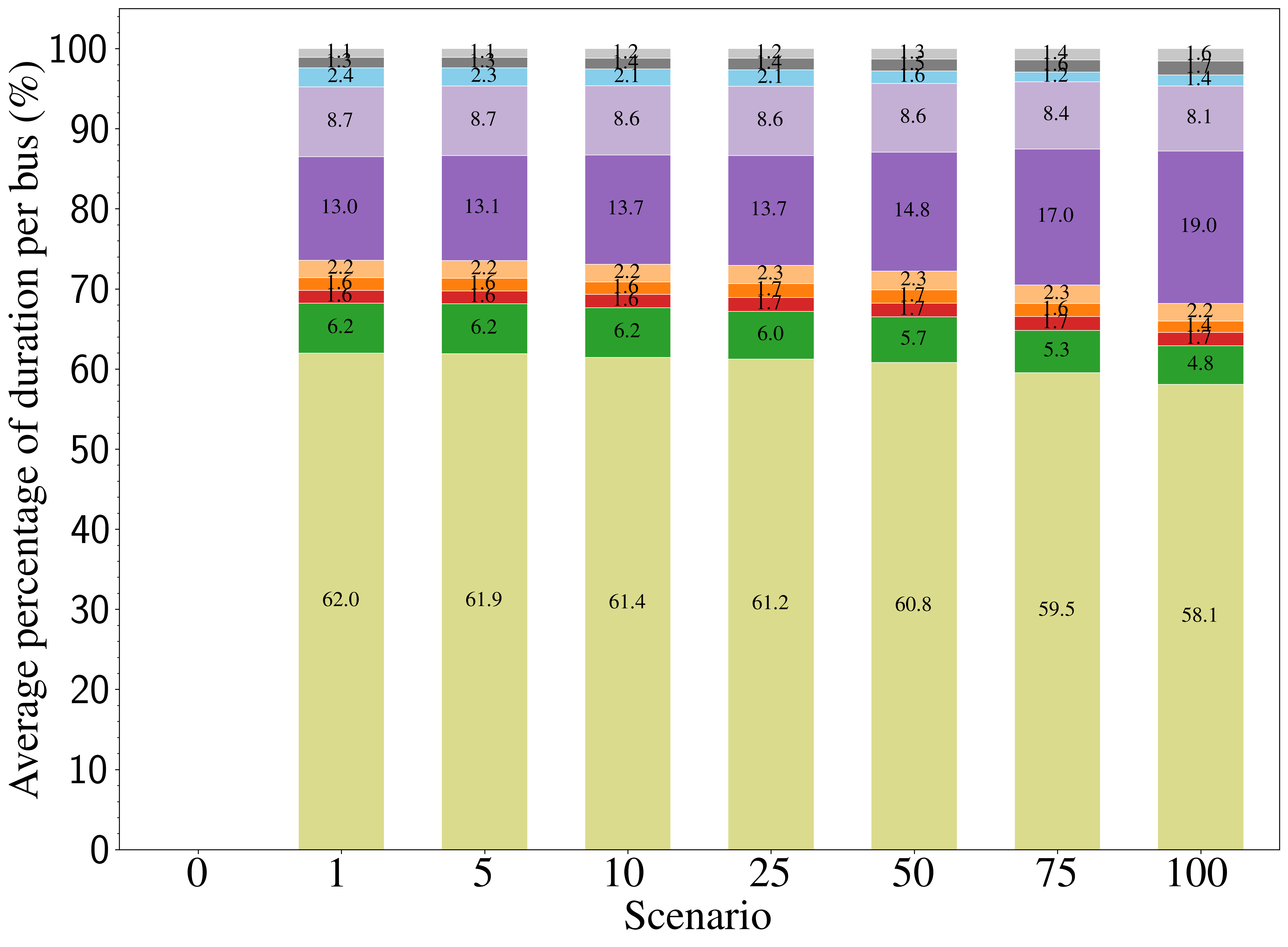} 
        \caption{BEBs}
        \label[fig]{activity_BEB}
    \end{subfigure}
    
    \begin{subfigure}[t]{\textwidth}
        \centering        \includegraphics[width=\textwidth]{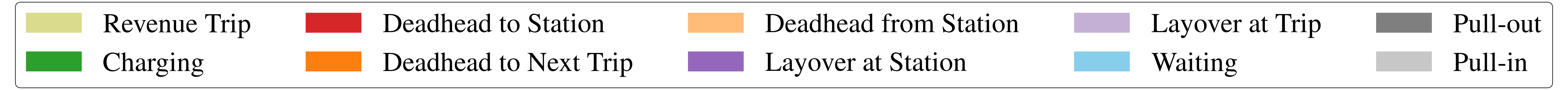} 
    \end{subfigure}
    
    \caption{Average percentage of operational time distributed across different activities.}
    \label[fig]{activity_percent_decomposition}
\end{figure}

As the fleet transitions to alternative powertrains, the allocation of vehicle time shifts. \cref{activity_percent_decomposition} breaks down the average duration spent on various activities, displaying these metrics for DBs and BEBs, respectively. For BEBs, the operational profile remains relatively stable as the penetration rate $A^\nu$ increases. The percentage of time spent on revenue trips hovers around 58-62\% indicating that the schedule optimization efficiently integrates electric buses without significantly eroding service productivity. This revenue share is maintained by balancing necessary station layovers and charging times. Notably, the time spent actually charging remains low (approximately 5--6\%, and the time spent \textit{waiting} for a charger is negligible (below 2.5\%). This validates the efficacy of the proposed adaptive queuing strategy, which successfully schedules charging events to avoid congestion. Conversely, the data for DBs in the 100\% scenario (\cref{activity_DB}) appears skewed; this is because the few remaining DBs are assigned to specific, isolated long-haul trips that lack the typical chain of connections, resulting in disproportionately high pull-in/pull-out times relative to their service time.

\subsection{Sensitivity analyses}\label[sec]{sensitivity}
To assess the robustness of the optimal solutions and understand the drivers of cost and efficiency, this section performs sensitivity analyses on critical model parameters. By systematically varying factors such as battery range limits, charger power output, and varying cost structures, we quantify their impact on the total objective function, vehicle and operational cost components, and the resulting fleet composition.

We consider three problem sizes for the analyses: $|\sI|\in\{25,50,100\}$. Following a similar process as in \cref{computational_experiments}, we randomly draw a sample from a randomly chosen garage. We generate 100 problem instances for each problem size. We solve each instance considering eight scenarios. For each instance, we provide a time limit of 120, 300, and 600 seconds for problem sizes of $|\sI|=25$, $|\sI|=50$, and $|\sI|=100$, respectively.

The scenarios and their specific lever designs are motivated by evolving technologies, market uncertainties, and operational constraints. In the sensitivity figures presented later, the x-axis \textit{Lever} indices correspond to the sequential non-baseline alternatives listed in each scenario, and the baseline in each scenario is Lever 0. Scenarios are defined as follows:

\begin{enumerate}[leftmargin=*, noitemsep, topsep=0pt]
    \item \textit{Battery capacity:} Recognizing that battery technology is evolving and higher capacities are expected in the future, we multiply the default value (440 kWh as noted in \cref{param_values1}) by factors of 1, 1.2, 1.4, 1.6, 1.8, 2, 2.2, 2.4, and 2.6. Factor of 1 is considered to be the baseline (Lever 0) in this scenario.
    
    \item \textit{Battery range:} While the default configurations assume a range of 20\%--80\% (refer to \cref{param_values1}), real-world operations may allow for condensed or stretched windows. We test nine specific windows: 10\%--90\%, 20\%--90\%, 30\%--90\%, 10\%--80\%, 20\%--80\%, 30\%--80\%, 10\%--70\%, 20\%--70\%, and 30\%--70\%. Here, we treat 10\%--90\% as the baseline scenario. 
    
    \item \textit{Charger layout:} To assess the impact of infrastructure density and type, we compare the existing all-fast-charger (Lever 0) default configuration against seven alternatives: (1) 25\% reduced charger count, (2) 50\% reduced charger count, (3) 75\% reduced charger count, (4) all slow chargers (same count as the default configuration), (5) slow chargers at garages with fast chargers at other locations, (6) fast chargers only at garages (no chargers elsewhere), and (7) slow chargers only at garages (no chargers elsewhere). For the reduced count levers, we rounded up the number of chargers at locations while reducing the numbers based on their corresponding percentages.
    
    \item \textit{Charger power:} Acknowledging that high-power chargers are expected to be used in the future, we multiply the baseline fast charger power (450 kW) by 1 (Lever 0), 1.25, 1.5, 1.75, 2, 2.25, 2.5, 2.75, and 3.

    \item \textit{Diesel cost:} Fuel prices fluctuate significantly, making future forecasting difficult. We test changes to the referenced price of 3.7 \$/gal (in \cref{param_values1}) by multiplication factors of 0, 0.25, 0.50, 0.75, 1, 1.25, 1.50, 1.75, and 2. This range encompasses scenarios from free fuel to double the current market rates. These instances use the default setting of $A^\nu=100\%$. Here, Lever 0 is considered to be the one with no diesel cost.
    
    \item \textit{Diesel cost ($A^\nu=0$):} Since changes in diesel cost should not considerably impact the total system cost in fully electric scenarios ($A^\nu=100\%$), we run a separate set of scenarios setting $A^\nu=0$ to properly observe the sensitivity. The levers remain the same as the standard diesel cost scenario.
    
    \item \textit{Electricity cost:} Similar to diesel, we apply the same multiplication factors to the referenced electricity cost of 0.08 \$/kWh (in \cref{param_values1}). As in its diesel counterpart, Lever 0 in this scenario is no electricity cost.
    
    \item \textit{Vehicle cost:} As BEB technology emerges, costs are expected to decrease and potentially converge with DBs. We create levers by applying multipliers to the referenced (in \cref{param_values1}) costs of BEBs and DBs respectively: (1, 1), (0.9, 0.99), (0.85, 0.975), (0.8, 0.95), (0.75, 0.925), (0.7, 0.9), (0.65, 0.85), (0.6, 0.8), and (0.5, 0.77). The first lever is treated as the baseline, and the last lever considers equal vehicle costs.
\end{enumerate}

\begin{figure}[!htbp]
    \centering
    \begin{subfigure}[t]{.9\textwidth}
        \centering
        \includegraphics[width=\linewidth]{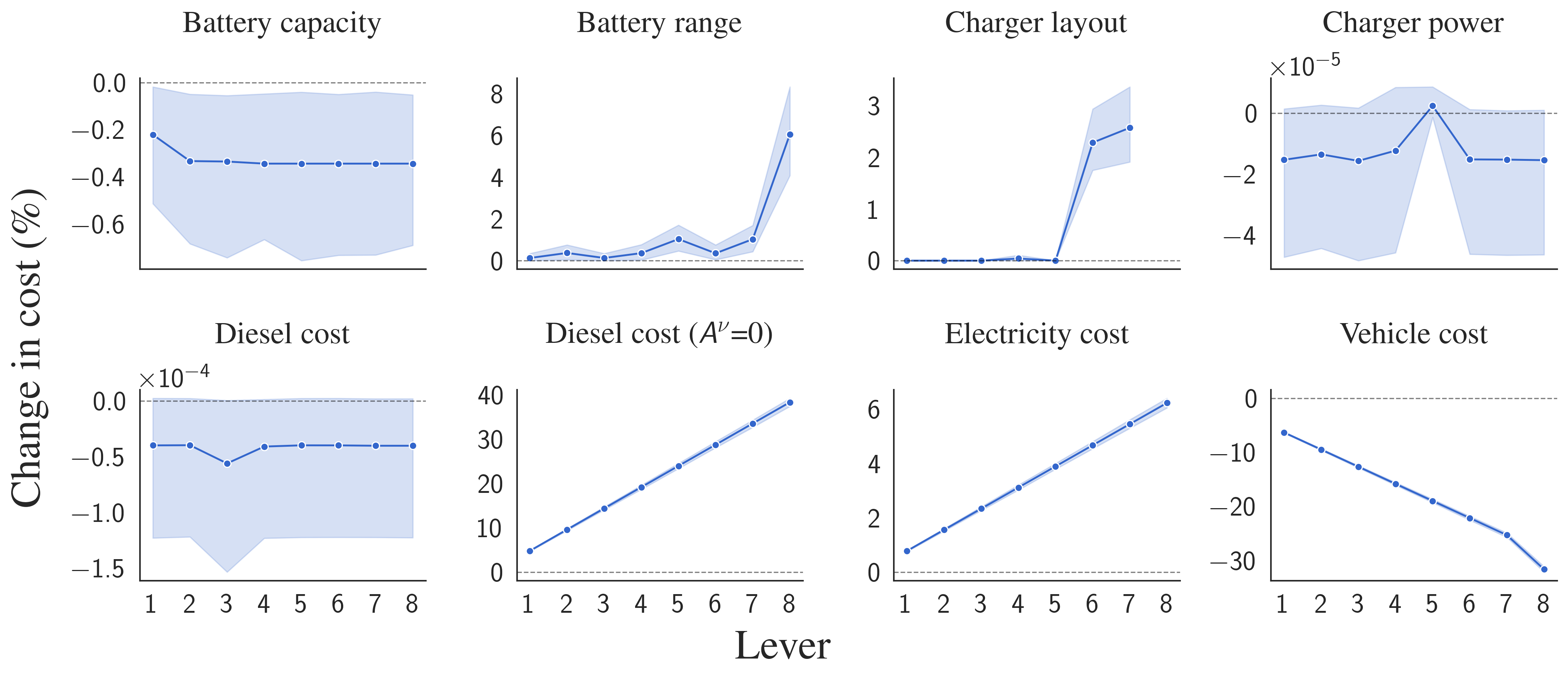} 
        \caption{$|\sI|=25$}
        \label[fig]{sens_cost_25}
    \end{subfigure}
    \begin{subfigure}[t]{.9\textwidth}
        \centering
        \includegraphics[width=\linewidth]{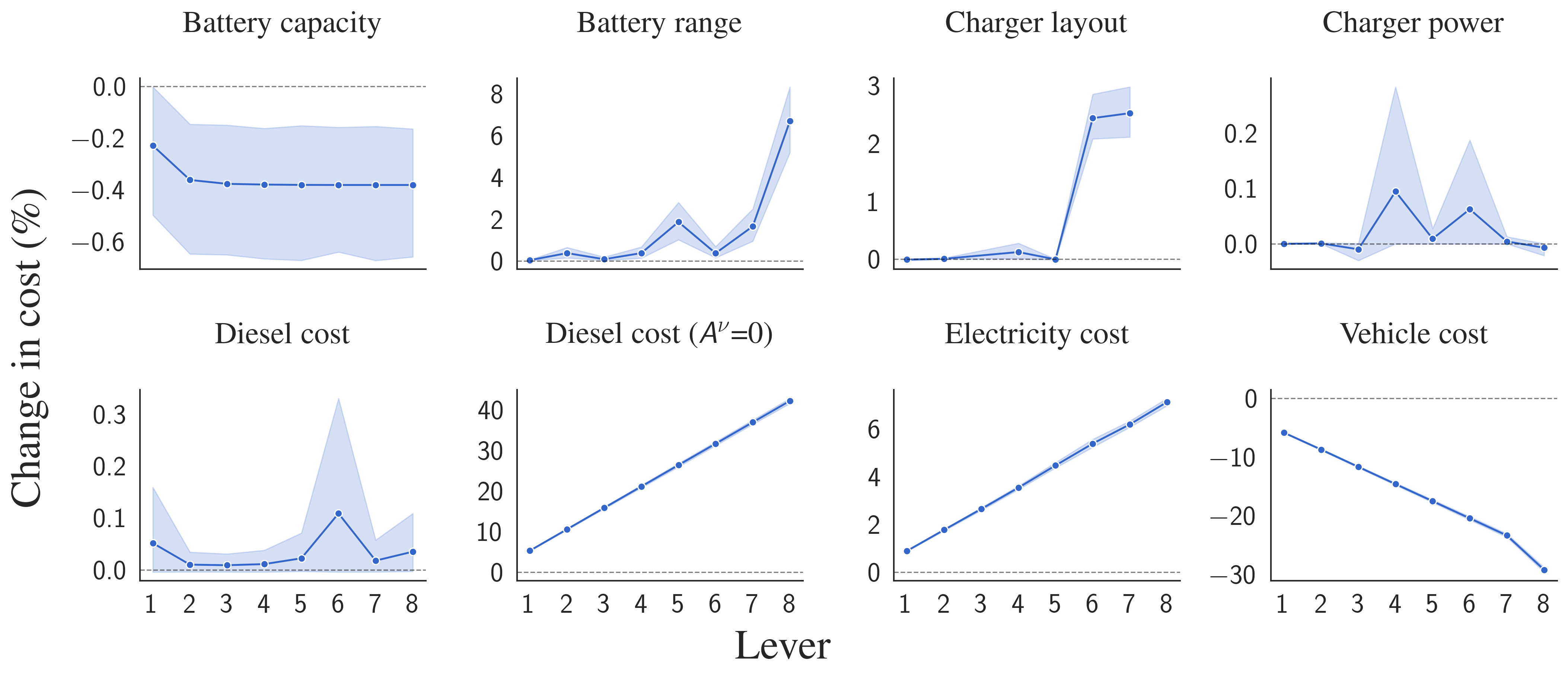} 
        \caption{$|\sI|=50$}
        \label[fig]{sens_cost_50}
    \end{subfigure}
    \begin{subfigure}[t]{.9\textwidth}
        \centering
        \includegraphics[width=\linewidth]{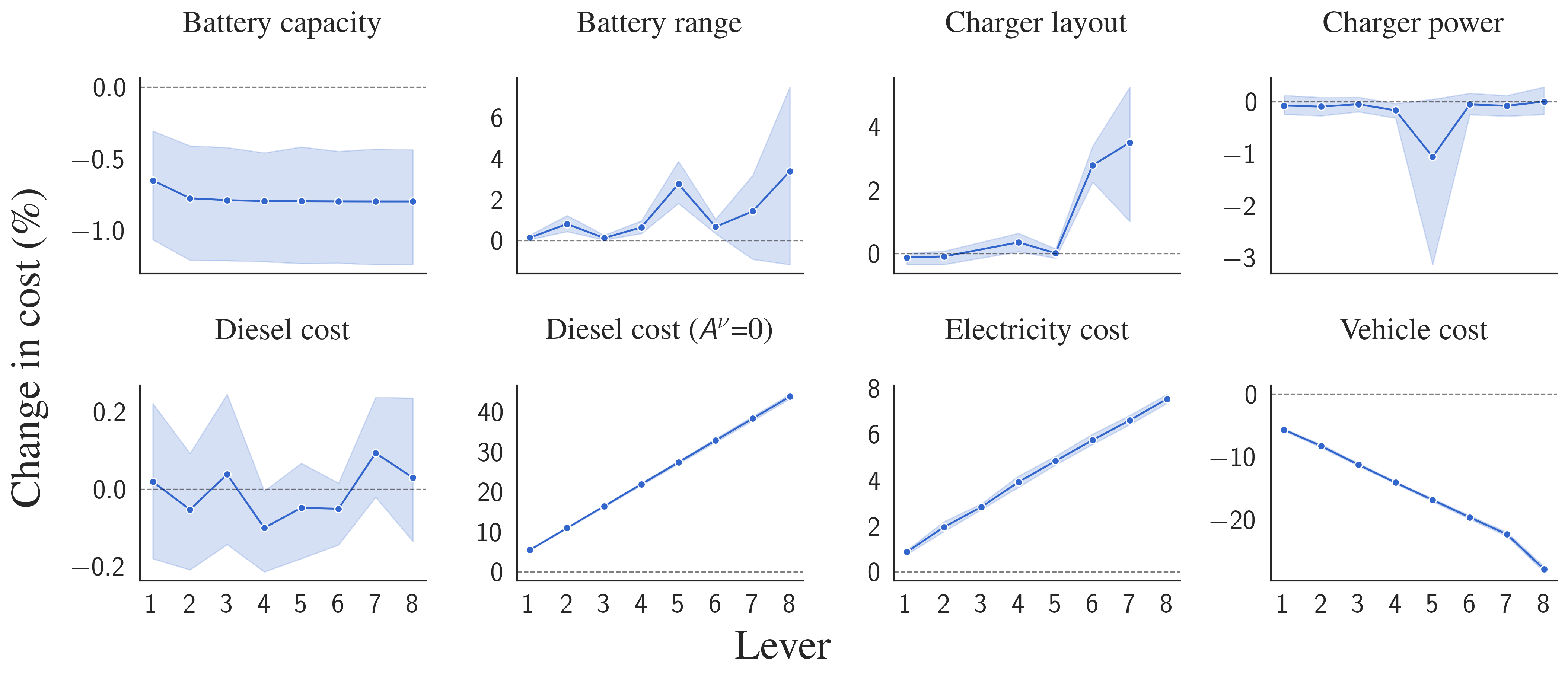} 
        \caption{$|\sI|=100$}
        \label[fig]{sens_cost_100}
    \end{subfigure}
    \caption{Percentage change in total cost of scenarios with respect to their baseline levers.}
    \label[fig]{sensitivity_cost}
\end{figure}

\cref{sensitivity_cost} shows the average percentage change (with 95\% confidence interval via the shaded region) in total cost per scenario and lever with respect to its specific baseline for three different problem sizes. Each point on these graphs correspond to average percent change value of 100 instances. We observe that a drop in vehicle costs and an increase in diesel costs (in the $A^\nu=0$ case) considerably impact the total system cost. Battery range and electricity cost also have a considerable impact. In \cref{sens_cost_25}, we notice that charger power does not have a significant impact on total system cost. Slight changes in total cost in \cref{sens_cost_50} and \cref{sens_cost_100} for this scenario are mainly due to sub-optimality in the solutions. The standard diesel cost scenario (where $A^\nu=100\%$) demonstrates minimal random spikes; this is expected because these instances do not utilize DBs in many cases, rendering the stochastic noise inherent to the solution method.

Regarding the charger layout scenario, the \textit{only garage} settings result in a 3--5\% increase in total cost, whereas simply reducing the number of chargers or swapping to slow chargers at garages does not show a high impact. This is because the charging infrastructure was optimized, and the majority of charging events take place at other locations (terminals). Similar to total cost, we also present the average impact of these scenarios on the number of BEBs, as well as vehicle and operational cost components, in \ref{app:sensitivity_plots}.

\section{Conclusion}\label[sec]{conclusion}
This study presented an integrated optimization framework for simultaneously determining optimal scheduling and charging strategies for mixed fleets consisting of BEBs and DBs. The developed MILP model uniquely combines decisions regarding optimal fleet composition, scheduling, and flexible charging strategies, effectively addressing both slow and fast charger types located at terminal stops and garages. By allowing partial recharging and implementing a novel adaptive queueing strategy, the model significantly enhanced scheduling flexibility, reduced operational delays, and optimized charging infrastructure utilization.

Numerical experiments using a real-world transit network in the Chicago area confirmed the efficacy of this approach. The results demonstrated that BEBs, despite higher upfront acquisition costs, serve a disproportionately higher share of trips due to lower energy expenditures. Consequently, while full BEB fleet ($A^\nu=100\%$) results in a modest increase in total system costs compared to a baseline DB fleet, mixed-fleet scenarios (e.g., $A^\nu=25\%$) were shown to reduce total costs, indicating an immediate pathway for a partial BEB fleet. 

Methodologically, we developed and presented a detailed CG framework to solve this complex problem. The methodology decomposes the problem into an RMP for schedule selection and two distinct pricing subproblems to generate candidate columns for both DB and BEB fleets. A key contribution of this work is the formulation of the BEB subproblem as an SPPRC, which models the intricate operational rules of BEBs, including SoC management and charging schedules. Furthermore, we addressed practical challenges inherent to the CG procedure, particularly the issue of algorithmic stalling caused by degeneracy in the RMP. 


Overall, the developed optimization approach provides valuable insights for transit agencies transitioning to alternative powertrain fleets. It highlights the strategic trade-offs between initial fleet investment and ongoing operational savings, demonstrating that a carefully managed transition to a mixed fleet can be cost-competitive. The integration of flexible charging strategies and adaptive scheduling methods presented in this work lays a strong foundation for future studies and practical implementations, paving the way for more efficient, resilient, and cost-effective public transit systems.

\section*{Acknowledgments}
This material is based on work supported by the U.S. Department of Energy, Office of Science, under contract number DE-AC02-06CH11357. This report and the work described were sponsored by the U.S. Department of Energy (DOE) Vehicle Technologies Office (VTO) under the Pathways to Affordable, Convenient, and Efficient Regional Mobility, an initiative of the Energy Efficient Mobility Systems (EEMS) Program. Erin Boyd, a DOE Office of Critical Minerals and Energy Innovation (CMEI) manager, played an important role in establishing the project concept, advancing implementation, and providing guidance.

\bibliographystyle{elsarticle-harv}
\bibliography{cas-refs}

@misc{gurobi,
  author = {{Gurobi Optimization, LLC}},
  title = {{Gurobi Optimizer Reference Manual}},
  year = 2024,
  url = "https://www.gurobi.com"
}

@misc{pace,
  author = {{Pace Suburban Bus}},
  title = {About Pace},
  year = {2025},
  url = {"https://www.pacebus.com/about"}
}

@misc{GTFS,
author = {{GTFS}},
title = {General Transit Feed Specification Reference},
howpublished = {Available at \url {https://gtfs.org/schedule/reference/} accessed on May 17, 2024},
year = {2022}
}

@article{ALVO2021102528,
title = {An exact solution approach for an electric bus dispatch problem},
journal = {Transportation Research Part E: Logistics and Transportation Review},
volume = {156},
pages = {102528},
year = {2021},
issn = {1366-5545},
doi = {10.1016/j.tre.2021.102528},
author = {Matías Alvo and Gustavo Angulo and Mathias A. Klapp},
}

@INPROCEEDINGS{HU20219564617,
  author={Hu, Hao and Du, Bo and Perez, Pascal},
  booktitle={2021 IEEE International Intelligent Transportation Systems Conference (ITSC)}, 
  title={Integrated optimisation of electric bus scheduling and top-up charging at bus stops with fast chargers}, 
  year={2021},
  volume={},
  number={},
  pages={2324-2329},
  doi={10.1109/ITSC48978.2021.9564617}
}

@article{HU2022103732,
title = {A joint optimisation model for charger locating and electric bus charging scheduling considering opportunity fast charging and uncertainties},
journal = {Transportation Research Part C: Emerging Technologies},
volume = {141},
pages = {103732},
year = {2022},
issn = {0968-090X},
doi = {10.1016/j.trc.2022.103732},
author = {Hao Hu and Bo Du and Wei Liu and Pascal Perez},
}

@article{HE2023103653,
title = {Joint optimization of electric bus charging infrastructure, vehicle scheduling, and charging management},
journal = {Transportation Research Part D: Transport and Environment},
volume = {117},
pages = {103653},
year = {2023},
issn = {1361-9209},
doi = {10.1016/j.trd.2023.103653},
author = {Yi He and Zhaocai Liu and Ziqi Song},
}

@article{GAIROLA2023103697,
title = {Optimization framework for integrated battery electric bus planning and charging scheduling},
journal = {Transportation Research Part D: Transport and Environment},
volume = {118},
pages = {103697},
year = {2023},
issn = {1361-9209},
doi = {10.1016/j.trd.2023.103697},
author = {Pranav Gairola and N. Nezamuddin},
}

@article{HE2023103587,
title = {Battery electricity bus charging schedule considering bus journey’s energy consumption estimation},
journal = {Transportation Research Part D: Transport and Environment},
volume = {115},
pages = {103587},
year = {2023},
issn = {1361-9209},
doi = {10.1016/j.trd.2022.103587},
author = {Jia He and Na Yan and Jian Zhang and Tao Wang and Yan-Yan Chen and Tie-Qiao Tang},
}

@article{BAO2023120512,
title = {An optimal charging scheduling model and algorithm for electric buses},
journal = {Applied Energy},
volume = {332},
pages = {120512},
year = {2023},
issn = {0306-2619},
doi = {10.1016/j.apenergy.2022.120512},
author = {Zhaoyao Bao and Jiapei Li and Xuehan Bai and Chi Xie and Zhibin Chen and Min Xu and Wen-Long Shang and Hailong Li},
}

@article{XIE2023103551,
title = {Collaborative optimization of electric bus line scheduling with multiple charging modes},
journal = {Transportation Research Part D: Transport and Environment},
volume = {114},
pages = {103551},
year = {2023},
issn = {1361-9209},
doi = {10.1016/j.trd.2022.103551},
author = {Dong-Fan Xie and Ya-Peng Yu and Guang-Jing Zhou and Xiao-Mei Zhao and Yong-Jun Chen},
}

@Article{Jovanovic202114206610,
AUTHOR = {Jovanovic, Raka and Bayram, Islam Safak and Bayhan, Sertac and Voß, Stefan},
TITLE = {A GRASP Approach for Solving Large-Scale Electric Bus Scheduling Problems},
JOURNAL = {Energies},
VOLUME = {14},
YEAR = {2021},
NUMBER = {20},
ARTICLE-NUMBER = {6610},
ISSN = {1996-1073},
DOI = {10.3390/en14206610}
}

@article{TANG2023103652,
title = {Optimization of single-line electric bus scheduling with skip-stop operation},
journal = {Transportation Research Part D: Transport and Environment},
volume = {117},
pages = {103652},
year = {2023},
issn = {1361-9209},
doi = {10.1016/j.trd.2023.103652},
author = {Chunyan Tang and Hudi Shi and Tao Liu},
}

@article{HE2023128227,
title = {Time-dependent electric bus and charging station deployment problem},
journal = {Energy},
volume = {282},
pages = {128227},
year = {2023},
issn = {0360-5442},
doi = {doi.org/10.1016/j.energy.2023.128227},
author = {Yi He and Zhaocai Liu and Yiming Zhang and Ziqi Song},
}

@article{perumal2022electric,
  title={Electric bus planning \& scheduling: A review of related problems and methodologies},
  author={Perumal, Shyam SG and Lusby, Richard M and Larsen, Jesper},
  journal={European Journal of Operational Research},
  volume={301},
  number={2},
  pages={395--413},
  year={2022},
  publisher={Elsevier},
  doi={10.1016/j.ejor.2021.10.058}
}

@article{PERUMAL2021105268,
title = {Solution approaches for integrated vehicle and crew scheduling with electric buses},
journal = {Computers \& Operations Research},
volume = {132},
pages = {105268},
year = {2021},
issn = {0305-0548},
doi = {10.1016/j.cor.2021.105268},
author = {Shyam S.G. Perumal and Twan Dollevoet and Dennis Huisman and Richard M. Lusby and Jesper Larsen and Morten Riis},
}

@article{LIU2023120483,
title = {Data-driven simulation-based planning for electric airport shuttle systems: A real-world case study},
journal = {Applied Energy},
volume = {332},
pages = {120483},
year = {2023},
issn = {0306-2619},
doi = {10.1016/j.apenergy.2022.120483},
author = {Zhaocai Liu and Qichao Wang and Devon Sigler and Andrew Kotz and Kenneth J. Kelly and Monte Lunacek and Caleb Phillips and Venu Garikapati},
}

@article{NAJAFI2025104664,
title = {Integrated optimization of charging infrastructure, electric bus scheduling and energy systems},
journal = {Transportation Research Part D: Transport and Environment},
volume = {141},
pages = {104664},
year = {2025},
issn = {1361-9209},
doi = {10.1016/j.trd.2025.104664},
author = {Arsalan Najafi and Kun Gao and Omkar Parishwad and Georgios Tsaousoglou and Sheng Jin and Wen Yi},
}

@article{bazarnovi2024problem,
  title={Problem of Locating and Allocating Charging Equipment for Battery Electric Buses under Stochastic Charging Demand},
  author={Bazarnovi, Sadjad and Cokyasar, Taner and Verbas, Omer and Mohammadian, Abolfazl Kouros},
  journal={European Journal of Operational Research},
  year={2025},
doi={10.1016/j.ejor.2025.07.064}
}

@misc{authority2022charging,
  title={{C}harging {F}orward: {C}{T}{A} bus electrification planning report},
  author={{CTA}},
  howpublished={Available at \url{https://www.transitchicago.com/assets/1/6/Charging\_Forward\_Report\_2-10-22\_(FINAL).pdf}},
note={Accessed on Mar. 4, 2025.},
  year={2022}
}

@misc{eiadieselprice,
  title={{Weekly Retail Gasoline and Diesel Prices}},
  author={{U.S. EIA}},
  howpublished={Available at \url{https://www.eia.gov/dnav/pet/pet_pri_gnd_dcus_nus_w.htm}},
note={Accessed on Feb. 20, 2025.},
  year={2025}
}

@misc{blsinflation,
  title={{Consumer Price Index, Chicago-Naperville-Elgin area – January 2025}},
  author={{U.S. BLS}},
  howpublished={Available at \url{https://www.bls.gov/regions/midwest/news-release/consumerpriceindex_chicago.htm}},
note={Accessed on Feb. 20, 2025.},
  year={2025}
}

@article{DAVATGARI2024953,
title = {Electric vehicle supply equipment location and capacity allocation for fixed-route networks},
journal = {European Journal of Operational Research},
volume = {317},
number = {3},
pages = {953-966},
year = {2024},
issn = {0377-2217},
doi = {10.1016/j.ejor.2024.04.022},
author = {Amir Davatgari and Taner Cokyasar and Anirudh Subramanyam and Jeffrey Larson and Abolfazl (Kouros) Mohammadian},
}

@misc{sundin2018scheduling,
  title={Scheduling of Electric Buses with Column Generation},
  author={Sundin, Daniel},
  year={2018},
  url={https://www.diva-portal.org/smash/get/diva2%3A1283727/FULLTEXT01.pdf}
}

@article{lin2016column,
  title={A column generation algorithm for the bus driver scheduling problem},
  author={Lin, Dung-Ying and Hsu, Ching-Lan},
  journal={Journal of Advanced Transportation},
  volume={50},
  number={8},
  pages={1598--1615},
  year={2016},
  publisher={Wiley Online Library},
  doi = {doi.org/10.1002/atr.1417}
}

@article{XU2025111138,
title = {Enhancing column generation by reinforcement learning-based hyper-heuristic for vehicle routing and scheduling problems},
journal = {Computers \& Industrial Engineering},
volume = {206},
pages = {111138},
year = {2025},
issn = {0360-8352},
doi = {10.1016/j.cie.2025.111138},
author = {Kuan Xu and Li Shen and Lindong Liu},
}

@article{DUAN2023104175,
title = {Integrated optimization of electric bus scheduling and charging planning incorporating flexible charging and timetable shifting strategies},
journal = {Transportation Research Part C: Emerging Technologies},
volume = {152},
pages = {104175},
year = {2023},
issn = {0968-090X},
doi = {10.1016/j.trc.2023.104175},
author = {Mengyuan Duan and Feixiong Liao and Geqi Qi and Wei Guan}
}

@article{GERBAUX2025106848,
title = {A machine-learning-based column generation heuristic for electric bus scheduling},
journal = {Computers \& Operations Research},
volume = {173},
pages = {106848},
year = {2025},
issn = {0305-0548},
doi = {10.1016/j.cor.2024.106848},
author = {Juliette Gerbaux and Guy Desaulniers and Quentin Cappart}
}

@article{gkiotsalitis2023exact,
  title={An exact approach for the multi-depot electric bus scheduling problem with time windows},
  author={Gkiotsalitis, Konstantinos and Iliopoulou, Christina and Kepaptsoglou, K},
  journal={European Journal of Operational Research},
  volume={306},
  number={1},
  pages={189--206},
  year={2023},
  publisher={Elsevier},
  doi = {10.1016/j.ejor.2022.07.017}
}

@article{desaulniers1998multi,
  title={Multi-depot vehicle scheduling problems with time windows and waiting costs},
  author={Desaulniers, Guy and Lavigne, June and Soumis, Francois},
  journal={European Journal of Operational Research},
  volume={111},
  number={3},
  pages={479--494},
  year={1998},
  publisher={Elsevier},
  doi = {10.1016/S0377-2217(97)00363-9}
}

@article{dell1993heuristic,
  title={Heuristic algorithms for the multiple depot vehicle scheduling problem},
  author={Dell'Amico, Mauro and Fischetti, Matteo and Toth, Paolo},
  journal={Management Science},
  volume={39},
  number={1},
  pages={115--125},
  year={1993},
  publisher={INFORMS},
  doi={10.1287/mnsc.39.1.115}
}

@article{freling2001models,
  title={Models and algorithms for single-depot vehicle scheduling},
  author={Freling, Richard and Wagelmans, Albert PM and Paix{\~a}o, Jos{\'e} M Pinto},
  journal={Transportation Science},
  volume={35},
  number={2},
  pages={165--180},
  year={2001},
  publisher={INFORMS},
  doi={10.1287/trsc.35.2.165.10135}
}

@article{rinaldi2020mixed,
  title={Mixed-fleet single-terminal bus scheduling problem: Modelling, solution scheme and potential applications},
  author={Rinaldi, Marco and Picarelli, Erika and D'Ariano, Andrea and Viti, Francesco},
  journal={Omega},
  volume={96},
  pages={102070},
  year={2020},
  publisher={Elsevier},
  doi={10.1016/j.omega.2019.05.006}
}

@inproceedings{cokyasar2023solving,
  title={Solving the electric vehicle scheduling problem at large-scale},
  author={Cokyasar, Taner and Verbas, Omer and Davatgari, Amir and Mohammadian, Abolfazl Kouros},
  booktitle={2023 IEEE 26th International Conference on Intelligent Transportation Systems (ITSC)},
  pages={1134--1139},
  year={2023},
  organization={IEEE},
  doi={10.1109/ITSC57777.2023.10422441}
}

@inproceedings{cokyasar2023electric,
  title={Electric vehicle scheduling problem with tour combinations},
  author={Cokyasar, Taner and Verbas, Omer and Auld, Joshua},
  booktitle={Procedia Computer Science, The 14th International Conference on Ambient Systems, Networks and Technologies Networks (ANT) and The 6th International Conference on Emerging Data and Industry 4.0 (EDI40)},
  volume={220},
  pages={413--420},
  year={2023},
  publisher={Elsevier},
  doi={10.1016/j.procs.2023.03.053}
}

\newpage

\renewcommand{\appendixname}{Supplementary Material}
\appendix

\pretocmd{\section}{\FloatBarrier}{}{}

\setcounter{table}{0}
\renewcommand{\thetable}{A\arabic{table}}

\setcounter{figure}{0}
\renewcommand{\thefigure}{A\arabic{figure}}

\section{Charger network layout}\label{app:charger_network_layout}

\begin{figure}[!ht]
    \centering
    \includegraphics[width=1\linewidth]{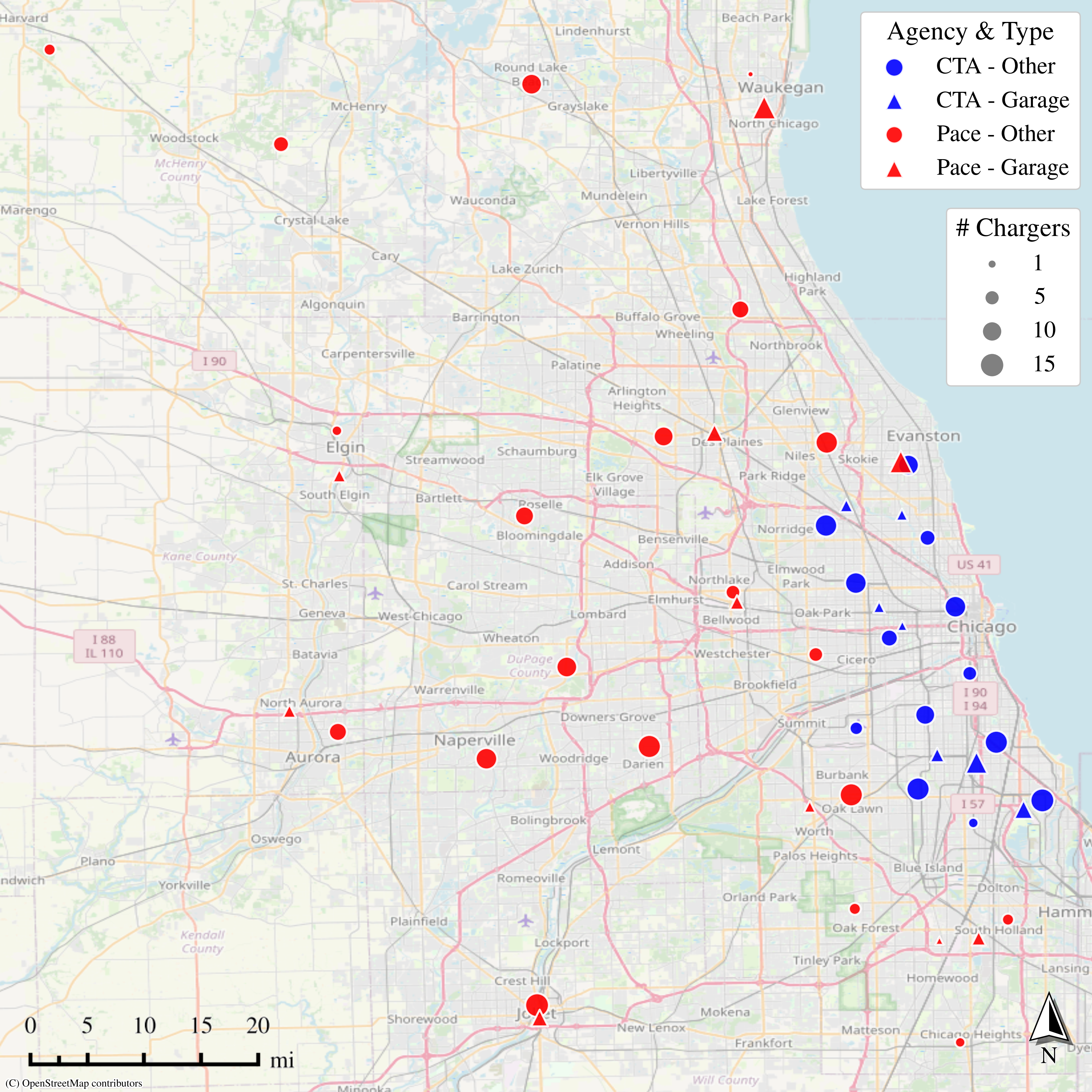}
    \caption{Charger network layout displaying the count of chargers at garage and non-garage locations.}
    \label[fig]{charger_layout}
\end{figure}

\section{More computational analysis}\label{app:more_computational}
To evaluate the stability of the CG heuristic, we analyzed the deviation in solution time across the 10 replicas. \cref{comp_time_dev}, using a \textit{violin} plot, visualizes the distribution of time deviations from the mean for each instance. For small instances ($|\mathcal{I}| \leq 25$, \cref{comp_time_dev_small}), the deviation is negligible (typically $\pm 4$ seconds), indicating deterministic behavior on tractable problems. A transition phase is observed for medium-to-large instances ($|\mathcal{I}| \in \{100, 250\}$), shown in \cref{comp_time_dev_large}. Here, the variance is highest, visualized by the elongated violin shapes. This suggests that the stochastic nature of the CG search tree leads to varying convergence paths before the time limit is reached. However, for the largest instances ($|\mathcal{I}| \in \{500, 1000\}$), the deviation drops to zero. This confirms that for massive problem sizes, the solver consistently hits the computational time limit across all replicas, utilizing the full budget in a uniform manner.

\begin{figure}[!ht]
    \centering
    \begin{subfigure}[t]{0.49\textwidth}
        \centering
        \includegraphics[width=\linewidth]{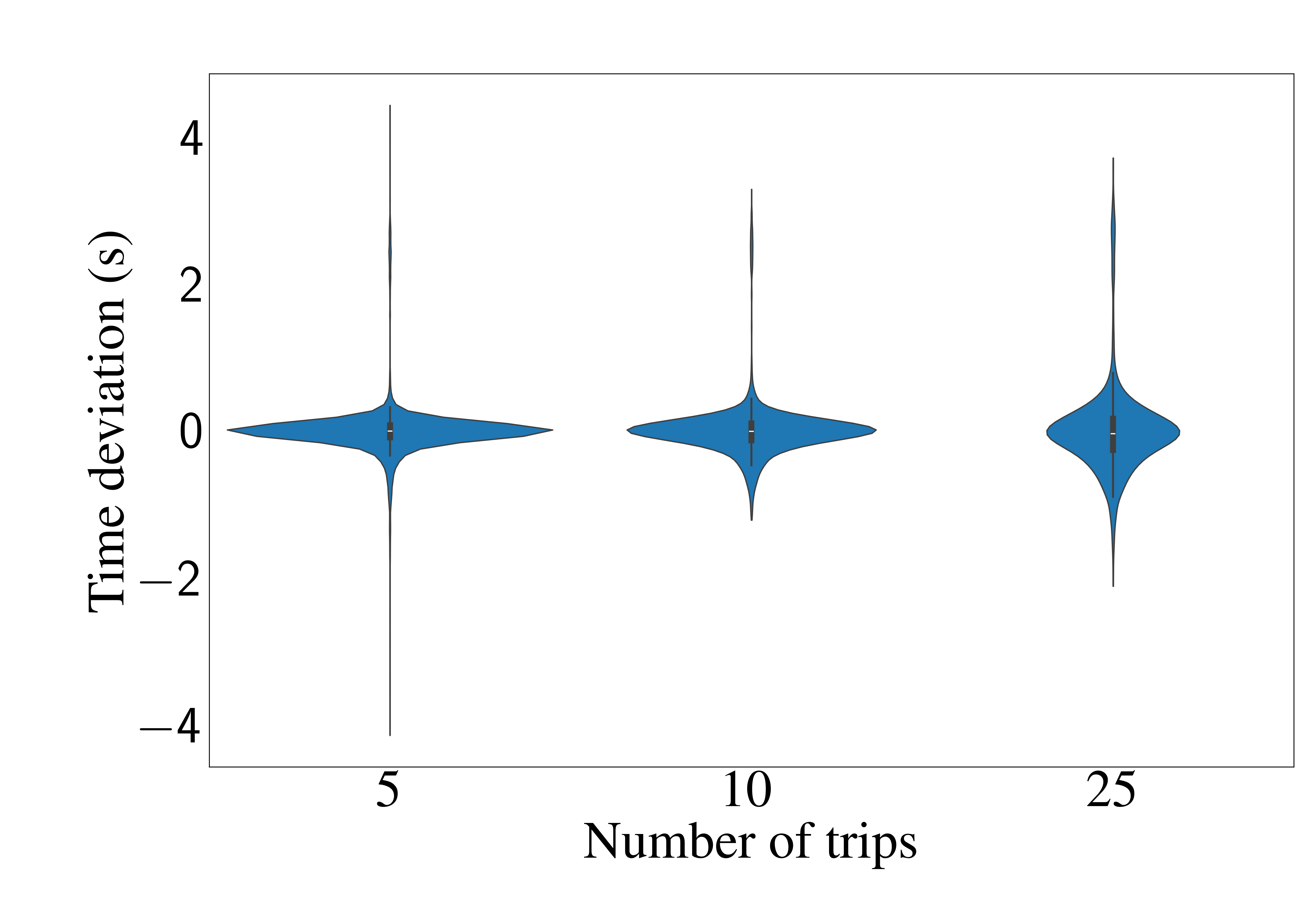} 
        \caption{Scenarios with $|\mathcal{I}| \in \{5, 10, 25\}$}
        \label[fig]{comp_time_dev_small}
    \end{subfigure}
    \hfill
    \begin{subfigure}[t]{0.49\textwidth}
        \centering
        \includegraphics[width=\linewidth]{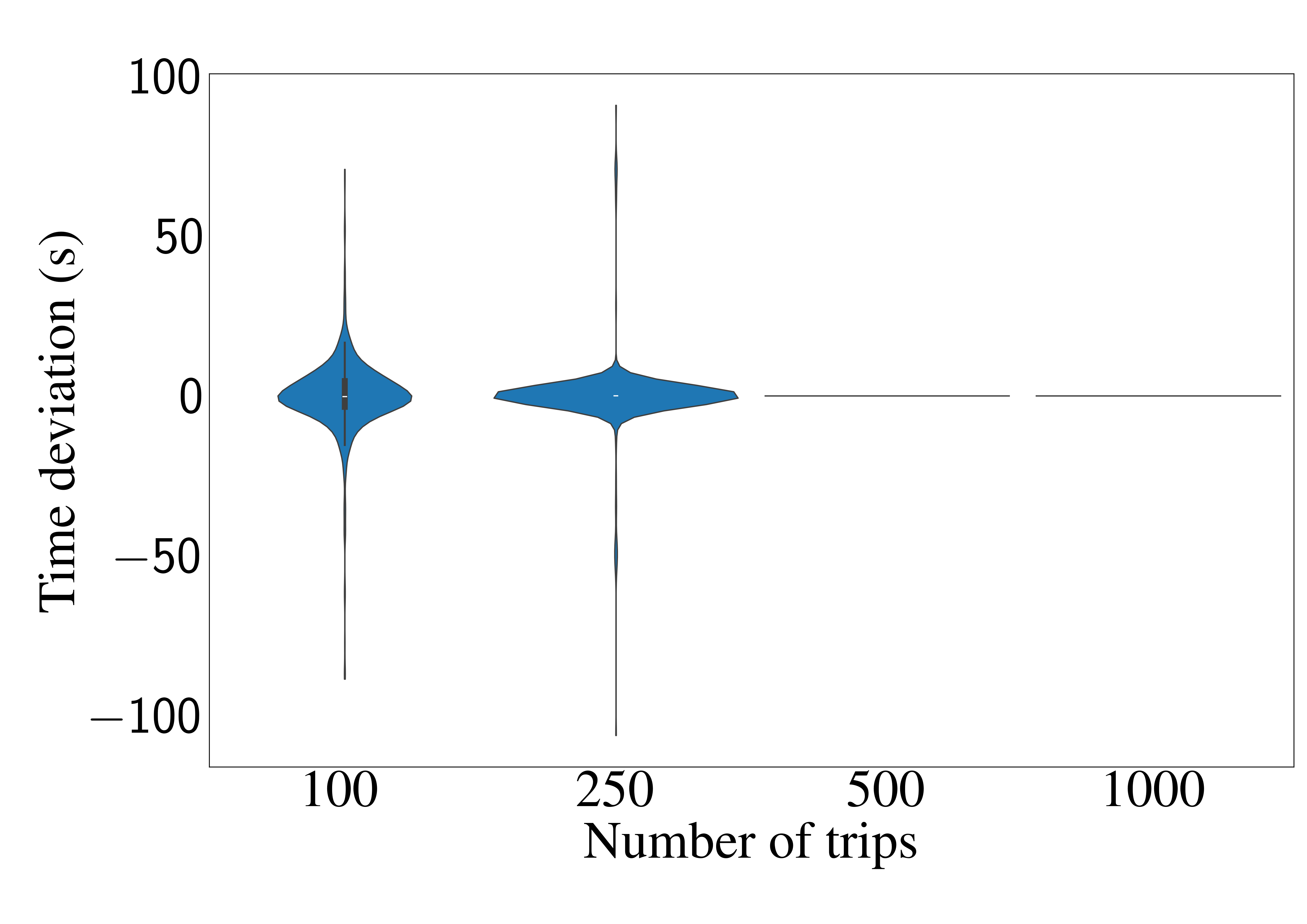} 
        \caption{Scenarios with $|\mathcal{I}| \in \{100, 250, 500, 1000\}$}
        \label[fig]{comp_time_dev_large}
    \end{subfigure}
    \caption{Distribution of computational time deviations across CG replicas.}
    \label[fig]{comp_time_dev}
\end{figure}

In the exact method formulation, the parameter $T^\Delta$ governs the discrete resolution used to track charger availability. A smaller $T^\Delta$ offers higher temporal precision but increases the size of the time-expanded graph, potentially impacting computational tractability. To quantify this trade-off, we conducted a sensitivity analysis on 840 problem instances.

The experimental design spanned six problem sizes, $|\mathcal{I}| \in \{5, 10, 25, 50, 75, 100\}$, and seven discretization levels, $T^\Delta \in \{60, 150, 300, 600, 1800, 3600, 7200\}$ seconds. For each combination of $|\mathcal{I}|$ and $T^\Delta$, 20 random instances were solved using the exact method. Time limits were adjusted based on problem size: 10, 30, 120, 300, 480, and 600 seconds, respectively.

\cref{comp_T_delta} illustrates the solution time distributions. We observe a strictly monotonic relationship: increasing $T^\Delta$ significantly reduces computational effort across all problem sizes. This effect is most pronounced in larger instances. For $|\mathcal{I}| \geq 75$ (\cref{comp_T_delta_75} and \cref{comp_T_delta_100}), a fine resolution of $T^\Delta=60$ consistently forces the solver to hit the time limit (600s). However, increasing the step size to $T^\Delta \geq 1800$ reduces the solution time to under 100 seconds, effectively rendering intractable problems solvable.

\begin{figure}[!ht]
    \centering
    \begin{subfigure}[t]{0.49\textwidth}
        \centering
        \includegraphics[width=\linewidth]{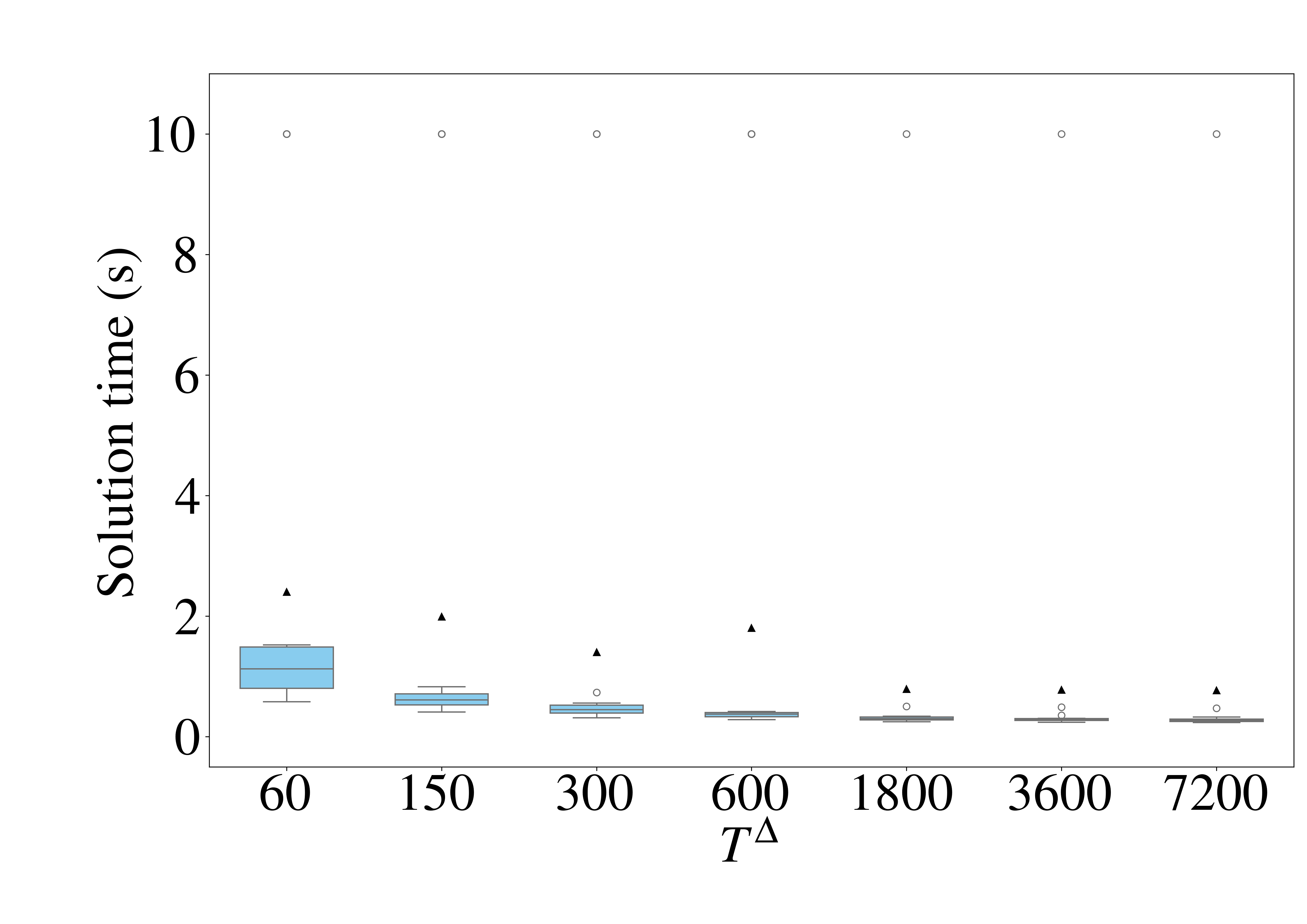} 
        \caption{Scenarios with $|\mathcal{I}|=5$}
        \label[fig]{comp_T_delta_5}
    \end{subfigure}
    \hfill
    \begin{subfigure}[t]{0.49\textwidth}
        \centering
        \includegraphics[width=\linewidth]{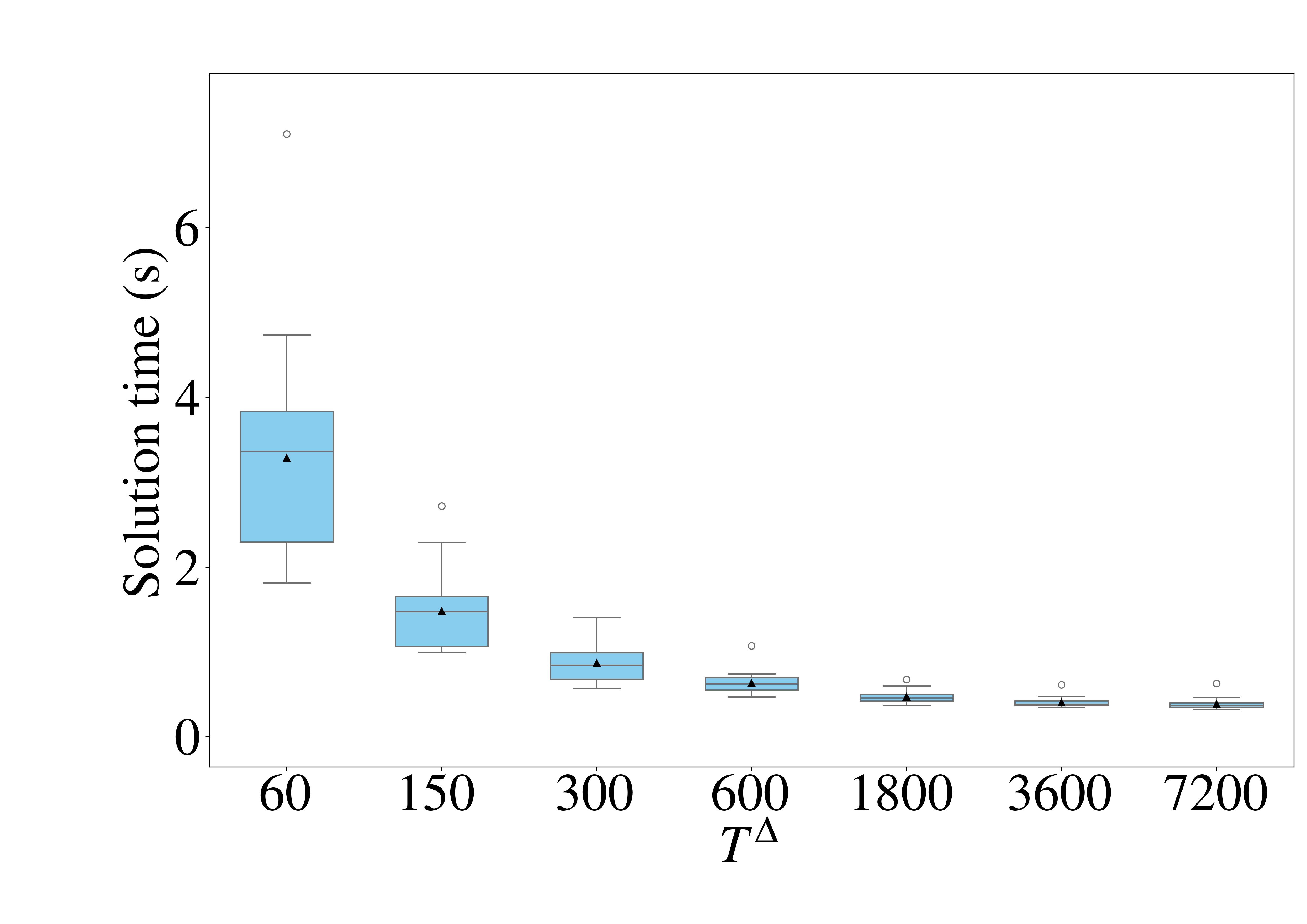} 
        \caption{Scenarios with $|\mathcal{I}|=10$}
        \label[fig]{comp_T_delta_10}
    \end{subfigure}
    
    \begin{subfigure}[t]{0.49\textwidth}
        \centering
        \includegraphics[width=\linewidth]{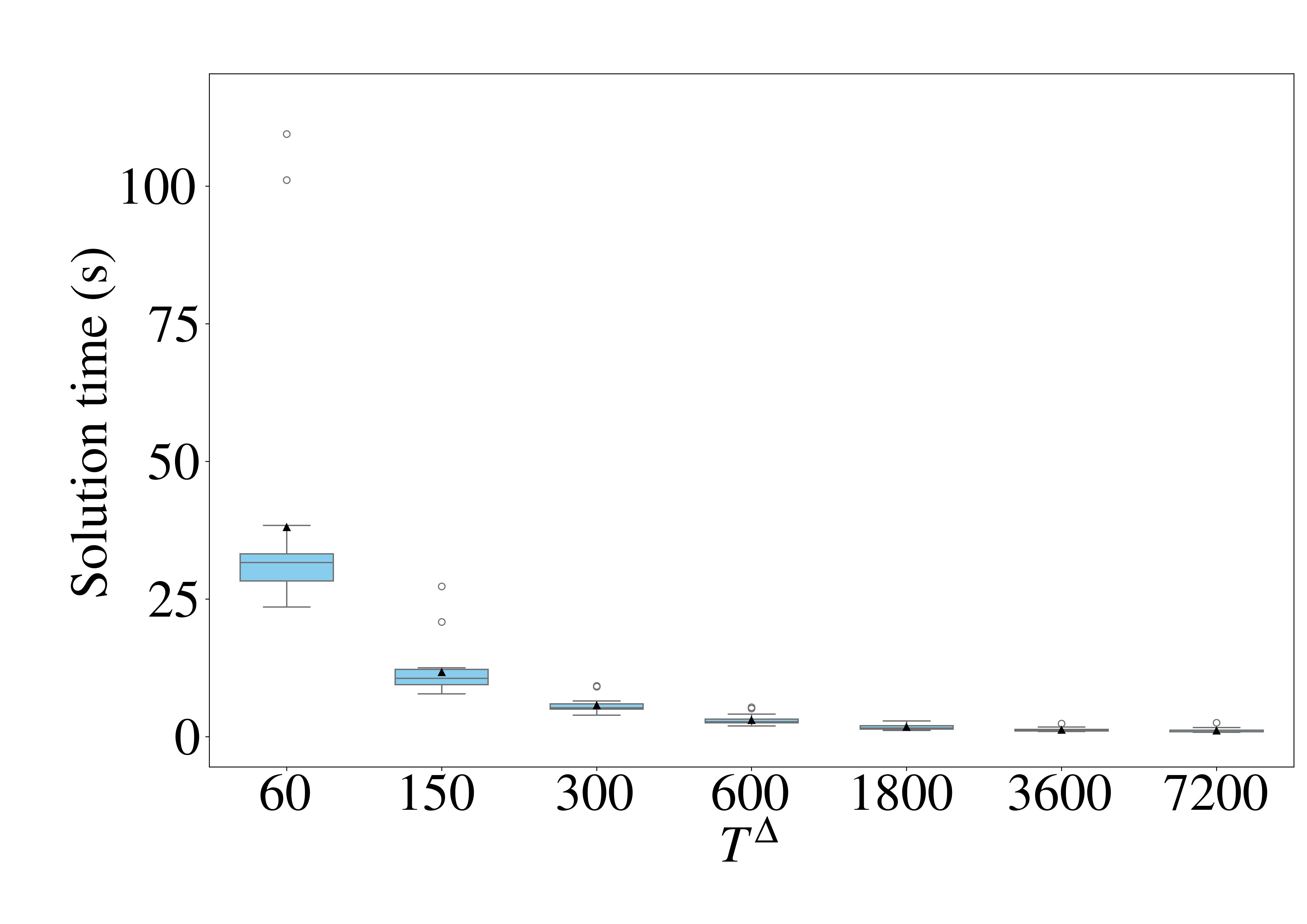} 
        \caption{Scenarios with $|\mathcal{I}|=25$}
        \label[fig]{comp_T_delta_25}
    \end{subfigure}
    \hfill
    \begin{subfigure}[t]{0.49\textwidth}
        \centering
        \includegraphics[width=\linewidth]{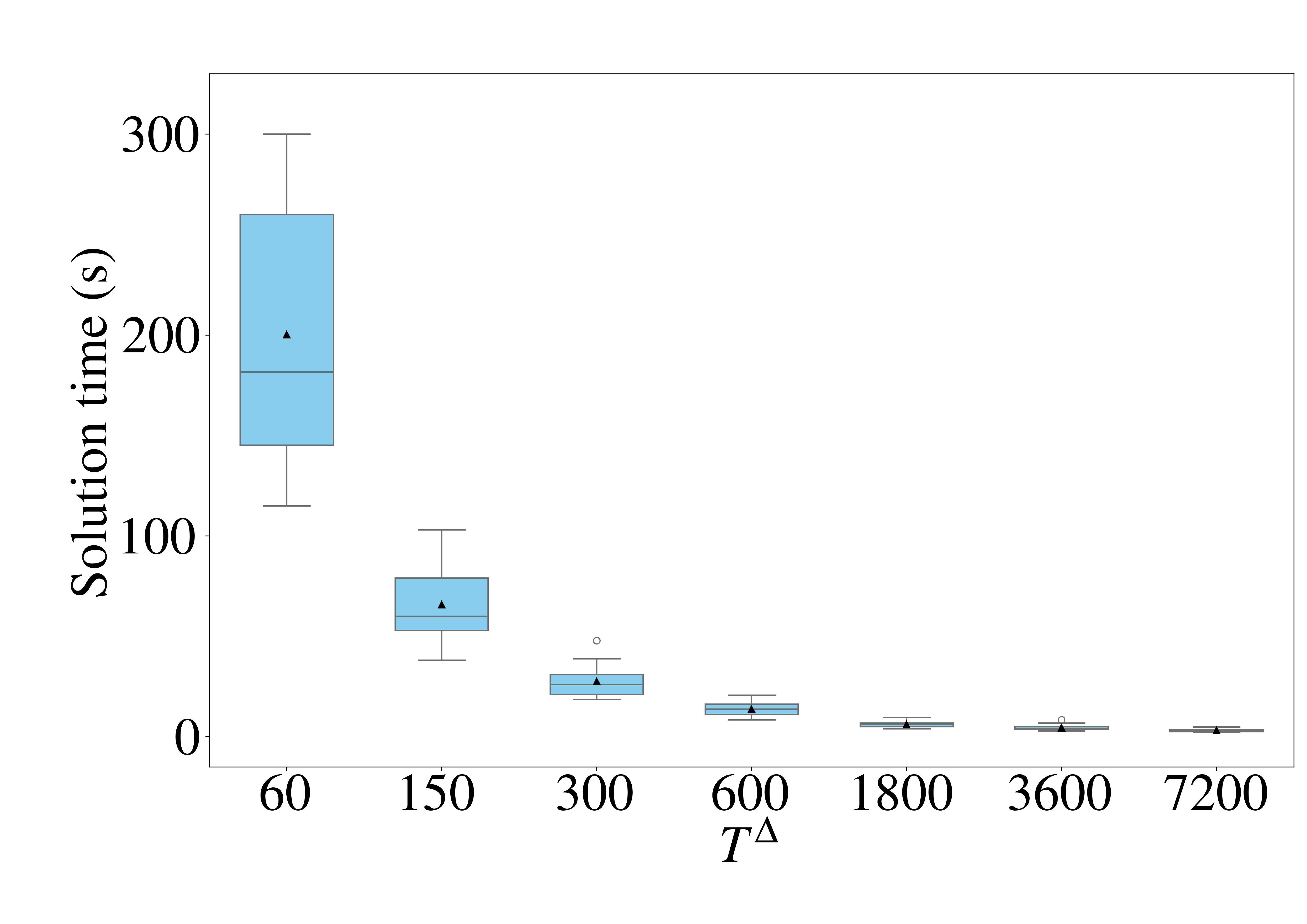} 
        \caption{Scenarios with $|\mathcal{I}|=50$}
        \label[fig]{comp_T_delta_50}
    \end{subfigure}
    
    \begin{subfigure}[t]{0.49\textwidth}
        \centering
        \includegraphics[width=\linewidth]{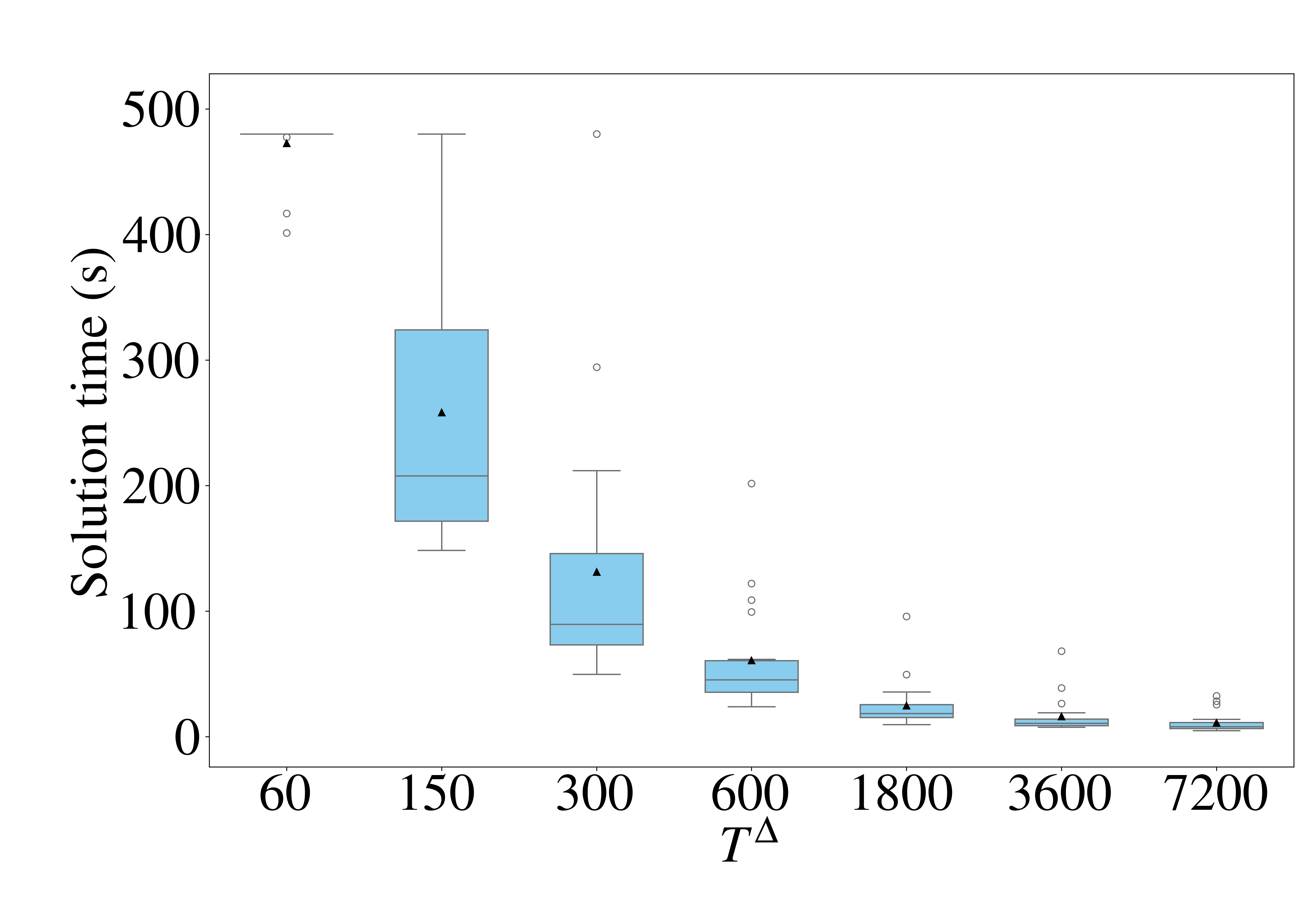} 
        \caption{Scenarios with $|\mathcal{I}|=75$}
        \label[fig]{comp_T_delta_75}
    \end{subfigure}
    \hfill
    \begin{subfigure}[t]{0.49\textwidth}
        \centering
        \includegraphics[width=\linewidth]{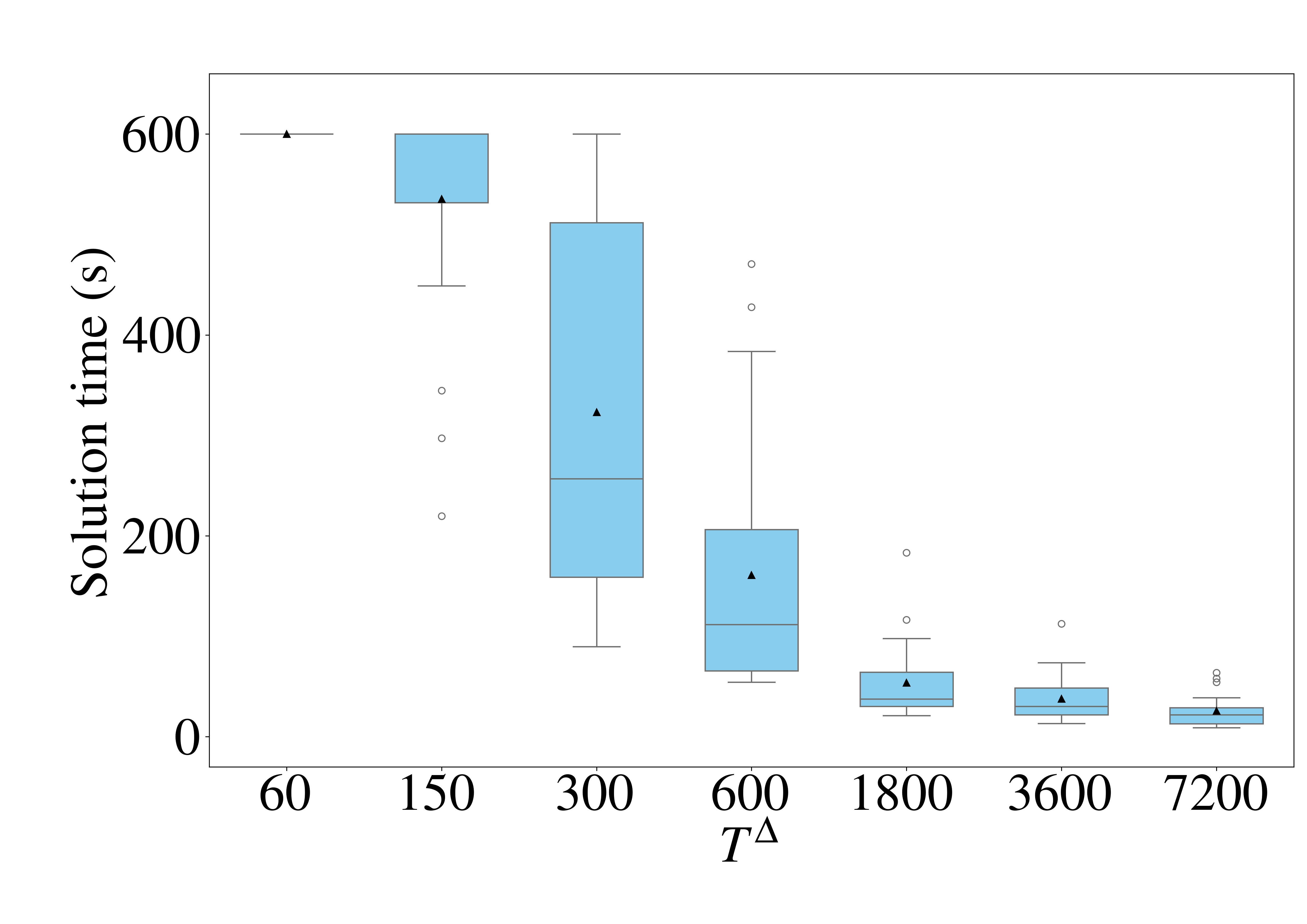} 
        \caption{Scenarios with $|\mathcal{I}|=100$}
        \label[fig]{comp_T_delta_100}
    \end{subfigure}
    \caption{Impact of time discretization parameter $T^\Delta$ on solution time across varying problem sizes.}
    \label[fig]{comp_T_delta}
\end{figure}

Furthermore, we analyzed the impact of discretization on solution quality. Comparing the optimal objective function values across different $T^\Delta$ settings for the same instances, we observed no significant deviation. We attribute this stability to the nature of the experimental design, where the number of chargers at designated locations is relatively high compared to the number of trips ($|\mathcal{I}|$). Consequently, competition for charging slots was minimal, preventing the granularity of time tracking from becoming a binding constraint.

It is important to acknowledge that this observation may not generalize to scenarios with highly constrained infrastructure. In cases with fewer chargers or locations, a coarser discretization (large $T^\Delta$) could artificially prolong charger occupancy status, causing the model to view available chargers as occupied and leading to suboptimal schedules. However, for the capacity levels examined in this study, larger $T^\Delta$ values provided substantial computational speedups, often by orders of magnitude, without practically compromising the quality of the charging schedule.

\section{More sensitivity analysis}\label[app]{app:sensitivity_plots}

This section provides the decomposed sensitivity analysis results. While the main text discusses the aggregate impact on total system cost, the figures below isolate the effects on vehicle capital costs, operational expenses, and fleet sizing requirements.

\cref{sensitivity_vehicle_cost} illustrates the sensitivity of total vehicle capital costs. A notable trend is observed in the \textit{Battery Capacity} scenario: increasing capacity from the baseline (Lever 0) to approximately 1.4x (Lever 2) yields a reduction in vehicle costs, as buses with extended range can complete longer blocks without needing replacement vehicles. Beyond this point, the benefit saturates. Conversely, the \textit{Charger Layout} scenario highlights a steep penalty for Levers 6 and 7 (garage only charging). Restricting charging to garages forces the fleet size to expand to maintain service frequency while buses are out of service for charging, directly inflating capital costs. Undoubtedly, reduction in vehicle cost parameters considerably reduce the vehicle capital costs. The parameters of charge power, electricity cost and diesel do not impact the vehicle cost component of the objective function in scenarios $|\sI|=25$. In larger problem sizes, these parametric choices have a slight impact.

\cref{sensitivity_operational_cost} presents the impact on daily operational expenses. The \textit{Charger Layout} scenario again reveals a critical insight: garage only settings result in a massive spike in operational costs (exceeding 4\% in some instances). This is driven by the significant non-revenue deadheading required to return to garages for mid-day charging, confirming the operational efficiency of distributed terminal charging. Additionally, the \textit{Electricity Cost} and \textit{Diesel Cost ($A^\nu=0$)} scenarios display the expected linear relationship with operational expenditures.

Finally, \cref{sensitivity_num_beb} depicts the percentage change in the required number of BEBs. The results closely mirror the vehicle cost findings. Tighter constraints in \textit{Battery Range} (moving from Lever 1 to Lever 8) generally force an increase in fleet size. Similarly, the restrictive garage only layout levers necessitate a larger fleet (an increase of roughly 2--6\%) to compensate for the inefficiency of garage-centric charging operations.

\begin{figure}[!ht]
    \centering
    \begin{subfigure}[t]{.9\textwidth}
        \centering
        \includegraphics[width=\linewidth]{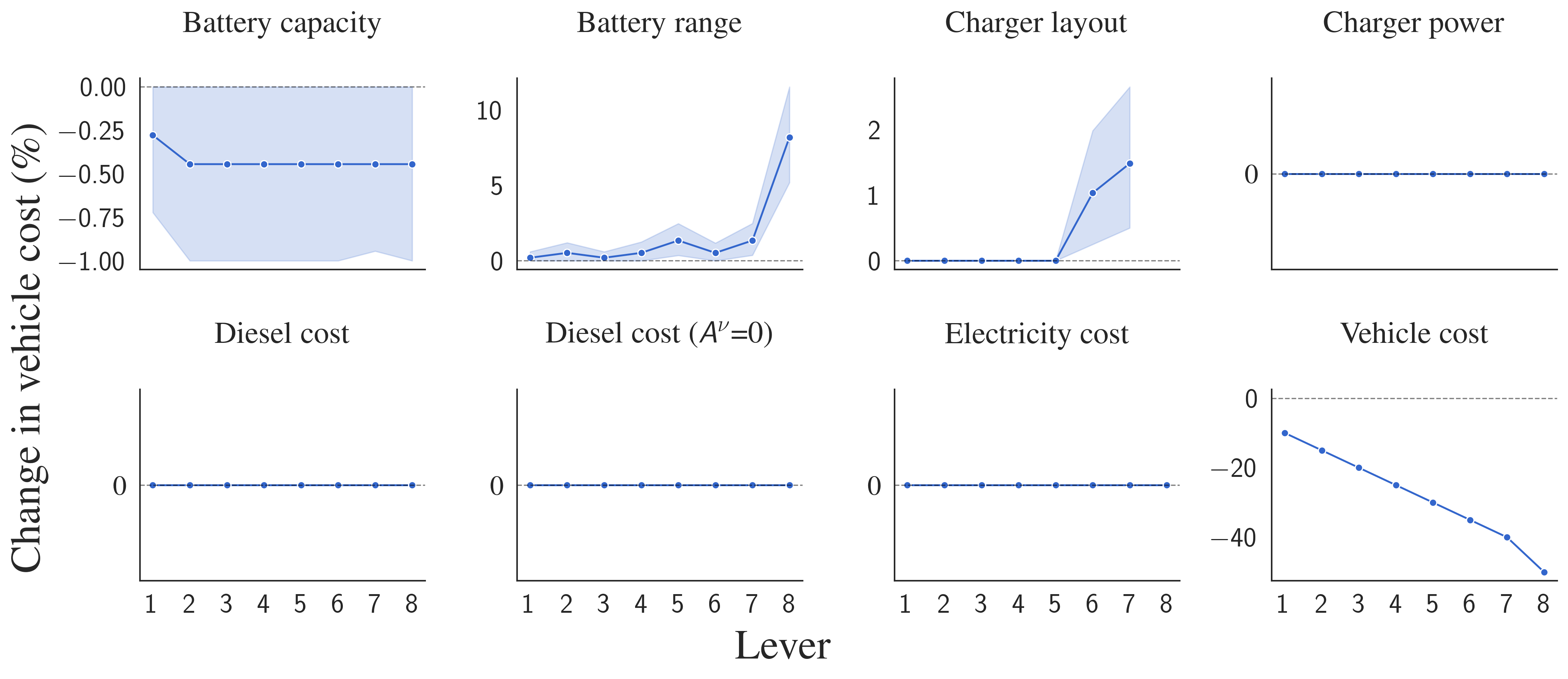} 
        \caption{$|\sI|=25$}
    \end{subfigure}
    \begin{subfigure}[t]{.9\textwidth}
        \centering
        \includegraphics[width=\linewidth]{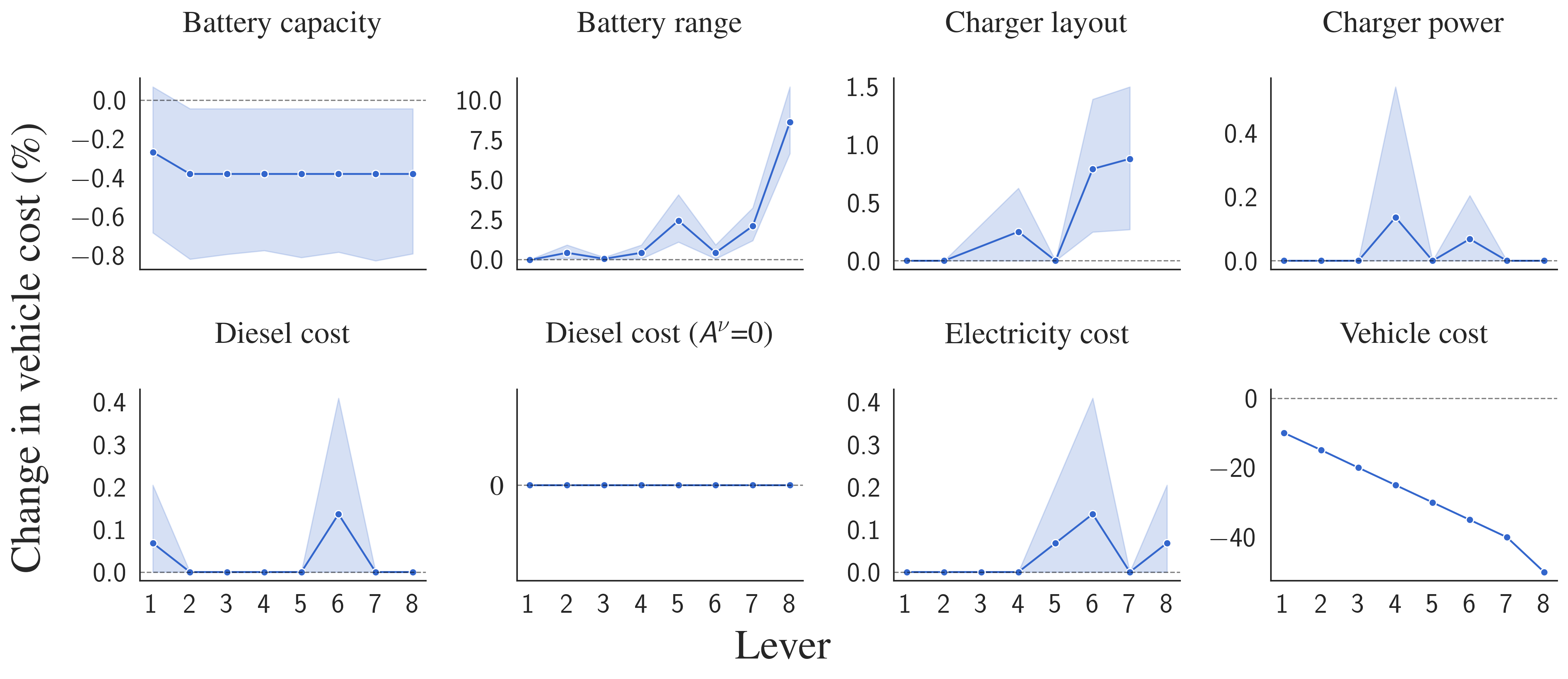} 
        \caption{$|\sI|=50$}
    \end{subfigure}
    \begin{subfigure}[t]{.9\textwidth}
        \centering
        \includegraphics[width=\linewidth]{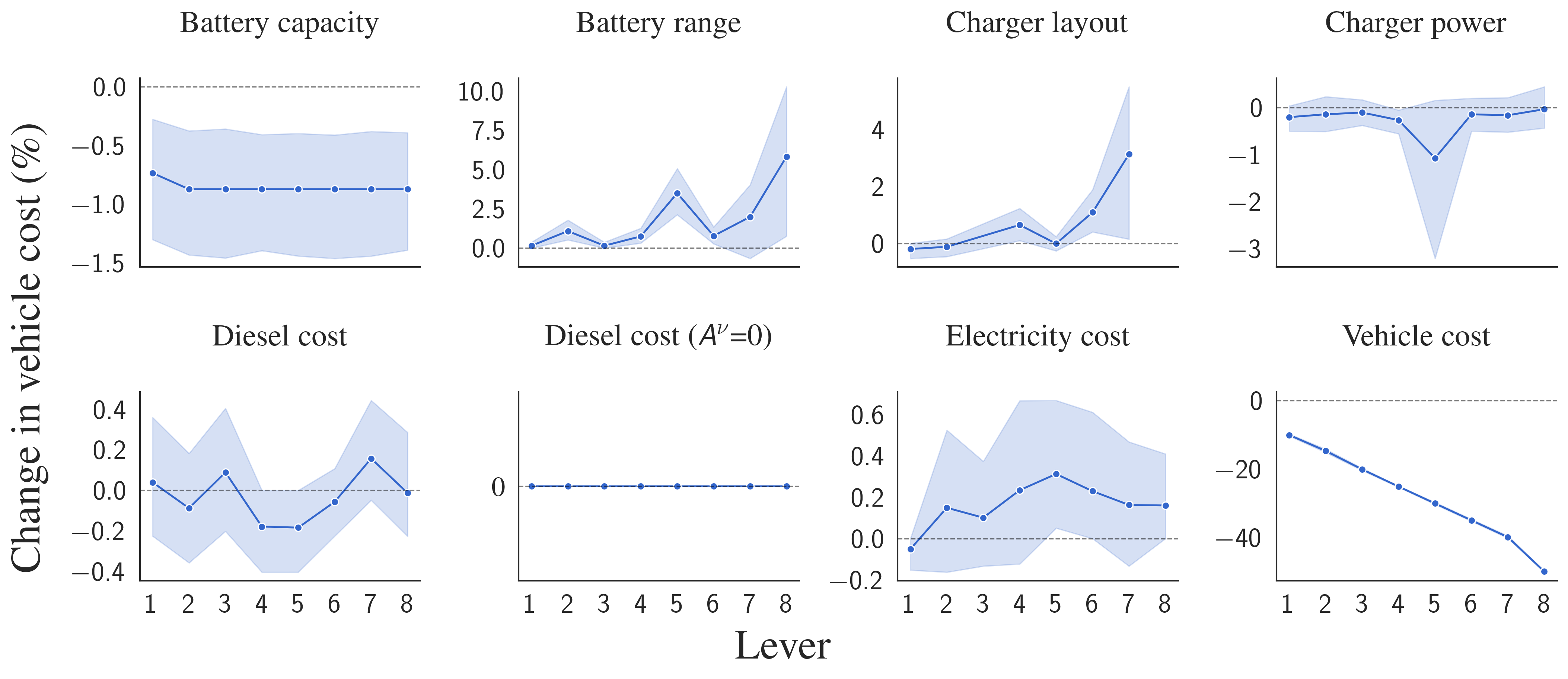} 
        \caption{$|\sI|=100$}
    \end{subfigure}
    \caption{Percentage change in total vehicle cost of scenarios with respect to their baseline levers.}
    \label[fig]{sensitivity_vehicle_cost}
\end{figure}

\begin{figure}[!ht]
    \centering
    \begin{subfigure}[t]{.9\textwidth}
        \centering
        \includegraphics[width=\linewidth]{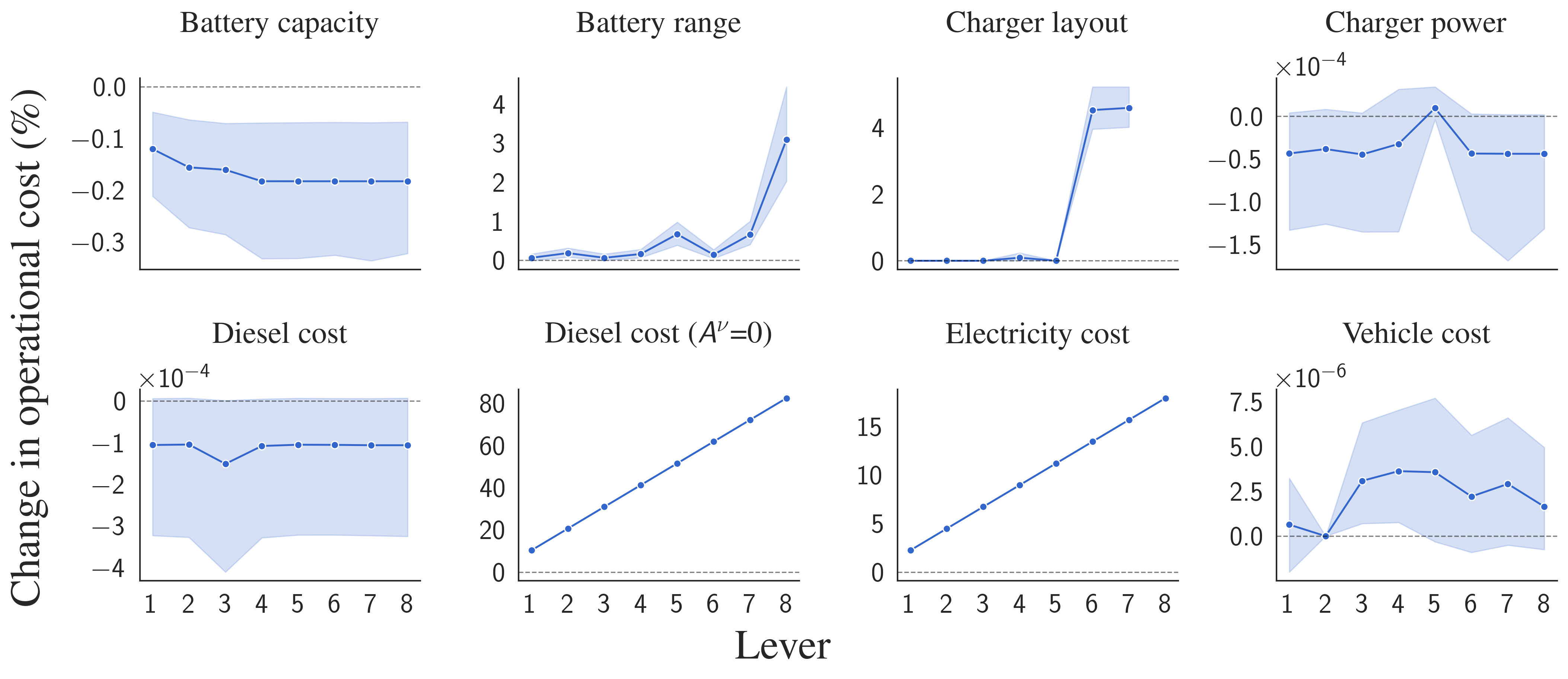} 
        \caption{$|\sI|=25$}
    \end{subfigure}
    \begin{subfigure}[t]{.9\textwidth}
        \centering
        \includegraphics[width=\linewidth]{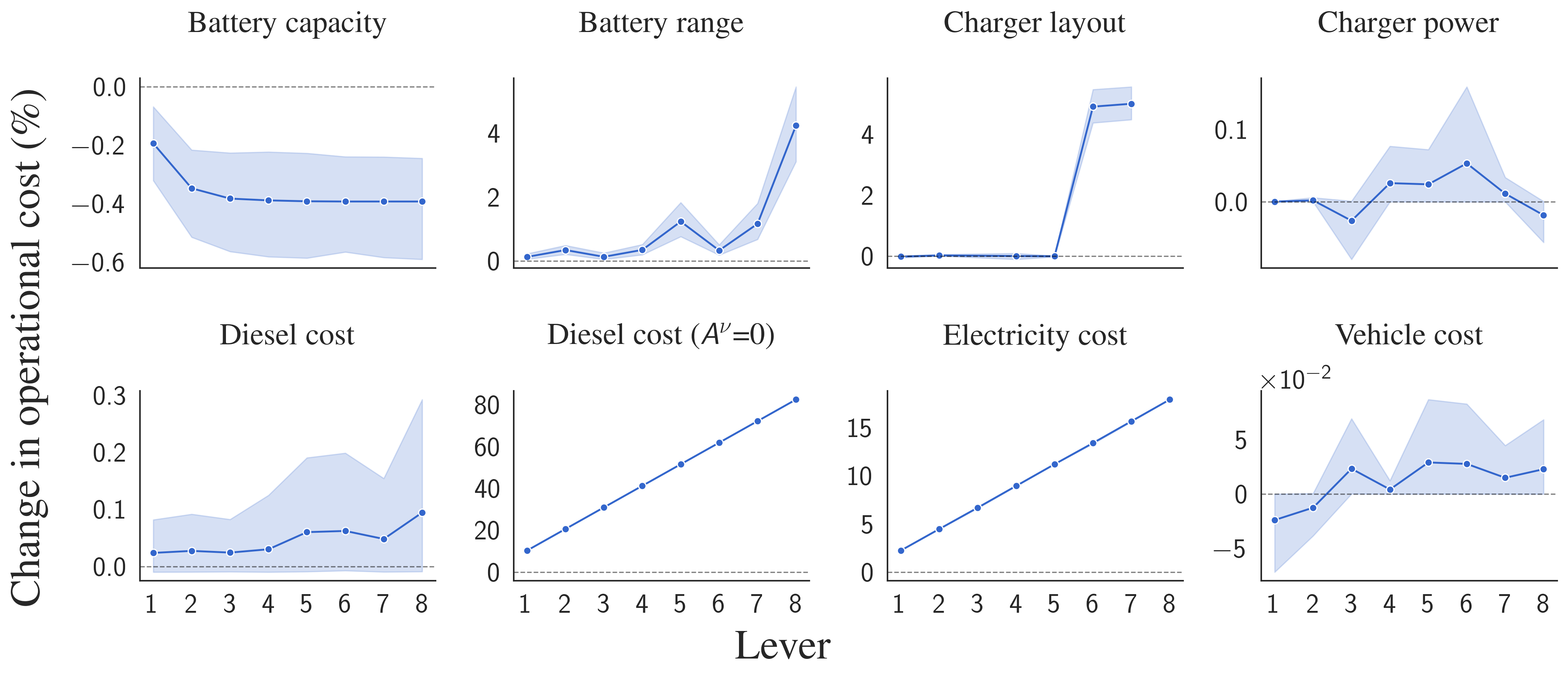} 
        \caption{$|\sI|=50$}
    \end{subfigure}
    \begin{subfigure}[t]{.9\textwidth}
        \centering
        \includegraphics[width=\linewidth]{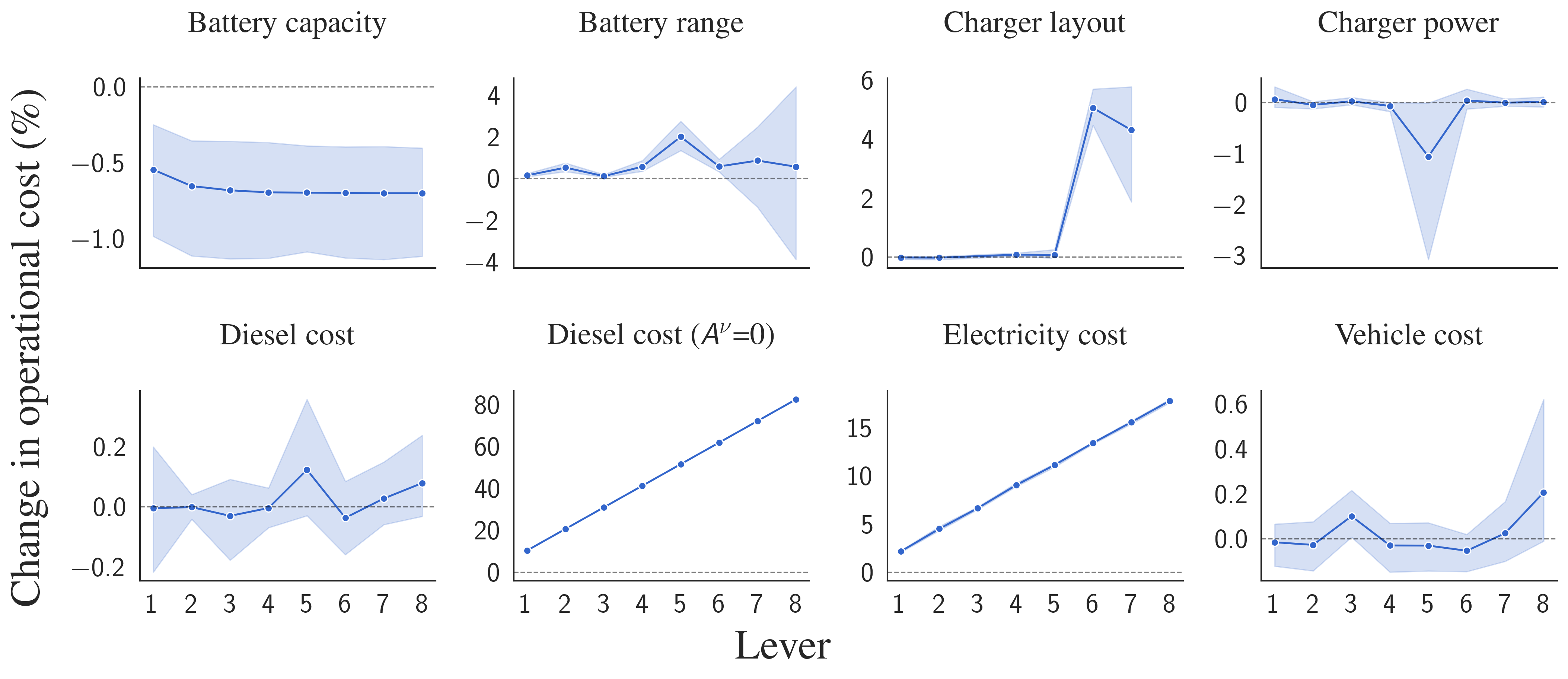} 
        \caption{$|\sI|=100$}
    \end{subfigure}
    \caption{Percentage change in operational cost of scenarios with respect to their baseline levers.}
    \label[fig]{sensitivity_operational_cost}
\end{figure}

\begin{figure}[!ht]
    \centering
    \begin{subfigure}[t]{.9\textwidth}
        \centering
        \includegraphics[width=\linewidth]{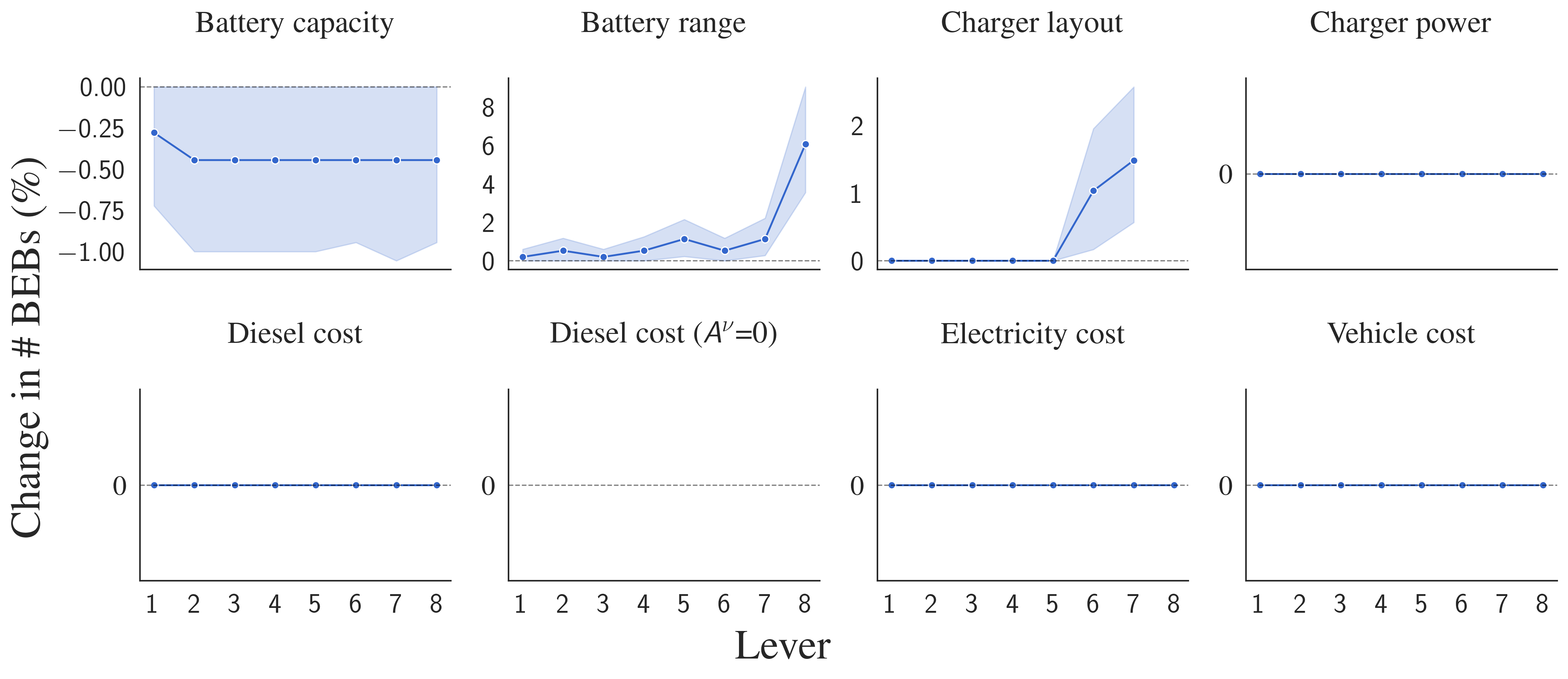} 
        \caption{$|\sI|=25$}
    \end{subfigure}
    \begin{subfigure}[t]{.9\textwidth}
        \centering
        \includegraphics[width=\linewidth]{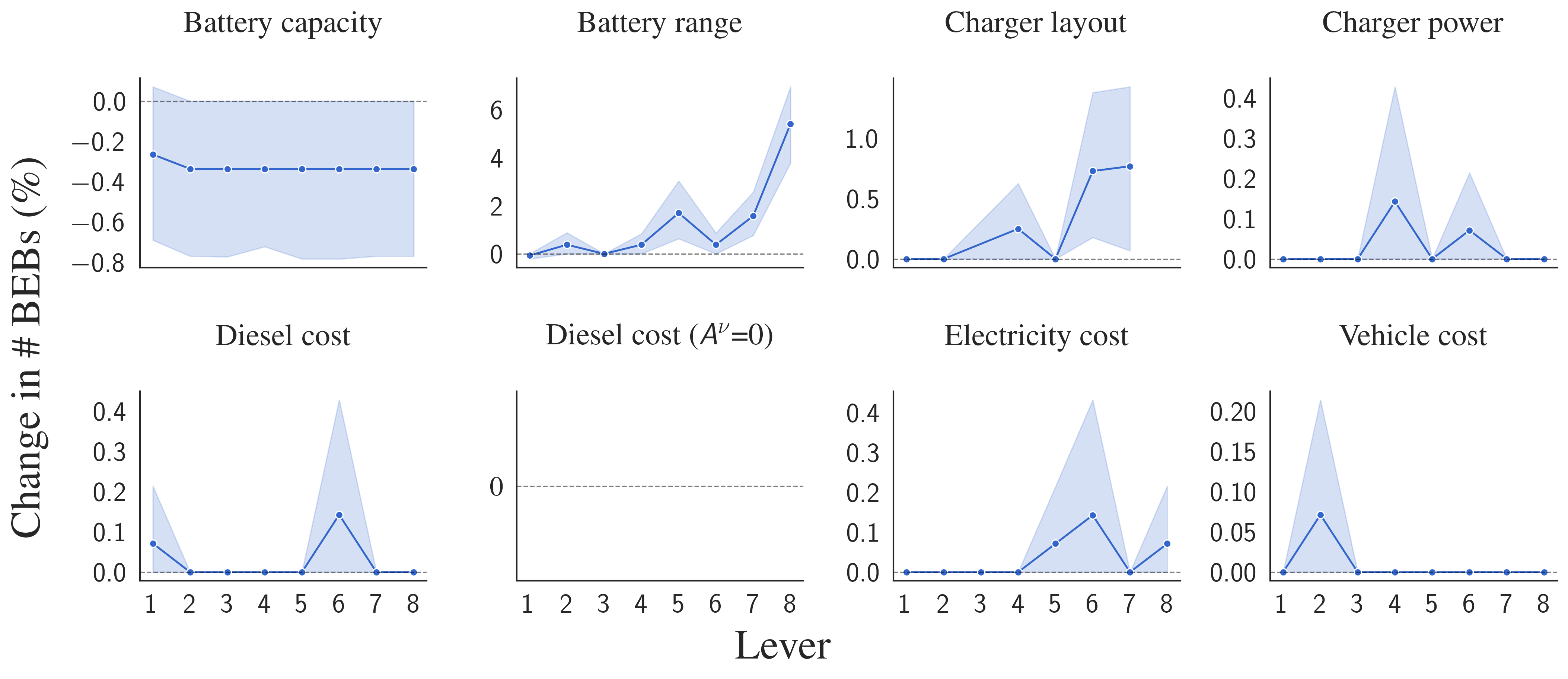} 
        \caption{$|\sI|=50$}
    \end{subfigure}
    \begin{subfigure}[t]{.9\textwidth}
        \centering
        \includegraphics[width=\linewidth]{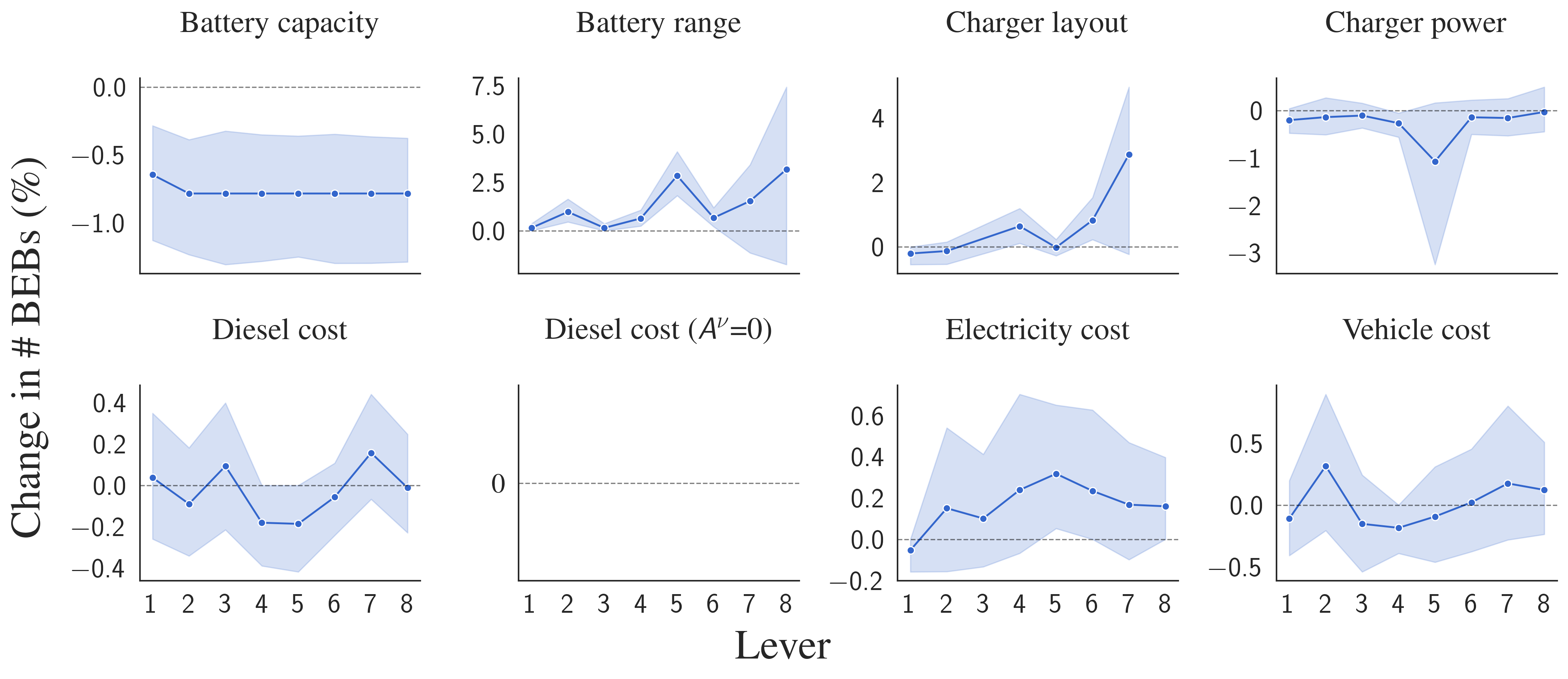} 
        \caption{$|\sI|=100$}
    \end{subfigure}
    \caption{Percentage change in number of BEBs of scenarios with respect to their baseline levers.}
    \label[fig]{sensitivity_num_beb}
\end{figure}

\vfill
\framebox{\parbox{.90\linewidth}{\scriptsize The submitted manuscript has been created by
        UChicago Argonne, LLC, Operator of Argonne National Laboratory (``Argonne'').
        Argonne, a U.S.\ Department of Energy Office of Science laboratory, is operated
        under Contract No.\ DE-AC02-06CH11357.
        The U.S.\ Government retains for itself,
        and others acting on its behalf, a paid-up nonexclusive, irrevocable worldwide
        license in said article to reproduce, prepare derivative works, distribute
        copies to the public, and perform publicly and display publicly, by or on
        behalf of the Government.
        The Department of Energy will provide public access
        to these results of federally sponsored research in accordance with the DOE
        Public Access Plan \url{http://energy.gov/downloads/doe-public-access-plan}.}}
\end{document}